\input amstex
\input xy
\xyoption{all}
\documentstyle{amsppt}
\document
\magnification=1200
\NoBlackBoxes
\nologo
\hoffset1.5cm
\voffset2cm
\vsize15.5cm
\def\F{\bold{F}}

\def\C{\bold{C}}
\def\P{\bold{P}}
\def\Q{\bold{Q}}

\def\Z{\bold{Z}}
\def\R{\bold{R}}
\def\F{\bold{F}}
\def\N{\bold{N}}

\def\bA{\bold{A}}
\def\bL{\bold{L}}
\def\bY{\bold{Y}}
\def\bG{\bold{G}}
\def\bS{\bold{S}}
\def\bW{\bold{W}}

\def\cA{\Cal{A}}

\def\cC{\Cal{C}}
\def\cO{\Cal{O}}

\def\cM{\Cal{M}}
\def\cV{\Cal{V}}
\def\cS{\Cal{S}}
\def\cZ{\Cal{Z}}
\def\cF{\Cal{F}}
\def\cL{\Cal{L}}
\def\cQ{\Cal{Q}}
\def\cX{\Cal{X}}
\def\cP{\Cal{P}}
\def\cU{\Cal{U}}
\def\cY{\Cal{Y}}
\def\cB{\Cal{B}}
\def\cD{\Cal{D}}
\def\cH{\Cal{H}}

\bigskip



\bigskip

\centerline{\bf HOMOTOPY SPECTRA}

\medskip

\centerline{\bf AND  DIOPHANTINE EQUATIONS}

\medskip

\centerline{\bf  Yuri~I.~Manin,\quad  Matilde~Marcolli}

\bigskip

{\it \  To Xenia and Paolo, from Yuri and Matilde, with all our love and gratitude.}

\bigskip

{\it ABSTRACT.} Arguably, the first bridge between vast, ancient,
 but disjoint domains of mathematical knowledge, --
topology and number theory, -- was built only during the last fifty years. This
bridge is {\it the theory of spectra in the stable homotopy theory.}
\smallskip

In particular, it connects $\Z$, the initial object in the theory of commutative rings,
with the sphere spectrum $\bold{S}$: see [Sc01] for one of versions of it.
This connection poses the challenge: discover new information in 
number theory using the independently developed machinery of homotopy
theory. (Notice that a passage in the reverse direction  has already generated
results about computability in homotopy theory: see [FMa20]
and references therein.)

\smallskip

In this combined research/survey paper we suggest to apply homotopy spectra to the problem of
distribution of rational points upon algebraic manifolds.
\smallskip

\medskip

Subjects: {\bf Algebraic Geometry (math.AG); Number Theory (math.NT); Topology (math.AT)}

Comments: 87 pages.

MSC-classes: 16E35, 11G50, 14G40, 55P43, 16E20, 18F30.

\medskip

{\it CONTENTS}

\medskip

0. Introduction and summary
\smallskip

1. Homotopy spectra: a brief presentation

\smallskip

2. Diophantine equations: distribution of rational points
on algebraic varieties

\smallskip

3. Rational points, sieves, and assemblers

\smallskip

4.  Anticanonical heights and points count

\smallskip

5. Sieves ``beyond heights'' ?

\smallskip

6. Obstructions and sieves
\smallskip

7. Assemblers and spectra for Grothendieck rings with exponentials

References

\bigskip

\centerline{{\bf 0. Introduction and summary}} 
\bigskip

{\bf 0.0. A brief history.} For a long stretch of time in the history of mathematics, 
 number theory and topology formed vast, but disjoint domains
of mathematical knowledge. Emmanulel Peyre reminds us in [Pe19] that  
the Babylonian clay tablet Plimpton 322 
(about 1800 BC) contained a list
of integer solutions of the ``Diophantine'' equation $a^2+b^2=c^2$:
archetypal theme of number theory, named after Diophantus of Alexandria 
(about 250 BC). 

\smallskip

``Topology'' was born much later, but arguably, its cousin -- modern 
measure theory -- goes back  to Archimedes, author of {\it Psammites} (``Sand Reckoner''),
who was approximately a contemporary of Diophantus. 
\smallskip
In modern language, Archimedes measures the volume of observable universe
by counting the number of small grains of sand necessary to fill this volume.
Of course, many {\it qualitative} geometric models and {\it quantitative} estimates 
of the relevant distances known in a narrow small world of scientists of his times
 precede his calculations. Moreover, since the estimated
numbers of grains of sands are quite large (about $10^{64}$), Archimedes
had to invent and describe a system of notation for large numbers
going far outside the possibilities of any of the standard ancient systems.

\smallskip
The construction of the first bridge between number theory and topology  was accomplished 
only about fifty years ago (around 1970): it is
it is 
\smallskip
\centerline{\it the theory of spectra in stable homotopy theory}
\bigskip

In this paper, after a brief survey of  relevant  setups in homotopy theory 
and number theory,  we focus upon the ongoing research dedicated to
the applications of spectra  to the problems of distribution
of rational/algebraic points on algebraic varieties.

\bigskip

{\bf 0.1. Integers and finite sets.} Below, we will briefly describe the contents of our paper
appealing to the intuition of a reader who is already accustomed to {\it categorical}
reasoning.

\smallskip
 
Intuitively, the set of non--negative integers $\bold{N}$ can be imagined as
embodiment of {\it ``sizes''} (formally, cardinalities) of all finite sets  (including the empty set). 
\smallskip
Already to compare sizes of two disjoint sets, that is
 to interpret the
inequality $m\le n$ in $\bold{N}$, one needs to look at the {\it category $FSets$}:
if $m$ is the size of $M$ and $n$ is the size of $N$, then $m\le n$ means that
in $FSets$ there is a {\it monomorphism} $M\to N$, or equivalently, an {\it epimorphism} $N\to M$.

\smallskip

A categorical interpretation of multiplication of integers uses {\it direct products of finite sets}.
Associativity of multiplication is a reflection of the {\it monoidal structure} on $FSet$.

\smallskip

Addition in $\bold{N}$ already requires more sophisticated constructions in $FSet$:
if two sets $X, Y$ are disjoint, then the cardinality of $X\cup Y$ is the sum of
cardinalities  of $X$ and $Y$, but far from all pairs $(X,Y)$ are disjoint. So an additional
formalism must be developed.

\medskip
{\bf 0.2. Integers and topological spaces.}  Much less popular is the background
vision of a natural number $n\in \bold{N}$ as the {\it dimension of a manifold,}
say the sphere $S^n$. 

\smallskip
Here already the definition of {\it equality} of numbers requires quite
sophisticated technicalities, if we want to lift it: it involves functorial definitions of classes of morphisms
that we will declare {\it invertible ones} after passing to the category of categories.
It is at this step that the idea of {\it homotopy}, homotopically equivalent maps,
and homotopy classes of spaces and maps, acquires the key meaning.

\smallskip

While dimension is not a homotopy invariant, 
the natural numbers have another topological manifestation, 
in terms of certain categories, which we review below, 
describing the combinatorics of simplices. 

\smallskip
But as soon as we accomplish this, addition and multiplication in $\bold{N}$
become liftable to an appropriate category and, after some more work, 
passage from $\bold{N}$ to $\bold{Z}$ becomes liftable as well.

\medskip

{\bf 0.3. From finite sets to topological spaces, and back.} 
The critically important bridging of finite sets and  topological spaces
is furnished by the machinery of {\it simplicial sets}, in particular 
simplicial sets associated to {\it coverings}.

\smallskip

Roughly speaking, this passage
replaces a finite set $X$ by the simplex $\sigma(X)$ of which $X$ is the set of vertices,
and an admissible set theoretical map $X\to Y$ by its extension by convex combinations from $\sigma(X)$
to $\sigma(Y)$. 
Admissibility condition can be described in terms of the simplex category recalled below.

\smallskip

In the reverse direction, the passage from a {\it topological space with its covering}
(say, by a family of open subsets) $\{U_a\,| a\in A\}$ to a simplicial
set replaces each $U_a$ by a vertex of a simplex, and associates
simplices $\Delta_J$ with non--empty intersections $\cap_{i\in J} U_{a_i}$.
\smallskip
For a detailed treatment, see [GeMa03]. Here we will only recall
the combinatorial setup (see [GeMa03] Section~II.4(a) p.~58)
with the following properties.

\smallskip

The simplex category $\Delta$ has objects the totally ordered sets $[n] := \{0,1,\dots ,n\}$
and morphisms are nondecreasing maps $f: [n] \to [m]$. 

\smallskip

The set of all morphisms is generated by two classes of maps: ``$i$--th face'' $\partial^i_n$
and ``$i$--th degeneration'' $\sigma^i_n$. The $i$--th face is the increasing injection $[n-1] \to [n]$
not taking the value $i$, and $\sigma^i_n$ is the nondecreasing surjection taking the value $i$ twice.

\smallskip

All the relations between faces and degenerations are generated by the relations
$$
\partial^j_{n+1}   \partial^i_{n} = \partial^i_{n+1}  \partial^{j-1}_{n}  \quad \roman{for} \quad i < j ;
$$
$$
\sigma^j_{n-1}   \sigma^i_{n+1} =   \sigma^i_n  \sigma^{j+1}_{n+1}  \quad \roman{for}\quad i \le j ;
$$
$$
\sigma^j_{n-1}\partial^i_n =
\cases 
\partial^i_{n-1}\sigma^{j-1}_{n-2}  \quad \roman{for}\quad  i<j ;\\
id_{[n-1]} \quad \roman{for}\quad  i\in \{j,j+1\} ; \\
\partial^{i-1}_{n-1}\sigma^{j}_{n-2} \quad \roman{for}\quad  i>j+1 .
\endcases
$$
\smallskip

A {\it simplicial set} $X_{\bullet}$ (or simply $X$) can be then defined as a functor from the
simplex category $\Delta$ to some category $Set$ of sets. 

\smallskip

Thus, a simplicial set is the structure consisting of a family of sets $(X_n)$, $n=0,1,2, \dots$
and a family of maps $X(f): X_n\to X_m$, corresponding to each nondecreasing map
$f: [m] \to [n]$, such that $X(id_{[n]}) = id_{X_n}$, and $X(g\circ f) =  X(f)\circ X(g)$.

\smallskip 

Restricting ourselves by consideration of only $X(\partial^i_n)$ and $X(\sigma^i_n)$,
and the respective relations, we get the convenient and widely used description
of simplicial sets.

\smallskip

In order to pass from combinatorics to geometry, we must start with geometric simplices.

\smallskip
The {\it $n$--dimensional simplex} $\Delta_n$ is the topological space embedded into
 real affine $n+1$--dimensional space endowed with coordinate system:
$$
\Delta_n := \left\{(x_0,\dots ,x_n\} \in \bold{R}^{n+1}\cap [0,1]^{n+1} \, \biggr|  \, \sum_{i=0}^n x_i = 1\right\}.
$$


Starting with a simplicial set $(X_n)$, we construct its {\it geometric realization} $|X|$.
It is the set, endowed with a topology,
$$
|X| := \cup_{n=0}^{\infty} (\Delta_n\times X_n)/R ,
$$ 
where $R$ is the weakest equivalence relation identifying boundary subsimplices
in $\Delta_n\times X_n$ and  $\Delta_m\times X_m$ corresponding to nondecreasing
maps $f: [m] \to [n]$: see [GeMa03], pp. 6--7. In particular, the boundary of $\Delta_{n+1}$ for $n\ge 2$
is the geometric realisation of the union of $n+1$ boundary $(n-1)$--dimensional simplices:
this is a simplicial model of $S^{n-1}$. 
\smallskip

One popular path in the opposite direction produces a simplicial set
from a topological  space $\Cal{Y}$ endowed with a covering by open
(or closed) sets $(U_{\alpha})$, $\alpha \in A$.

\smallskip

Put for $n\ge 1$,
$$
X_n := \{ (\alpha_1, \dots ,\alpha_n)\,  | \, \cap_{i=1}^{n}U_{\alpha_i} \ne \emptyset\},
$$
$$
X(f)(\varphi ) := \varphi \circ \Delta_f, \quad \roman{where}\  f: [m] \to [n], \, \,\,\, \Delta_f:\Delta_m\to \Delta_n .
$$

The role of homotopy appearing at this point finds  a beautiful expression in the
following fact: {\it  the geometric realisation $|X|$ is homotopically equivalent to $\Cal{Y}$,
if $U$ is a locally finite covering, and all nonempty finite intersections 
$\cap_{i=1}^{n}U_{\alpha_i}$ are contractible.}

\smallskip

These constructions form the background for the much more sophisticated
machinery of {\it homotopical spectra} sketched below in Sec. 1. 

\smallskip

A note on terminology: in this paper we will adopt the somewhat unconventional 
terminology ``homotopical spectra" or ``homotopy spectra" 
to simply refer to ``spectra as meant in the context of homotopy theory" (as opposed to all
the other conflicting notions of ``spectrum" that exist in mathematics). While we understand
that this terminology is nonstandard, it is a convenient shorthand, especially when other notions
of spectrum may also be present.

\medskip

{\bf 0.4. Structure of the paper.}
Homotopy theory methods have found some recent significant applications in arithmetic geometry,
for instance in the work of Corwin--Schlank [CorSch20] on higher obstructions to the existence of rational points.
We argue in this paper that homotopy theoretic methods, especially centered on the notions
of assemblers and spectra, can also provide a novel viewpoint on various aspects of Manin's problem on 
the distribution of rational points of bounded height and some related ideas, including the case of the
motivic height zeta functions studied by Chambert-Loir and Loeser [ChamLoe15], and Bilu [Bilu18]. 
The main new results in this paper consist of the construction of various relevant sieves and assemblers 
(Lemma~3.5.1 and Propositions~3.6.2 and 4.4.1, as well as Propositions~7.2, 7.5, 7.6, 7.7, and Theorem~7.4), 
and the proof that the motivic Fourier transform and the motivic Poisson summation formula, stated
as relations in a Grothendieck ring of varieties with exponentials, lift to the categorical level in the form 
of covering families in an assembler (Theorems~7.12 and 7.14.1).  
We present various related questions concerning the categorification
of aspects of classical and motivic height zeta function. 

\medskip

{\bf 0.4.1 Plan of exposition.} 
Section 1 is dedicated to a brief,
but sufficiently formal presentation of homotopy-theoretic spectra.
In particular, we would like to draw the attention of a reader to the category of
$\Gamma$--spaces that has several different applications to
the construction and study of new geometries: for the context of geometries 
``in characteristic 1'', cf. [CoCons16],  [MaMar18], [LMM19], and references therein.
The main use of $\Gamma$-spaces in this paper will be in Sections~3.3 and 7.9.1.

\smallskip
Section 2 introduces an approach to distributions of rational/algebraic
points on algebraic varieties based upon counting of points
of bounded sizes. Here the central role is played by the notion
of {\it heights}: ways to measure sizes of rational/algebraic points
on algebraic varieties taking into account their positions in a
projective space to which the variety is/can be embedded.

\smallskip

The pure categorical Section 3 is dedicated to the basic machinery of {\it assemblers}
that we later use to bridge diophantine geometry and homotopical algebra, following
[Za17a], [Za17b], [Za17c].
Its central intuitive notion is similar to one that lies in the background
of the definition of functor $K_0$ of an abelian category:
we replace each object of a category by the formal sum of
its ``irreducible'' pieces, neglecting ways in which these pieces are
assembled together.

\smallskip

The new results of Section~3 are contained in Sections~3.5 and 3.6, where we
construct sieves (Lemma~3.5.1) and an assembler $\cC_U$ (Proposition~3.6.2)
with Grothendieck group $K_0(\cC_U)$ that
detects decompositions of $U$ into a closed subvariety $V$ and its open
complement $W=U\smallsetminus V$ where $W\subset U$ satisfies a strict
inequality $0< \sigma(W,L_W)< \sigma(U,L_U)$ of the abscissas of convergence
of the height zeta function with respect to a line bundle $L_U$.  Arithmetic
stratifications in the sense of [BaMa90] give examples of disjoint covering families 
in this assembler. 

\smallskip

In Section 4 we  introduce arithmetical/geometric environments
that are more restricted, but better describable by homotopy means.
In particular Section 4.2 recalls specific classes of varieties (such as the Fano complete intersections
discussed in Section~4.2.3 and the generalized flag manifolds discussed in Section~4.2.4) for
which asymptotic formulae are known.  Section~4 is based mostly on [BaMa90] and [FrMaTsch89].
The only new results are in Section~4.4 (Proposition~4.4.1).

\smallskip

In Section~4.4 we show that the sieves and assemblers constructed in Section~3.5 and 3.6,
based on conditions on the abscissa of convergence of the height zeta function can be
generalized to other similar constructions, involving asymptotic conditions on the
number of rational points of bounded height. We construct in Proposition~4.4.1 assemblers
with associated Grothendieck group that detects splittings of $U$ into a closed subvariety $Y$
and its open complement, where $Y$ is a strongly (or weakly) accumulating subvariety. 

\smallskip

 Section 5 considers possibilities to avoid point count
on algebraic varieties,
replacing it by a machinery of measurable sets in adelic spaces.
This allows us to include consideration of Brauer--Manin
obstructions and more intuitive notions of equidistribution
of rational points. We pose the question of a possible formulation
in the language of assemblers. Section~5 is mostly based on [Pe18] and [Sa20]. 

\smallskip

Section 6 is a survey of some very  recent (2020) results  of Corwin--Schlank  [CorSch20],
and we highlight the unexpected appearance of assemblers
in the setup of varieties with empty sets of rational points,
and obstructions for their existence. 

\smallskip

Finally, in Section 7 we obtain categorifications of some parts of the recent work  
by Chambert--Loir and Loeser  [ChamLoe15], dedicated  to lifting
 arithmetic height zeta--functions 
(main generating functions for point counting) to  motivic height zeta functions.
We present the problem of how to lift these motivic height zeta functions
to the categorical level of assemblers and homotopy-theoretic spectra.
We discuss various aspects of this question, providing as initial result in this
direction the construction of categorical lifts 
of the motivic Fourier transform and the motivic Poisson summation formula.  

\smallskip

The first part of Section~7, from Section~7.1 to Proposition~7.7, introduces and extends
various generalizations of the Grothendieck ring of varieties $K_0(\cV_K)$. There are
three generalizations to begin with that have been variously considered in previous literature:

\medskip

(1) {\it from varieties to stacks} (see [BeDh07], [Ek09]): the Grothendieck ring of 
stacks $K_0(\cS_K)$ is obtained as a localization $\cM_K$ of the Grothendieck ring of varieties $K_0(\cV_K)$,
obtained by inverting the classes $[Gl_n]$;

\medskip

(2) {\it from varieties to exponential sums} (see [CluLoe10], [ChamLoe15]): classes $[X]$
of varieties in $K_0(\cV_K)$ are replaced by classes $[X,f]$ of varieties with a morphism $f: X \to \bA^1$
in the Grothendieck ring with exponentials $KExp_K$,
with relations designed to reflect the properties of exponential sums $\sum_{x\in X(\F_q)}\chi(f(x))$, for
$\chi: \F_q\to \C^*$ a character.

\medskip

(3) {\it from Grothendieck rings to assemblers} (see [Za17a], [Za17b], [Za17c]): the Grothendieck
ring of varieties in $K_0(\cV_K)$ is realized as $\pi_0 K(\cC_{\cV_K})$ of the spectrum $K(\cC_{\cV_K})$
associated to an assember $\cC_{\cV_K}$.

\medskip

A common generalization of (1) and (2) is introduced in Section~7.3 where we obtain a
relation between Grothendieck rings summarized by the diagram
$$ \xymatrix{ & K_0(\cV_K)\ar[r] \ar[d] \ar[dl] & KExp_K \ar[d] \ar[dr] & \\
K_0(\cS_K) & \ar[l]_{\simeq} \cM_K \ar[r] & Exp\cM_K \ar[r]^{\simeq} &  K_0(Exp\cS_K) ,  } $$
where $Exp\cM_K$ is the localization of $KExp_K$ at the classes $[GL_n,0]$ and 
$K_0(Exp\cS_K)$ is the Grothendieck ring of stacks with exponentials.

\smallskip

A common generalization of (2) and (3) is obtained in Theorem~7.4, where it is shown
that there is a simplicial assembler $\cC^{KExp}_K=\cC/\Phi$, obtained as the cofiber of an endomorphism of
an assembler, with $\pi_0 K(\cC/\Phi)=KExp_K$. This generalization fits into a diagram 
$$ \xymatrix{ K_0(\cV_K) \ar[r] & KExp_K \\ K(\cC_{\cV_K}) \ar[u] \ar[r] & K(\cC/\Phi), } $$
where the bottom line is a map of connective spectra and the upward maps are maps from
the spectral to their $\pi_0$'s. 

\smallskip

The common simultaneous generalization of (1), (2), and (3) is then obtained in Proposition~7.5,
with the construction of a simplicial assembler $\cC^{KExp\cS}_K =\cC^\cS/\Phi$ with
$\pi_0 K(\cC^{KExp\cS}_K)=Exp\cM_K=K_0(Exp\cS_K)$, the Grothendieck ring of
stacks with exponentials.

\smallskip

Propositions~7.6 and 7.7 describe the relative version of the previous construction, where one starts
with the Grothendieck ring $K_0(\cV_S)$ of varieties over a base scheme $S$ and
considers the corresponding generalizations of the other notions. In particular, Proposition~7.6 
obtains the relative version of the generalization of (2) and (3) and Propositions~7.7 gives the
relative form of the generalization of (1), (2), and (3). 

\smallskip

The main reason for developing these common generalizations of the Grothendieck
rings of varieties, stacks, and varieties with exponentials, with their associated assemblers,
is the fact that the Grothendieck rings $KExp_S$ and $Exp\cM_S$ (the $S$-relative versions of the 
Grothendieck ring with exponentials and its localization) are the settings where the motivic Fourier
transform and the Hrushovski--Kazhdan motivic Poisson summation formula naturally
live, and the corresponding assemblers $\cC^{KExp}_S$ and $\cC^{KExp\cS}_S$ constructed in Propositions~7.6 and 7.7
are where their categorical lifts will take place (Theorem~7.12 and Theorem~7.14.1). 

\smallskip

In Section~7.8 we discuss a natural class of motivic measures
on the Grothendieck rings of varieties with exponentials, in the
case over a finite field, given by the exponential sums, and the
associated zeta function. We present in Proposition~7.8.1 an
Euler product expansion generalizing the usual case of the
Hasse--Weil zeta function. 

\smallskip

In Section~7.9 and 7.10 we consider a different approach
to the categorification of Grothendieck rings, through the
construction of associated categories of motives through
Nori's formalism. In Section~7.9 we recall the general
functioning of Nori diagrams and Nori motives. As we
recall in Section~7.10, a motivic category related to the 
Grothendieck ring of varieties with exponentials is provided
by the category $MotExp(K)$ of exponential motives of Fres\'an and Jossen [FreJo20],
which has a unique ring homomorphism 
$$ \chi: KExp_K \to K_0(MotExp(K)). $$
As we discuss in Section~7.9.1, the Nori approach  also leads to
the construction of a homotopy-theoretic spectrum, through the
construction of a $\Gamma$-space associated to the Nori category
and an associated category of summing functors [Carl05], [Se73]. The spectrum
obtained in this way is a delooping of the infinite loop space
associated to the $K$-theory of the Nori category. 

\smallskip

In Section~7.11 we return to our main goal of Section~7, which is to
investigate categorical aspects of the motivic height zeta functions, and
we recall the setting and the main properties of the motivic Fourier transform,
following [ChamLoe15]. 

\smallskip

In Theorem~7.12 we prove that the motivic Fourier transform lifts to
a morphism of assemblers, and the identity in $KExp_V$ (with $V$
a finite dimensional $K$-vector space) satisfied by the square of
the Fourier transform is induced by a covering family in the assembler
$\cC^{KExp\cS}_V$ constructed in Proposition~7.6.

\smallskip

In Section~7.13 we recall the Hrushovski--Kazhdan motivic Poisson 
summation formula, following [ChamLoe15] and [HruKaz09], and in
Section~7.14 we present its categorification. In Theorem~7.14.1 we
prove that the motivic Poisson summation formula, seen as a
relation in $Exp\cM_k$, is induced by a covering family in a corresponding
assembler. 

\smallskip

The remaining part of Section~7 introduces the question of  a possible
categorification of the motivic height zeta function itself, and discusses
some more specific aspects of this question. The motivic height
zeta function is recalled in Section~7.15, along with the role of motivic Fourier 
transform and Poisson summation formula in establishing its properties,
following [ChamLoe15] and [Bilu18]. General methods for categorical
liftings of zeta functions using Witt rings are discussed in Section~7.16,
following [LMM19]. While these methods would apply to zeta functions
associated to good (exponentiable) motivic measures on Grothendieck
rings of varieties with exponentials, they cannot be applied to the motivic
height zeta function itself. We propose, however, that the multivariable
versions introduced in [Bilu18] on the basis of multivariable versions
of the Kapranov motivic zeta function and motivic Euler product
decomposition, may be more suitable for the purpose of a possible
categorification via assemblers.

\bigskip

\centerline{\bf 1. Homotopy spectra: a brief presentation} 

\medskip

 The notion of {\it spectra} in homotopy theory (which we refer to
 as ``homotopy spectra") evolved 
through several stages. Below we will briefly describe two
of them: {\it sequential spectra}, which can be considered as the
first stage (fundamental definition), and { \it $\Gamma$--spaces}, the construction method
which we will mostly use.
The latter is a powerful method (introduced
by Segal in [Se74]) for generating spectra from data consisting of 
categories with a zero-object and a categorical sum. 

\medskip

{\bf 1.1. Sequential spectra}  (see [SpWh53]).  In order
to define them, we need to introduce the {\it smash product} $\wedge$  in the category of based
(or pointed) sets. 
It can be defined by
$$
(X,*_ X)\wedge (Y,*_ Y) = ((X\times Y)/(X\times *_Y)\cup (*_X\times Y), \bullet ),
$$
where the base point $\bullet$ in the smash product is the image of the union
of two ``coordinate axes'' $(X\times *_Y)\cup (*_X\times Y)$ after the contraction
of this union. Similarly, the smash product of pointed simplicial sets is the quotient
of $X\times Y$ that collapses the pointed simplicial subset 
$(X\times *_Y)\cup (*_X\times Y)$.

\smallskip
A {\it sequential spectrum} $E$ is a sequence of based
simplicial sets $E_n$, $n=0, 1, 2, \dots$ (see 0.3 above), and the structure maps
$$
\sigma_n: \Sigma E_n:= S^1\wedge E_n \to E_{n+1}.
$$

The {\it sphere (sequential) spectrum $\bold{S}$} consists of simplicial sets $S^n:= S^1\wedge \dots \wedge S^1$ 
and identical $\sigma_n$'s.

\smallskip

The smash product can be extended to the 
category of sequential spectra itself, however there it loses the commutativity property
(it no longer forms a symmetric monoidal category). One can form a commutative and 
associative smash product only after formally inverting stable equivalences, that is, 
by passing to the homotopy category. We will not review the construction here, but
we refer the reader to [EKMM97], [HoShiSmi00], [Sc12]. We just remark
that, when one uses $\Gamma$-spaces as a construction method for spectra, as 
we discuss in the next subsection, the smash product has a very transparent 
description as shown in [Ly99].

\medskip
 {\bf 1.2. $\Gamma$--spaces.} $\Gamma$--spaces were introduced by Graham Segal in [Se74]. 
 Here we reproduce arguments of Lydakis in [Ly99] and mostly keep his notation. Let $\Gamma$ 
 ($\Gamma^{op}$ in [Ly99])
 be the category of based finite sets
$n^+ := \{0,1, \dots , n\}$ (formerly called $[n]$), with base preserving maps as morphisms. So here
we {\it do not restrict morphisms by nondecreasing maps}.
\smallskip
A  {\it $\Gamma$--space} $E$ is a functor from $\Gamma$ to based simplicial sets sending
$0^+$ to the point.
\smallskip

$\Gamma$--spaces themselves are objects of the category denoted $\Cal{G}\Cal{S}$
in [Ly99], morphisms in which are natural functors (we omit a precise description).

\smallskip
Following Section~2 of [Ly99], we will call based simplicial sets
simply {\it spaces}. A space is called {\it discrete} if its simplicial set is constant.
\medskip

{\bf 1.3. Theorem.} {\it One can define a functor    $\wedge : \Cal{G}\Cal{S} \times \Cal{G}\Cal{S} \to \Cal{G}\Cal{S}$
such that there is a canonical isomorphism of functors in three variables
$$
\Cal{G}\Cal{S} (F\wedge F^{\prime}, F^{\prime\prime})  \cong
\Cal{G}\Cal{S} (F, \, \roman{Hom} \, (F^{\prime}, F^{\prime\prime}))
$$
where we denote by $\Cal{G}\Cal{S} (F, F^{\prime})$ the set of morphisms $F\to F^{\prime}$
in $\Cal{G}\Cal{S}$.

The category of\quad $\Gamma$--spaces endowed with smash--product $\wedge$ is a symmetric monoidal category.}

\medskip

{\bf Sketch of Proof.} Consider the category of $\Gamma\times\Gamma$--spaces
$\Cal{G}\Cal{G}\Cal{S}$. Its objects are pointed functors 
$\Gamma\times\Gamma \to \Cal{S}$. Define also the external smash product
$F \overline{\wedge} F^{\prime}$ of two $\Gamma$--spaces $F, F^{\prime}$ as the functor
that sends  $(m^+, n^+)$ to $F(m^+)\wedge F^{\prime}(n^+)$.

\smallskip

Then one can check that
$$
\Cal{G}\Cal{S}(F, \roman{Hom}\,(F^{\prime}, F^{\prime\prime})) \cong
\Cal{G}\Cal{G}\Cal{S}(F \overline{\wedge} F^{\prime}, RF^{\prime\prime}),
$$
where $R$ is a functor $\Cal{G}\Cal{S} \to \Cal{G}\Cal{G}\Cal{S}$, determined
by the map $\Gamma\times\Gamma\to\Gamma$ sending $(m^+,n^+)$ to $m^+\wedge n^+$,
by setting $RF''(m^+,n^+)=F''(m^+\wedge n^+)$. 

\smallskip

The $q$-simplices of
$\roman{Hom}\,(F^{\prime}, F^{\prime\prime})(m^+)$ are given by
$\Cal{G}\Cal{S}(F\wedge (\Delta^q)^+, F'(m^+\wedge))$, where
$F'(m^+\wedge): n^+\mapsto F'(m^+\wedge n^+)$.

\smallskip

After that one can prove that $R$ has a left adjoint functor 
$L: \Cal{G}\Cal{G}\Cal{S} \to \Cal{G}\Cal{S}$. Namely, 
$LF^{\prime\prime\prime}$        is the colimit
 of  $F^{\prime\prime\prime}(i^+,j^+) $
over all morphisms $i^+\wedge j^+ \to n^+$.
\smallskip

This is essentially the statement of Theorem 1.3.

\smallskip
For many more details, see [Ly99], especially Theorem~2.2.
\medskip
{\bf 1.4. Definition.} {\it The sphere spectrum $\bold{S}$ is  the unit object 
 in the symmetric monoidal category of $\Gamma$--spaces.}
 
 \medskip
 
 This is the functor $\bold{S}$ given by the inclusion of the
 category $\Gamma$ of pointed sets in the category of pointed
 simplicial sets.

\medskip

Here is a more detailed description of  $\bold{S}$. For any $n^+$ we can define
the representable $\Gamma$--space $\Gamma^n$ by 
$$
\Gamma^n(m^+) := \Cal{G}\Cal{S} (n^+,m^+) .
$$
This satisfies $\Cal{G}\Cal{S}(\Gamma^n\wedge (\Delta^q)^+, F)=F(n^+)$.
From this it follows that $\bold{S}$ is canonically isomorphic to $\Gamma^1$, and
is in fact the sphere sequential spectrum defined in 1.1 above (see Definitions~2.5 and 2.7,
Lemma~2.6 and Proposition~2.8 of [Ly99]).

\medskip
In particular, this construction suggests to consider homotopical enrichments of arithmetics passing through 
\smallskip

\centerline{{ \it ring of integers} $\bold{Z}$  $\Longrightarrow$ {\it the sphere spectrum} $\bold{S}$}

\smallskip

This is in particular justified by the fact that $\bold{S}$ is a ring spectrum and 
$\bold{Z}$, viewed as the Eilenberg-MacLane spectrum $H \bold{Z}$, is a
ring spectrum over $\bold{S}$. This idea leads naturally to one of the approaches
to $\bold{F}_1$ geometry (geometry ``below $\roman{Spec}(\bold{Z})$") as
based on the sphere spectrum $\bold{S}$, as in [ToVa09].

\smallskip

This also poses a challenge: discover a new information in number theory
using the independently-developed machinery of homotopy theory.
\medskip

{\bf 1.5. $\Gamma$--spaces: from categories to spectra.}
Our main reason for introducing the properties of $\Gamma$-spaces lies in the
fact that these provide a very useful general method for constructing spectra.

\smallskip

A $\Gamma$-space $F$, which is a functor from the category $\Gamma$ of finite 
pointed sets to that of pointed simplicial sets, extends to an endofunctor of the category of 
pointed simplicial sets, by defining $F(K)$ as the coend
$$ F(K)= \int^{n^+\in \Gamma} K_n \wedge F(n^+), $$
with $K_n$ the $n$-th skeleton of $K$. There are 
natural assembly maps $K\wedge F(K')\to
F(K\wedge K')$. 

\smallskip

Thus, one obtains a spectrum $\bS(F)$ associated to the $\Gamma$ space $F$ 
by applying this functor to spheres, $\bS(F)_n:=F(S^n)$,
with structure maps $S^1\wedge F(S^n)\to F(S^{n+1})$.

\smallskip

In particular,  as originally described by Segal in [Se74], $\Gamma$-spaces can
be used to associate a spectrum to a category $\cC$ that has a categorical
sum and a zero object, via the construction of a $\Gamma$-space
$F_{\cC}$ obtained in the following way (see [Carl05]). 

\smallskip

Let $(X,x_0)$ be a finite pointed set. Consider the category $P(X)$ with objects
all the pointed subsets $A\subset X$ and morphisms given by pointed inclusions.
A {\it summing functor} $\Phi: P(X) \to \cC$ is a functor satisfying
$\Phi(A\cup A')=\Phi(A)\oplus \Phi(A')$ for $A,A'$ in  $P(X)$ with $A\cap A'=\{ x_0 \}$,
and $\Phi(\{ x_0 \})=0$. Let $\Sigma_\cC(X)$ denote the category of summing
functors with morphisms given by invertible natural transformations.

\smallskip

The reason for using natural isomorphisms instead of arbitrary natural
transformations as morphisms in $\Sigma_\cC(X)$ is in order to ensure
that the resulting nerve (see below) has interesting topology: using all
natural transformations would lead to a contractible space (the category
of summing functors has an initial object, since $\cC$ has a zero object).

\smallskip

Recall that the nerve $\Cal{N}$ of a category $\Cal{A}$
is the simplicial set $\Cal{N}\Cal{A}$ whose vertices are
indexed by objects $A$ of $\Cal{A}$; 1--simplices by morphisms $A_1 \to A_2$ of $\Cal{A}$;
and generally, $n$--simplices by sequences of morphisms $A_1\to \dots \to A_{n+1}$.
Faces and degenerations (see 0.3 above) are determined via composition
of morphisms in a pretty obvious way.

\smallskip

The $\Gamma$--space $F_{\cC}$ is then determined by assigning to a finite 
pointed set $(X,x_0)$ the pointed simplicial set given by the nerve 
$\Cal{N}(\Sigma_\cC(X))$ of $\Sigma_\cC(X)$, with associated spectrum $\bS(F_{\cC})$.

\smallskip

The intuition behind this construction is the following: the category of summing
functors provides a {\it delooping} of the infinite loop space given by (a completion
of) the classifying space of $\cC$, see [Carl05] for a more detailed discussion. 

\smallskip

While this construction of spectra via $\Gamma$-spaces only gives
rise to connective spectra, it is shown in [Tho95] that all connective spectra can be
obtained in this way.

\smallskip

In this paper, the main applications 
of this general method for the construction of spectra via $\Gamma$-spaces will
appear in Section~3.3 and in Section~7.9.1. The first is the application of
$\Gamma$-spaces that occurs within Zakharevich's theory of  
assemblers (assembler categories), where the $\Gamma$-spaces method provides
a key step in the construction of the corresponding spectrum. We will review this in
Section~3.3 below and we will apply it then several times to different settings, in 
Sections~6 and then more explicitly in Proposition~7.2, Theorem~7.4, Propositions~7.5 and 7.6,
Theorems~7.12 and 7.14.1. The second occurrence of $\Gamma$-spaces will
be an application of the construction of spectra from categories, as recalled in this subsection, 
applied to Nori diagrams in Section~7.9.1.

\bigskip

\centerline{\bf 2.  Diophantine equations:}
\smallskip
\centerline{\bf distribution of rational points on algebraic varieties} 

\medskip

{\bf 2.1. Diophantine equations and heights.} We will be studying here 
how fast the number of solutions of a system of equations can grow,
if one first restricts the counting to solutions of {\it height} $\le H$, and then
lets   $H$ grow.

In order to define heights over general algebraic number fields, we need
the following preparations.

Let $K$ be a number field, $\Omega_K = \Omega_{K,f}\sqcup\Omega_{K,\infty}$
the set of its places $v$ represented as the union of finite and infinite ones, and 
let $K_v$ denote the respective completion of $K$.

 For  $v\in \Omega_{K,f}$, denote by
  $\Cal{O}_v$, resp.~$m_v$,
the ring of integers of $K_v$, resp.~its maximal ideal. 
For a uniformizer $\pi_v$ of $\Cal{O}_v$, $|\pi_v|_v^{-1}$ is 
the size of the residue field $\Cal{O}_v/m_v$.
The Haar measure $dx_v$ on $K_v$
is normalised in such a way that the measure of $\Cal{O}_v$ becomes 1. Moreover,
for an archimedean $v$, the Haar measure will be the usual Lebesgue measure, if $v$ is real,
and for complex $v$ it will be induced by Lebesgue measure on $\bold{C}$, for which
the unit square $[0,1] + [0,1]i$  has volume 2. 

Let the map $|.|_v : K_v \to \R^*_{\ge 0}$ be defined by the condition
$d(\lambda x)_v = |\lambda |_v dx_v$. In particular, 
$|.|_v$ is the usual
absolute value for $v$ real, and its square for $v$ complex.

Then for any $\lambda \in K^*$ we have the following product formula:
$\prod_v |\lambda |_v = 1$.

Let $\bold{P}^n$ be a projective space with a chosen system of
homogeneous coordinates $(x_0: x_1: $\dots$  :x_{n-1}:x_n)$, $n\ge 1$.
Then we can define {\it the exponential Weil height} of a point 
$p = ( x_0(p): \dots :x_n(p)) \in \bold{P}^n(K)$ as
$$
h(p) := \prod_{v\in \Omega_K} \roman{max} \left\{ |x_0(p)|_v, \dots , |x_n(p)|_v \right\} .
$$

Because of the product formula, the height does not change if we replace coordinates
$(x_0:\dots : x_n)$
by $(\lambda x_0:\dots : \lambda x_n)$, $\lambda \in K$.

\smallskip

Notational remark: it is more common in the literature to use a different notation
$H(p)$ for the exponential height and $h(p)$ for the logarithmic height, but
we will be only considering the exponential height and we will simply 
use $h(p)$ for all instances of height functions in this paper. 

\medskip

 {\bf 2.2. Height zeta functions.} Let now $(U,L_U)$ be a pair consisting of a projective
variety $U$ over $K$ and an ample line bundle $L_U$ on it. 
Let $V\subset U$ be a locally closed subvariety of $U$, also defined over $K$.
Then we can define
the height function $h_{L_V}(p)$ on $p\in V(K)$ using the same formula as above,
but this time interpreting $(x_i)$ as a basis of sections in $\Gamma (V, L_V)$.

If $h^{\prime}_{L_V}$ is another height, corresponding to a different choice of the basis
of sections, then there exist two positive real constants $C, C^{\prime}$ such that for all $x$,
$$
C h_{L_V}(x) \le h^{\prime}_{L_V}(x) \le C^{\prime} h_{L_V}(x).
$$

Now define the height zeta--function 
$$
Z(V, L_V, s) : = \sum_{x\in V(K)} h_{L_V}(x)^{-s} .
$$ 

\medskip

{\bf 2.3.  Convergence boundaries}. 
Denote by $\sigma (V, L_V) \in \R$ the greatest lower bound of the set of
positive reals $\sigma$ for which
$Z(V,L_V,s)$ absolutely converges if $Re\, s\ge \sigma$. 

We will call
 $\sigma (V, L_V)$ the respective {\it convergence boundary}.
  
  Clearly, it is
 finite  (because this is so for projective spaces), and non--negative whenever 
 $V(K)$ is infinite.
 
 For example, one has $\sigma(\P^n,\Cal{O}(m))=(n+1)/m$ (Section~1.3 of
 [BaMa90]). 

Intuitively, we may say that $V$ contains ``considerably less'' $K$--points than $U$,
if 
$$
\sigma (V,L_V) < \sigma (U,L_U).
$$
and ``approximately the same'' amount of  $K$--points, if
$$
\sigma (V,L_V) = \sigma (U,L_U),
$$
in the sense that the rate of growth of the counting of such
points (weighted by height) is lower or the same.

\smallskip

 {\bf 2.3.1. Example: accumulating subvarieties}. 
Let $V\subset U$ be a Zariski closed subvariety over $K$. 
As in Section~3.2.3 of [Cham10], we say that $V$ is a 
{\it strongly accumulating subvariety} if the fraction
$$ \frac{\roman{card}\left\{ x\in V(K)\, \biggr| \, h_{L_V}(x) \leq H \right\} }
{\roman{card}\left\{ x\in U(K)\, \biggr| \, h_{L_U}(x) \leq H \right\} } \rightarrow 1 $$
as $H\to \infty$, and a {\it weakly accumulating subvariety} if
$$ \roman{liminf}_H \frac{\roman{card}\left\{ x\in V(K)\, \biggr| \, h_{L_V}(x) \leq H \right\} }
{\roman{card}\left\{ x\in U(K)\, \biggr| \, h_{L_U}(x) \leq H \right\} } >0. $$

\smallskip

Clearly, then  $\sigma (V,L_V) = \sigma (U,L_U)$. 

\smallskip

We will return to these notions of strongly and weakly accumulating subvarieties in
Proposition~4.4.1, where we construct assemblers whose associated Grothendieck
groups detects the presence of such subvarieties through its 
scissor-congruence relations.

\medskip

We will now describe a categorical environment appropriate for
describing various versions of accumulation and connecting
distributions with spectra.

\bigskip

 \centerline{\bf 3. Rational points, sieves, and assemblers}
 \medskip
 
 {\bf 3.1.  Grothendieck topologies, sieves, and assemblers.} (See [Za17a], [Za17b], [MaMar18]).
  Let $\Cal{C}$ be a category with a unique initial object $\emptyset$. Two morphisms
  $f_1 : U_1\to U$ and $f_2 : U_2 \to U$ are called {\it disjoint},  if $U_1\times_U U_2$
  exists and is $\emptyset.$
  \smallskip
  
  Notice that if $U_1\times_U U_2$ exists, then $U_2\times_U U_1$ also exists,
  and these two relative products are canonically isomorphic, through a
  unique isomorphism.
  
  \smallskip
A {\it sieve} in $\Cal{C}$ is a full subcategory
$\Cal{C}^{\prime}$ such that if $f : V \to U$ is a morphism in $\Cal{C}$,
and $U$ is an object of $\Cal{C}^{\prime}$, then $V$
 is also an object of  $\Cal{C}^{\prime}$.
 
 \smallskip
 
 Fixing an object $U$ in $\Cal{C}$, can apply this notion also to the category of morphisms 
 $f: V\to U$, or in other words to the ``overcategory"  $\Cal{C}/U$. Functors between 
 the overcategories  $f_*: \Cal{C}/W  \to  \Cal{C}/V$ induced by composition with $f: W\to V$
 induce functors (in the opposite direction) between sieves in the respective overcategories.
 In fact, any functor $f: \Cal{C}\to \Cal{D}$ between two categories induces a pullback
 map on sieves: if $\Cal{S}\subseteq \Cal{D}$ is a sieve, then the full subcategory $f^{-1}(\Cal{S})$
 is a sieve in $\Cal{C}$.
 
\smallskip
This notion is convenient in order to define   a {\it Grothendieck topology} on $\Cal{C}$:
it is {\it a collection of sieves} $\Cal{J}(U)$ in $\Cal{C}/U$, one for each object $U$ of $\Cal{C}$,
satisfying three axioms: 
\smallskip

{\it 

(i) Any morphism $f: W\to V$ lifts to the map $f^*:  \Cal{J}(V) \to \Cal{J}(W)$.

\smallskip

(ii) The full overcategory \  $\Cal{C}/U$ belongs to    $\Cal{J}(U)$ for any object of $U$ of $\Cal{C}$.

\smallskip

(iii) Let $\Cal{S}\in \Cal{J}(U)$ and $\Cal{T}$ be a sieve in $\Cal{C}/U$. If $f^*(\Cal{T})\in \Cal{J}(V)$ for 
all $f: V \to U$ in $\Cal{S}$, then $\Cal{T}\in \Cal{J}(U)$.

}
\smallskip
   
For any object $U$ of a category with Grothendieck topology $\Cal{C}$ (called also Grothendieck site)
we can define the notion of
  {\it covering family}: it is a collection of morphisms
$\{f_i :U_i\to U\, |\, i\in I\}$ such that the full subcategory of $\Cal{C}/U$
containing all morphisms in $\Cal{C}$ factoring through one of the $f_i$'s
belongs to the initial collection of sieves $\Cal{J}(U)$. 

\smallskip

A family $\{ f_i: U_i \to U \}$ is disjoint if $f_i$ and $f_j$ are disjoint for all $i\neq j$.
As will become clear in the result of [Za17a] that we recall in Theorem~3.4 below,
heuristically one should think of disjoint covering families as the categorification of scissor-congruence
relations. 

\medskip

{\bf 3.2. Assemblers.} 
An  assembler (also referred to as assembler category) 
is a small category $\Cal{C}$ with a Grothendieck
topology with the following properties: 

\noindent (a) $\Cal{C}$ has an initial object $\emptyset$, for which the empty family
is a covering family; 

\noindent (b) all morphisms in $\Cal{C}$
must be monomorphisms;

\noindent (c) any two disjoint finite covering families
must admit a common refinement which is also a finite disjoint
covering family.
\smallskip
Assemblers themselves form a category, in which a morphism
is a functor continuous in the respective Grothenieck topologies, sending
initial object to initial object, and disjoint morphisms to disjoint morphisms.
\smallskip
Let $\Cal{C}$ be a Grothendieck site. Denote by $\Cal{C}^{\circ}$ its full
subcategory of non-initial objects. 
\smallskip

If we have a family of assemblers $\{ \Cal{C}_x\}$ numbered by
elements $x$ of a set $X$, we will denote by $\bigvee_{x\in X}  \Cal{C}_x$
the category whose non-initial objects are $\bigsqcup_{x\in X}\roman{Ob} \Cal{C}^{\circ}_x$
and to which  one initial object  is formally added.

\smallskip

A simplicial assembler is a functor from the simplex category to the category
of assemblers (a simplicial object in the category of assemblers), see Definition~2.14 of [Za17a].

\medskip

{\bf 3.3. From assemblers to $\Gamma$--spaces  and  spectra.} (See [Za17a], also reviewed in [MaMar18], Sec.~4.4).
Starting with an assembler $\Cal{C}$, we can construct the following category $\Cal{W}(\Cal{C})$:

\smallskip
(a) {\it An object} of $\Cal{W}(\Cal{C})$ is a map $I\to Ob\, (\Cal{C})$ where $I$ is a finite set,
and the map lands in non-initial objects. We may write it as $\{A_i\, |\, {i\in I}\}, A_i\in Ob\, (\Cal{C})$.

\smallskip

(b) {\it A morphism} $f: \{A_i\}_{i\in I} \to \{B_j\}_{j\in J}$ consists of a map of finite sets $f: I\to J$
and a family of morphisms $f_i: A_i \to B_{f(i)}$ such that $\{f_i : A_i\to B_j : i\in f^{-1}(j)\}$
for each $j\in J$ is a disjoint covering family.
\medskip

{\bf 3.3.1. Proposition. (See [Za17a], Prop. 2.11.)} {\it (i) All morphisms in $\Cal{W}(\Cal{C})$
are monomorphims.

\smallskip

(ii) If $\Cal{C}$ has all pullbacks (i.e.  is closed), then $\Cal{W}(\Cal{C})$ is closed as well.

\smallskip

(iii) Given a family of assemblers $\{\Cal{C}_x\,|\, x\in X\}$ indexed by elements of a set $X$,
denote by $\oplus \,\Cal{W} (\Cal{C}_x)$ the full subcategory of \  $\prod \,\Cal{W} (\Cal{C}_x)$
whose objects are  families of objects of $\Cal{C}_x$ for which all but finitely many of them
are indexed by $\emptyset$. Consider the functor 
$$
P: \Cal{W} (\bigvee_{x\in X} \Cal{C}_x) \to \prod_{x\in X} \Cal{W}( \Cal{C}_x)
$$
induced by the morphisms of assemblers $F_y : \vee_{x\in X}\, \Cal{C}_x \to \Cal{C}_y$
that map each $\Cal{C}_x$ to the initial object if $x\ne y$ and identically to $\Cal{C}_y$
if $x = y$.

\smallskip

Then $P$ induces an equivalence of categories
$$
 \Cal{W} (\bigvee_{x\in X} \Cal{C}_x) \to \bigoplus_{x\in X} \Cal{W}( \Cal{C}_x).
$$
}

\smallskip

{\bf Proof.} The proof is as in [Za17a], Proposition~2.11. We recall
it here for the reader's convenience. 
The main part of our argument is a detailed study of pullback
squares (see the diagram in the proof of Proposition~2.11 of [Za17a]).

\smallskip

Start with a morphism $f: \{A_i\}_{i\in I}  \to \{B_j\}_{j\in J}$, and two more
morphisms
$$
g,h :\   \{C_k\}_{k\in K}   \to    \{A_i\}_{i\in I}.
$$
In  order to show that $f$ is a monomorphism we must check that if $fg=fh$, then $g=h$.
Choose any $k\in K$, and look at the respective commutative square
with vertices $C_k, A_{g(k)}, A_{h(k)}, B_{fg(k)}$. We must have $g(k)= h(k)$,
because otherwise $f_{g(k)}$ and $f_{h(k)}$ are disjoint.

\smallskip

The remaining statements easily follow from these remarks.   $\blacksquare$

\medskip

{\bf 3.3.2.  Move to $\Gamma$--spaces, spectra, and $K$--theory. (Cf. Sec. 2.2 of [Za17a]).} 
To move from here to $\Gamma$--spaces,
start with the category $\Gamma^0$ whose objects are finite sets $n^+ : = \{0,1,2, \dots , n\}$
with base point $0$, and morphisms are all maps sending $0$ to $0$ as morphisms, 
as in Subsection~1.2 above.

\smallskip

 If $X$
is a pointed set, and $\Cal{C}$ an assembler, we can construct
the assembler $X\wedge \Cal{C} := \bigvee_{x\in X^{\circ}} \Cal{C}_x$,
where $X^{\circ} := X\setminus \{*\}$ .   


\smallskip

For each $n$, we have a simplicial assembler $S^n\wedge \cC$. Then 
$\Cal{W}(S^n\wedge \Cal{C})$ is a simplicial category, hence its nerve 
$\Cal{N}\Cal{W}(S^n\wedge \Cal{C})$ is a bisimplicial set. This is then
the $n$-th space of the associated spectrum,
with structure maps induced by $X\wedge \Cal{N}\Cal{W}(S^n\wedge \Cal{C})\to
 \Cal{N}\Cal{W}(X \wedge \Cal{C})$:
 $$
 X\wedge \Cal{N}\Cal{W} (\Cal{C}) \cong \bigvee_{X^{\circ}} \Cal{N}\Cal{W} (\Cal{C}) \to
 \Cal{N} ( \bigoplus_{X^{\circ}}     \Cal{W} (\Cal{C})) \equiv \Cal{N}\Cal{W} (X\wedge \Cal{C}) .
 $$
 \smallskip
 
 I.~Zakharevich defines the symmetric spectrum $K(\Cal{C})$ as
 spectrum of simplicial sets in which the $k$--th space
 is given by the diagonal of the bisimplicial set
 $[n] \to \Cal{N}\Cal{W} ((S^k)_n \wedge \Cal{C})$
 and writes
 $$
 K_i(\Cal{C}) := \pi_i K(\Cal{C})\, .
 $$
 Note that this is not just a sequential spectrum as described in Section~1.1, but
 a symmetric spectrum, where the sequence of spaces is endowed with compatible
 actions of the symmetric group, see Definition~1.2.1 of [HoShiSmi00].
 Symmetric spectra allow for a good notion of smash product but require more
 care when taking homotopy groups (we refer the reader to [HoShiSmi00] and [Sc12]).
 
\smallskip

The construction above of the spectrum $K(\Cal{C})$ can be seen as a
case of the construction of spectra $\bS(F)$ via $\Gamma$-spaces $F$ 
that we recalled in Section~1.5. 
Indeed, the assignment 
$$ F: n^+ \mapsto F(n^+)= \Cal{N}\Cal{W}(n^+ \wedge \Cal{C}) $$
defines a $\Gamma$-space and $K(\Cal{C})$ is the spectrum 
(called $\bS(F)$ in Section~1.5) determined by this $\Gamma$-space.
This fact is explicitly used in the proof of Theorem 2.13 in [Za17a], which
we recall below, to identify generators and relations
of $K_0(\cC)$ in terms of noninitial objects of $\cC$ and morphisms of $\Cal{W}(\cC)$,
using the identification $K_0(\cC)=\pi_0(\bS(F))$, see also Section 4 of [Se74].
  
\medskip
The following result (Theorem 2.13 in [Za17a]) furnishes the first justification
of the intuition encoded in the word ``assembler''.

\medskip

{\bf 3.4. Theorem.} {\it Let $\Cal{C}$ be an assembler. Then $K_0(\Cal{C})$
is canonically isomorphic to the abelian group generated by (isomorphism
classes of) objects $[A]$ of $\Cal{C}$ modulo the family of relations 
indexed by finite disjoint covering families $\{ A_i \to A\, | \, i\in I\}$:
each such family produces the relation
$$
[A] = \sum_{i\in I} [A_i].
$$
}
\medskip

{\bf Sketch of proof.} The calculation of generators and relations
for $K_0(\Cal{C})$ can be reduced to making explicit the composition
of functors
$$
P^{-1}:\ \Cal{W}(\Cal{C}) \oplus \Cal{W}(\Cal{C}) \to \Cal{W}(\Cal{C} \wedge \Cal{C})
$$
and
$$
\mu :\, \Cal{W}(\Cal{C} \vee \Cal{C}) \to \Cal{W}(\Cal{C}).
$$

Here $P^{-1}$ sends a pair of objects $(\{A_i \, |\, i\in I\},   \{B_j\, |\, j\in J\}) $
to the object  $\{C_k \,|\, k\in I\cup J\}$, where for $k\in I$, resp. $k\in J$,
we put $C_k=A_k$, resp. $C_k=B_j$: cf. Proposition 3.3 above.
The map $\mu^{-1}$ is induced by the folding map of assemblers.

\smallskip

Finally, it remains to remark that relations in $\pi_0$ are generated
by morphisms $\{ A_i \to A\, | \, i\in I\} \to \{A\}$.   $\blacksquare$

\medskip

In most significant applications, the Grothendieck group $K_0(\cC)$
is also endowed with a ring structure. This is induced by a symmetric
monoidal structure on the assembler $\cC$, which in turn determines
an $E_\infty$-ring structure on the spectrum $K(\cC)$, see [Za17a].

\medskip

For most of the applications to constructions of assemblers and
associated spectra that we consider in this paper, we will be
constructing some explicit disjoint covering families in an assembler 
$\cC$. These covering families will provide categorical lifts of certain 
relations between classes in the corresponding Grothendieck ring 
$K_0(\cC)$. This will be the case, for example, in  Theorems~7.12 and 7.14.1.

\medskip
We now move to our main preoccupation here:
constructing assemblers related to the distributions of rational points on algebraic varieties.

\smallskip
We will first of all define formally certain {\it sieves} via point distribution.

\medskip

{\bf 3.5. Categories $\Cal{C}(U,L_U)$.} Let $U$ be a projective variety over $K$ and $L_U$ an
ample rank 1 vector bundle on $U$ over $K$.
\smallskip
By definition, objects of $\Cal{C}(U,L_U)$ are locally closed subvarieties $V\subset U$
also defined over $K$, and morphisms are the structure embeddings $i_{V,U}$, or simply $i_V$: $V\to U$.
Here  we did not mention $L$ explicitly, but it is natural to endow each $V$ by $L_V:=i_V^*(L_U)$.
\smallskip
Structure embeddings are compatible with these additional data so that we have
in fact structure functors $\Cal{C}(V,L_V) \to \Cal{C}(U,L_U)$ which make of each $\Cal{C}(V,L_V)$
a full subcategory of  $\Cal{C}(U,L_U)$ closed under precomposition, that is, {\it a sieve}.
\smallskip
We will call such $\Cal{C}(V,L_V)$ {\it geometrical sieves}, and now introduce
the {\it arithmetical sieves}  $\Cal{C}^{ar}(V,L_V)$ in the following way.

\medskip

 {\bf 3.5.1. Lemma.} {\it The family of those morphisms $i_{V,U}$ 
 as above, together with their sources and targets, for which
$$
0 < \sigma (V,L_V) < \sigma (U,L_U) ,
$$
forms a sieve in $\Cal{C}(U,L_U)$ denoted $\Cal{C}^{ar}(U,L_U)$.}

\smallskip

{\bf Proof.} If we have a two--step ladder of locally closed embeddings
$W\subset V \subset U$ such that $0 < \sigma (V,L_V) < \sigma (U,L_U)$
and $0 < \sigma (W,L_W) < \sigma (V,L_V)$, then of course
$0 < \sigma (W,L_W) < \sigma (U,L_U)$, so that the composition of these embeddings
is also a morphism in $\Cal{C}^{ar}(U,L_U)$. $\blacksquare$

\smallskip

Notice, that if $V(K)$ is a finite set, then $\sigma (V,L_V) =0$, but the converse is not true:
$ \sigma (V,L_V) =0$ for any abelian variety $V/K$ and for many other classes of $V$.
A complete geometric description of this class of varieties seemingly is not known.

\smallskip

Note also that the arithmetic sieves $\Cal{C}^{ar}(V,L_V)$ behave differently from the sieves
$\Cal{C}(V,L_V)$, in the sense that in general 
$$ \Cal{C}^{ar}(U,L_U) \cap \Cal{C}(V,L_V) \neq \Cal{C}^{ar}(V,L_V). $$

\medskip

{\bf 3.6.  Arithmetic assemblers.} Using sieves  $\Cal{C}^{ar}(U,L_U)$, we can easily
introduce the respective arithmetic assemblers $\Cal{C}_U$ that we describe below.

\smallskip

We first need to identify a good class of morphisms to use for the
construction of the assembler $\Cal{C}_U$.

\medskip

{\bf 3.6.1. Lemma.} {\it If $V\hookrightarrow U$ is a locally closed subvariety
satisfying $0< \sigma(V,L_V)< \sigma(U,L_U)$, then the complement
$W=U\smallsetminus V$ satisfies $\sigma(W,L_W)=\sigma(U,L_U)$.}

\smallskip

{\it Proof.} It suffices to check this by considering height zeta functions as 
Dirichlet series, and normalize heights in such a way 
that, for $x\in U(K)$, the height $h(x)$ does not depend on other choices
(e.g.~does not depend on $V$ and $L_V$, if $x\in V(K)$). Then the corresponding
zeta functions satisfy $Z(U,s)= Z(V,s) + Z(U\smallsetminus V,s)$. If there is a real
point $s$ where one of the two summands in the right-hand-side diverges and 
another converges, then $Z(U,s)$ will diverge at this point, and if both summands
converge at $s$ then $Z(U,s)$ will also converge. $\blacksquare$

\medskip

This shows that, if we want to have an assembler with non-trivial disjoint covering 
families, we cannot just use morphisms in $\Cal{C}^{ar}(V,L_V)$. We can, however,
relax this requirement by considering morphisms that are arbitrary finite
compositions of open embeddings in $\Cal{C}^{ar}(V,L_V)$ and arbitrary closed 
embeddings, where we do not impose the strict inequality between the 
abscissas of convergence of the respective height zeta functions.  
 
\medskip

{\bf 3.6.2. Proposition.} {\it Let $\cC_U$ be the category whose objects are
locally closed subvarieties $V$ of $U$ and whose morphisms are 
finite composites of open embeddings and closed embeddings, where
the open embeddings are required to be in $\Cal{C}^{ar}(V,L_V)$,
namely $i_{W,V}: W\hookrightarrow V$ with $W$ Zariski open in $V$
and $0<\sigma(W,L_W)< \sigma(V,L_V)$. 
The Grothendieck topology generated by the disjoint covering families of
the form 
$\{ W\hookrightarrow U, U\smallsetminus W \hookrightarrow U \}$,
where $W\subset U$ is open with $0<\sigma(W,L_W)< \sigma(U,L_U)$.
The category $\cC_U$ is an assembler, and the associated $K_0(\cC_U)$
group is generated by isomorphism classes $[V]$ of subvarieties of $U$
with the relations $[U]=[V]+[W]$ where $W=U\smallsetminus V$ is
open with $0<\sigma(W,L_W)< \sigma(U,L_U)$.}

\smallskip

{\it Proof.} The empty scheme is the initial object. Both open and
closed embeddings are monomorphisms. 
The category has pullbacks as closed embeddings are preserved 
under pullbacks and open embeddings
$i_{W,V}: W\hookrightarrow V$ satisfying $0<\sigma(W,L_W)< \sigma(V,L_V)$
are also preserved, since the intersection $W=W_1\cap W_2$ of open
subvarieties satisfying $0<\sigma(W_i,L_{W_i})< \sigma(V,L_V)$ will also
satisfy the strict inequality $0<\sigma(W,L_W)< \sigma(V,L_V)$. This
ensures that any two finite disjoint covering families have a common refinement.
The finite disjoint covering families are determined by sequences 
$$ U \supset V_0 \supset V_1 \supset V_2 \cdots \supset V_n $$
of Zariski open subvarieties with $0<\sigma(V_{i+1},L_{V_{i+1}})< \sigma(V_i ,L_{V_i})$,
by taking $\{ f_i: X_i \hookrightarrow U \}$ with $X_0=U\smallsetminus V_0$ and
$X_i=V_i \smallsetminus V_{i+1}$, with $V_{n+1}=\emptyset$. 
The Grothendieck group $K_0(\cC_U)$ is then generated by the objects of $\cC_U$ 
with the scissor-congruence relations determined by the disjoint covering families.
This means that $K_0(\cC_U)$ is generated by the $[V]$ for locally closed subvarieties 
$V\subset U$ with relations $[U]=[V]+[W]$ for $W=U\smallsetminus V$ 
open with $0<\sigma(W,L_W)< \sigma(U,L_U)$.
$\blacksquare$

\medskip

{\bf 3.6.3. Example: arithmetic stratifications.} An arithmetic stratification of $U$ with respect to the line bundle $L_U$,
in the sense of [BaMa90], consists of a descending sequence of Zariski open subvarieties
$U\supset V_0 \supset V_1 \supset V_2 \cdots$, where $V_{i+1}$
is the maximal open subset of $V_i$ with  
$0<\sigma(V_{i+1},L_{V_{i+1}})< \sigma(V_i ,L_{V_i})$. Different choices of the line
bundle give rise to different arithmetic stratifications. An arithmetic stratification defines
a disjoint covering family in the assembler $\cC_U$ by taking $X_i=V_i \smallsetminus V_{i+1}$.

\smallskip

The remaining two sections sketch two diverging paths leading from distribution
of $K$--rational points $U(K)$ in $U$ to various further versions of the
arithmetic sieve  $\Cal{C}^{ar}(U,L_U)$.

\smallskip

In Section~4 we look at the more narrow class of varieties, Fano varieties,
for which more precise data about behaviour of some heights are known,
or at least, conjectured.

\smallskip

In Section~5 (restricted to $K=\bold{Q}$), we consider rational points as a subset of {\it adelic} points,
and try to go beyond heights by using new tools for studying the distribution
of adelic points    $U(A_{\bold{Q}})$ themselves.

\medskip

\bigskip

 \centerline{\bf 4. Anticanonical heights and points count}
 \medskip

{\bf 4.1. Anticanonical heights: dimension one.} Let $(U, L_U)$ be as above
a pair consisting of a variety and ample line bundle defined over $K$, $[K:\bold{Q}] < \infty$.
Choose a exponential height function $h_L$, and set for $B\in \R_+$
$$
N(U,L_U, B) := \roman{card} \left\{x\in U(K)\, \biggr| \, h_L(x) \le  B \right\}.
$$

In Section~3, we based the definition of an arithmetical sieve upon
an intuitive idea that $i_V: V \to U$ belongs to this sieve, {\it  if the number of $K$--points
on $U$ is ``considerably less" that such number on $V$.} To make this idea precise,
we used convergence boundaries.

Below, we will use a considerably more precise counting of points, in order to define
subtler sieves on a more narrow class of varieties $U$, using {\it counting functions themselves} $N(U,L_U, B)$
in place of convergence boundaries.

 Start with one--dimensional $U$.

If $U$ is a smooth irreducible curve of genus $g$, {\it with nonempty set $U(K)$}, we have
the following basic alternatives:
$$
g=0:  \quad    U= \bold{P}^1, \quad  L_U = -K_V, \quad   N(U,L_U, B) \sim cB.
$$
$$
g=1: \quad U=\text{ an elliptic curve with Mordell--Weil rank}\,\, r, \quad N(U, L_U. B) \sim c \,
 (\log B)^{r/2}.
$$
$$
g>1:  \quad N(U, L_U,B) = const,\text{ if $B$  is  big enough.}
$$

A survey of  expected typical behaviours of multidimensional analogs 
can be found in the Introduction to [FrMaTsch89].

Below, our attention will be focussed upon {\it Fano varieties}, that is, varieties
with ample anticanonical bundle  $\omega_V^{-1}=-K_V$, as a wide generalisation
of the one--dimensional case $g=0$.

\smallskip

In the case of $\bold{P}^n$ (the simplest Fano variety), it was shown in [Scha79] that
$$ N(\bold{P}^n, \Cal{O}_{\P^n}(1), B)\sim c' \, B^{n+1} \ \   \text{ and } \ \ 
 N(\bold{P}^n, -K_{\bold{P}^n}, B)\sim c \, B, $$
 where $-K_{\bold{P}^n}=\Cal{O}_{\P^n}(n+1)$.
 The constant $c'$ here is expressible in terms of basic invariants of the number field $K$:
 the class number, the value $\zeta_K(n+1)$ of the Dedekind zeta function, the number of 
 real and complex embeddings of $K$, the regulator, see [Scha79]. As mentioned in Section~2.3
 above, for the line bundle $\Cal{O}_{\P^n}(1)$, the abscissa of convergence is 
 $\sigma(\P^n,\Cal{O}_{\P^n}(1))=n+1$ (see Section~1.3 of [BaMa90] and
 Section~4.1.1 of [Cham10]).

\medskip

{\bf 4.2. Anticanonical heights for Fano varieties}. The  most precise conjectural asymptotic formula for Fano 
varieties (or Zariski open subsets of them) with dense $U(K)$ has the form
$$
N(U, -K_U, B) = c B (\log B)^t + O(B (\log B)^{t-1}),
\quad t:= \roman{rk}\, \roman{Pic}\, U -1 ,
$$
with $\roman{Pic}\, U$ the Picard group.

\smallskip

Note that in the case of $\P^n$ the Picard group is $\roman{Pic}\,\P^n\simeq \Z$, by
$m\mapsto \Cal{O}_{\P^n}(m)$, so that $t=0$, but the $(\log B)^t$ term already occurs
in the case of products of projective spaces. 

\smallskip

This conjectural asymptotic formula certainly is {\it wrong} for many 
subclasses of Fano varieties. On the other hand, it is 
\medskip
 (i) {\it stable under direct products; }
 \smallskip
 
 (ii) {\it compatible with predictions
of Hardy--Littlewood for complete intersections; }

\smallskip

(iii) {\it true for
quotients of semisimple algebraic groups modulo parabolic subgroups.}

\smallskip

We reproduce below some of the arguments of [FrMaTSCH89], Sections 1--2,
proving these statements.

\medskip

{\bf 4.2.1. Direct products.} Let $(U,L_U)$ and $(V,L_V)$ be data as above.
We must study the behaviour of $N(U\times V, L_{U\times V})$
where   $L_{U\times V} := pr_U^*(L_U)\otimes pr_V^*(L_V)$.
\smallskip
In this context, we may restrict ourselves by consideration
of exponential heights satisfying the exact equality
$$
h_{L_{U\times V}}(x,y) = h_{L_U}(x) h_{L_V}(y).
$$
This will be applied  to the case of anticanonical heights.

\smallskip

We will now change notation without referring anymore to the specical properties of heights
on Fano varieties etc, as in [FrMaTsch89].

\smallskip

Consider two infinite families of  nondecreasing real numbers indexed by $1,2, \dots$:  $\{\lambda_i\}$
and $\{\mu_j\}$. We allow each number to be repeated several times, so that they can have finite
mulltiplicities. We can then form a new family: $\lambda\mu := \{\lambda_i\mu_j\}$
again ordered nondecreasingly. Put $N_{\lambda} (B) := \roman{card}\,\{i\,|\, \lambda_i < B\}$,
and similarly for $N_{\mu}, N_{\lambda\mu}$.

\smallskip

Now ``stability under the direct product'' from 4.2 is a consequence of the following (see
Proposition~2 of [FrMaTsch89]).
\medskip
{\bf 4.2.2. Lemma.} {\it If
$$
N_{\lambda} (B) = c_{\lambda} B {\roman{log}}^s B    + O(B {\roman{log}}^{s-1} B),
$$
$$
N_{\mu} (B) = c_{\mu} B {\roman{log}}^r B    + O(B {\roman{log}}^{r-1} B),
$$
then
$$
N_{\lambda\mu} (B) =C\, c_{\lambda}c_{\mu} B {\roman{log}}^{r+s+1} B    + O(B {\roman{log}}^{r+s} B),
$$
where the constant $C=C(r,s)$ is the Euler beta--function.
}

\smallskip

{\bf Proof.} Directly from the definitions and assumptions, we get
$$
N_{\lambda\mu}(B)= \sum_{i=1}^{N_{\lambda}(B/\mu_i)} N_{\mu}(B/\lambda_i)=
c_{\mu}  \sum_{i=1}^{N_{\lambda}(B/\mu_i)} B/\lambda_i \,
\roman{log}^r (B/\lambda_i)
$$
$$
+ O(\sum_{i=1}^{N_{\lambda}(B/\mu_i)} B/\lambda_i \,
\roman{log}^{r-1} (B/\lambda_i))
\eqno(4.1)
$$
Since the error term has the same structure as the main one,
with $r$ replaced by $r-1$, we can apply this formula inductively,
and get the following expression for the main term
$$
c_{\lambda} c_{\mu} B \sum_{j=1}^N a(j) B/j \, \roman{log}^r (B/j)
\eqno(4.2)
$$
where
$$
a(j) := \roman{card} \, \{i \,|\, \lambda_1+j  \le  \lambda_i < \lambda_1+j+1\} ,\
N:= [B/\mu_1 - \lambda_1]  +1 .
$$
So the main term of (4.1) can be rewritten as
$$
c_{\lambda} c_{\mu} B \sum_{j=1}^N j\,\roman{log}^r j \int_j^{j+1} x^{-2} 
\roman{log}^s (B/x) dx .
$$
After approximating the sum by the integral, we get the expected result. $\blacksquare$

\medskip
{\bf 4.2.3. Hardy--Littlewood method and complete intersections.} Below we will sketch,
following [Ig78], 
methodology and results of applications of the Hardy--Littlewood method
in the setup of Fano complete intersections in projective spaces, 
explained in [FrMaTsch89]. Besides showing the  ``typical'' asymptotic
behaviour for some of them, they can serve as an example of study
of distribution of rational points in {\it adelic spaces}. Later we will
survey some recent developments in this direction, cf. Section 5 below.

\smallskip

Let $A_K$ be the ring of ad\`eles of a number field $K$ and let $X$ be an
affine variety over $K$. The ad\`elic Hardy--Littlewood formula is the
identity
$$ \sum_{x\in X(K)} \varphi(x) =\int_{X(A_K)} \varphi(x) \omega(x) + R(\varphi), $$
where $\omega$ is a measure on $X(A_K)$ and $\varphi$ is in a class of test functions 
on $X(A_K)$. The integral $P(\varphi)=\int_{X(A_K)} \varphi(x) \omega(x$ is the {\it singular series} 
and $R$ the error term, see Section~1.3 of [FrMaTsch89].

\smallskip

Let $\bold{P}^{n}$ be a projective space with a fixed 
system of homogeneous coordinates $(x_0:x_1:\dots :x_n)$ as in 2.1 above,
$K$ a base field of finite degree over $\bold{Q}$.
Consider $m\ge 1$ forms $f_i\in K[x_0,\dots ,x_n]$ and put
$d_i := \roman{deg}\,f_i$. Denote by $V$ the projective
variety over $K$ which is the nonsingular complete intersection of
hypersurfaces $f_i=0$. Then its anticanonical line bundle is
$$
-K_V=\omega_V^{-1} = \Cal{O}_V (n+1 - d_0 - \dots - d_m).
$$
Therefore $V$ is a Fano variety iff $n+1 > d_0+\dots +d_m.$

\smallskip

Let $W$ be a cone over the projective complete intersection $V$, 
For some $\tau>0$ and $\tilde\tau_v:=\tau^{[K_v:\Q_v]/[K:\Q]}$, and for a test
function $\varphi$ on $W(A_K)$, set $\varphi_\tau(x)=\varphi(\tilde\tau^{-1}x)$. If $P(\varphi)$
converges then
$$ P(\varphi_\tau)=\tau^{n+1-\sum_i d_i} P(\varphi_1), $$
as shown in Proposition~4 of [FrMaTsch89]. If $\varphi_1$ is the restriction to $W(A_K)$ of
the characteristic function of a product of balls at the infinite places and the $v$-adic integers 
at finite places, then the sum $\sum_{x\in W_0(K)} \varphi_1(x)$, with $W_0=W\smallsetminus \{ 0 \}$,
provides a counting of $K$-points of $W_0$, which in the case $K=\Q$ agrees with the counting
of rational points of $V$ with $\cO(1)$-height bounded above by $\tau$.

\smallskip

Finally, how summation over points may be approximated by integration
over adelic spaces, as in the case of Proposition~4 of [FrMaTsch89] recalled above,
is described in Section~5 below.

\medskip

{\bf 4.2.4. Generalised flag manifolds.} The last class of Fano varieties
we considered in [FrMaTsch89] consists of generalised flag manifolds
$V=P\backslash G$, where $G$ is a semisimple linear algebraic group, and $P$
is a parabolic subgroup, both defined over $K$. For convenience,
we assumed moreover that $P$ contains a fixed minimal parabolic
subgroup $P_0$. Cf. also [BlBrDeGa90].

\smallskip

Denote by $\pi : G\to V$ the canonical projection.
Let $\Cal{X}^*(P)$ be the group of characters  of $P$ defined over $K$.
Each character $\chi \in \Cal{X}^*(P)$ defines a line bundle on $V$
that we will denote $L_{\chi}$. Its local sections come from those
local functions on $G$ upon which the left multiplication
by $p\in P$ lifts to the multiplication by $\chi (p)$.

\smallskip

The anticanonical line bundle is located among the $L_{\chi}$. It is
denoted $L_{-2\rho}$ in Section~2 of [FrMaTsch89].

\smallskip

Now comes the main part of the construction: explicit
description of the anticanonical height. Let $A_K$ be the adele ring of $K$.
There exists a maximal compact subgroup $\Cal{K} = \prod_v \Cal{K}_v \subset A_K$
such that $G(A_K)  = P_0(A_K)\Cal{K}$, where $P_0$ is the fixed minimal
parabolic subgroup. Over an open set $W$ we write sections as
$$ \Gamma(W, L_\chi)=\left\{ f\in \Gamma(\pi^{-1}(W), \Cal{O}_G)\, \biggr| \, f(pg)=\chi(p) f(g), \,\,
\forall p\in P, \, g\in G \right\}. $$
If $s$ is a local section of $L_{\chi}$ on a neighborhood $W$ of a point $x\in V(K_v)$, 
we refer to $f$ as above as the corresponding local function.
Each bundle $L_{\chi}\otimes F_v$ on $V\otimes K_v$
is endowed by a canonical $K_v$--invariant
$v$--adic norm: if $s$ is a local section of  $L_{\chi}$ with local function $f$, 
we choose $k\in K_v$ with $\pi(k)=x$ and put $|s|_v = |f(k)|_v$. We define
the height function by
$$
h_{\chi} = h_{L_{\chi}} := \prod_v |s|^{-1}_v .
$$
Using this description, we can identify the anticanonical height zeta--function of $V$ with one of the
Eisenstein series from [La76], Appendix II. The analytic properties of these Eisenstein series are
then used to establish the asymptotic formula for the height zeta function.
We omit further details and we refer the reader to Section~2, Theorem~5 and subsequent
corollary in [FrMaTsch89].

\medskip
{\bf 4.3.  Heights with respect to more general line bundles $L_U$ on Fano varieties}.
In [BaMa90] it was suggested that only a slight generalisation of formula in 4.2
should be ``typical" (in the sense of being the expected behavior when the previous
formula fails, although in fact it is valid in a much more restricted set of cases):
$$
N(U, L_U, B) \sim c B^{\beta} (\roman{log}\,B)^t ,\quad t:= \roman{rk}\, \roman{Pic}\, U -1 .
$$
As was argued in [BaMa90], $\beta$ should be defined by the relative positions of $L_U$ and $-K_U$ in the 
cone of pseudo--effective divisors of $U$: see precise conjectures there.

Sh.~Tanimoto in [Ta19] provided arguments proving various inequalities for these numbers
related to the conjectures in [BaMa90].

Finally, the subtlest information about such asymptotics is given by several conjectures and proofs
regarding exact value of the constant $c$.

\medskip

{\bf 4.4. From asymptotic formulas to sieves}. If we restrict ourselves by those $V$ for which we can define a 
{\it Grothendieck
topology, whose objects satisfy strong asymptotic formulas discussed above}, or their
weaker versions,  as stated at the beginning of Section~4.2 and Section~4.3, then
we can try to define sieves in it by some inequalities weakening that in Section~3.6, such as
conditions given directly in terms of the number of points of bounded height, such as
$$
N(V, L_V, B)/ N(U, L_U, B) = o(1)
$$
or in terms of asymptotics, such as, for $\beta (U,L_U) =\beta (V,L_V)$, 
$$
t (V,L_V) < t(U,L_U),
$$
where $t(V,L_V)$ and $\beta (V,L_V)$  refer to the asymptotic formula mentioned at the
beginning of Section~4.3.

\smallskip

Along the lines of Section~3.6, we can then construct other assemblers $\cC_U$
with objects given by locally closed subvarieties $V$ and morphisms given by
inclusions $W\hookrightarrow V$, where the condition $0<\sigma(W,L_W)< \sigma(V,L_V)$
is now replaced by one of the conditions listed here above. 

\smallskip

Moreover, we can use conditions on the number of points of bounded height, such as
 those describing accumulating subvarieties in Example~2.3.1, to obtain an assembler
 whose Grothendieck group detects these subvarieties, as follows.
 
 \medskip
 
 {\bf 4.4.1. Proposition.} {\it 
 Let $\cC_U^s$ and $\cC_U^w$, be the categories where objects are, in both cases,
 the locally closed subvarieties $V\subset U$ and morphisms are finite compositions
 of open embeddings and closed embeddings, where, in the case of $\cC_U^s$
 the closed embeddings $i_{W,V}: W\hookrightarrow V$ satisfy
 $$ \frac{N(W,L_W,B)}{N(V,L_V,B)} \longrightarrow 1,\ \  \text{ for } B\to \infty, $$
 while in the case of $\cC_U^w$ the closed embeddings $i_{W,V}: W\hookrightarrow V$ satisfy
$$ \roman{liminf}_B\,\, \frac{N(W,L_W,B)}{N(V,L_V,B)} >0. $$
 The categories $\cC_U^s$ and $\cC_U^w$ are assemblers, and the associated
 Grothendieck groups are generated by classes $[V]$ of subvarieties of $U$ with
 the relations $[U]=[V]+[U\smallsetminus V]$ where
 $V$ is a strongly or weakly accumulating subvariety, in the case of $\cC_U^s$ or $\cC_U^w$,
 respectively.
 }
 
 \medskip
 
{\it Proof.} The argument is exactly the same as in Proposition~3.6.2. We just need to verify,
as in Lemma~3.5.1, that the composition $W'\subset W\subset V$ of closed embeddings that satisfy one of
the two conditions above still satisfies the same condition. This is clear both for the
first case where one still has convergence to $1$ of the product
$$ \frac{N(W',L_{W'},B)}{N(W,L_W,B)} \cdot  \frac{N(W,L_W,B)}{N(V,L_V,B)} $$
and for the second case where, since all the terms are non-negative, one has
$$ \roman{liminf_B} \frac{N(W',L_{W'},B)}{N(V,L_V,B)}\geq
\roman{liminf}_B \frac{N(W',L_{W'},B)}{N(W,L_W,B)} \cdot 
\roman{liminf}_B \frac{N(W,L_W,B)}{N(V,L_V,B)} > 0. $$
$\blacksquare$

\bigskip
\newpage

 \centerline{\bf 5. Sieves ``beyond heights'' ?}
 \medskip

{\bf 5.1.  Thin sets and Tamagawa measures}. Below we survey recent attempts to define the geome\-try of subsets
of rational points of $V(K)$, containing ``considerably less'' points than $V$.

From our viewpoint, these definitions should
also be tested on compatibility with our  philosophy of sieves and
assemblers.

\smallskip

Below, we adopt the framework of [Pe18], which was further studied by W.~Sawin ([Sa20],
under the additional restriction $K=\bold{Q}$). His paper starts with
Conjecture~1.1, called ``Modern formulation of Manin's conjecture''.
It involves both shifts from and extensions of the setup in our
previous Section~4.

\smallskip

We first need to define the following objects/quantities, in order 
to state this modern formulation:

\medskip

(i) {\it Summation} over points $x$ of  height $\le B$ is replaced by
the {\it averaging} of  the measures $\delta_x$ of the same points
embedded into the adelic space $V (A_{\bold{Q}})$.

\smallskip
(ii) The definition of the respective Tamagawa measure $\tau$ assumes
that $V$ is a geometrically integral smooth projective Fano variety,
with Picard group of rank $r$. Let $\Cal{V}$ be its proper integral
model over $\bold{Z}$. The $L$-function $L(s, Pic\, V_{\overline{\bold{Q}} })$ 
defined by the Picard group
and the respective local zetas $L_v$ are zetas of lattices
(in the sense that they are zeta functions associated to Artin (Galois) 
representations on lattices), see Construction~3.28 in [Pe18].
\smallskip

Then
$$
\tau := (\roman{lim}_{s\to 1} (s-1)^r L(s, Pic\, V_{\overline{\bold{Q}} }))
\prod_v L_v(s, Pic\, V_{\overline{\bold{Q}} })^{-1} \omega_v .
$$
Here $\omega_v$ is defined by the natural measure on local
non--archimedean points of $V$ or archimedean volume form.
\medskip
(iii) Finally, define the numbers $\alpha (V)$ and $\beta (V)$ by

$$
\alpha (V) := r\,\roman{vol}\, \{ y\in (({ \roman{Pic} (V)\otimes \bold{R})^{eff})^{\vee}} \,|\, 
K_V\cdot y \le 1\},
$$
and
$$
\beta (V) := \roman{card}\, H^1 (\roman{Gal} (\overline{\bold{Q}/\bold{Q}}),
\roman{Pic}\, V_{\overline{\bold{Q}}}) .
$$
 Then the modern formulation of the conjecture on the number of rational
 points of bounded height on a Fano variety, according to [Sa20]
 can be stated as follows.
 \medskip
 Let $f:\, U\to V$ be a morphism of geometrically integral smooth projective varieties.
 Call it  a {\it thin morphism}, if the induced map $U\to f(U)$ is generically finite
 of degree $\ne 1$.
 
 \smallskip
 
 The modern formulation of the Manin conjecture can then be stated
 in the following way (Conjecture~1.1 of [Sa20]).
 
 \medskip

 {\bf  5.2. Conjecture.} {\it There exists a finite set of thin maps $f_i: Y_i \to X$,
 with $W$ the complement of the union of the sets $f_i(Y_i(\Q))$,
 such that we have an exact formula for the weak limit of the form
 $$
 \lim_{B\to \infty}\, \frac{1}{B (\roman{log}B)^{r-1}}
 \sum \Sb x\in W(\bold{Q}) \\ H(x) < B \endSb \delta_x =
 \alpha (V) \beta (V) \tau^{Br} ,
 $$
where $\tau^{Br}$ is the restriction of the Tamagawa measure on the subset of
$V(A_{\bold{Q}})$ on which the Brauer--Manin obstruction vanishes.}

\medskip

We will review the Brauer--Manin obstruction, along with other
related material, in Section~6.1, where we also refine the general
question proposed in this section.

\smallskip
 
 For a class of varieties for which this conjecture is valid,
 we obtain interesting new possibilities for defining sieves.
 Thus, we can formulate this as a more general open ended question.
 
 \medskip

 {\bf 5.3. Question.} {\it
 Is there a natural construction of sieves and assemblers 
 based on subvarieties satisfying the conditions mentioned
 above? Such assemblers would be categorifications
 (with associated spectrifications) of these possible ways of
 expressing the notion of  having asymptotically 
 ``fewer rational points than the ambient variety", without
 relying directly on the height zeta function. What information
 can be obtained from the homotopy invariants of the spectra
 of such assemblers?
 }

 \medskip
 Regarding Brauer--Manin obstruction itself, see the monograph [CThSk19]
 and many references therein, in particular,  [Sko95], [Sko09].   This obstruction cannot explain many
 cases of missing rational points even when the initial obstruction for them vanishes,
 and it was generalised and/or modified in different ways many times: see in particular 
 [LeSeTa18], [Ta19], [DePi19]. It is worth mentioning, in particular, 
 Skorobogatov's \'etale Brauer--Manin obstruction, [CoPaSk16], [Sko99], [Sko09].

 \medskip
 The recent article [CorSch20] introduces interesting contexts
 that are very suitable for considering generalised obstructions
 from the viewpoint of sieves (see also [LiXu15] and many other references
 in these papers).
 We dedicate the next Section 6  to
 a description of these contexts.
 
 \medskip
 
 Finally, we should mention the related but not identical machinery
 for studying varieties with finite but non--empty sets of $K$--points.
 For this, we refer the reader to [Qu15].
 
\bigskip

\bigskip

 \centerline{\bf 6. Obstructions and sieves}

\medskip

{\bf 6.1. ``Invisible'' varieties.}  Let $V$ and $U$ be two objects
of a category with unique initial object $\emptyset$. Then we will
say that $V$ is {\it ``$U$--invisible''}, if  $\roman{Hom} (U,V) =\emptyset$.

\smallskip

The main setup we will have in mind is the category of $K$--schemes,
where $K$ is a number field.  
The $\roman{Spec} K$-invisible varieties over $K$ are exactly
those with no $K$-points. 
Then if $U = \roman{Spec}\, L$,
where $L$ is a commutative $K$--algebra,  in place
of ``$U$--invisibility'' we will sometimes speak about ``$L$--invisibility''.
Hopefully, this will not lead to a confusion. In particular,
we will often speak about $K$--invisible varieties over $K$,
and also $A_K$--invisible ones, where $A_K$ is the ring
of ad\`eles.

\smallskip

$K$--invisible varieties  are the objects of a category $\Cal{X}_K$,
whose morphisms are all morphisms of $K$--varieties (or a subfamily 
of them closed under composition.)

\smallskip

{\it The local--global principle} for a class of $K$--varieties $X$ is the statement that
 $X(A_K) \ne \emptyset$ implies $X(K) \ne \emptyset$, where $A_K$ is
 the ring of ad\`eles of $K$. Thus, such a variety can be $K$--invisible only if
 it has no ad\`elic points.  (Formerly this was called  {\it the Hasse principle}).
 
 \smallskip

 The {\it Brauer--Manin obstruction} to the local--global principle is
 provided by the definition of a set  $X(A_K)^{\roman{Br}}$
 such that $X(K) \subseteq X(A_K)^{\roman{Br}} \subseteq X(A_K).$ 
 One says that failure of the local--global principle for $X$ {\it is explained by this
 obstruction,} if $X(A_K) \ne \emptyset$, but  $X(A_K)^{\roman{Br}} = \emptyset$.
\smallskip

The first remarkable new result of [CorSch20] is the following theorem (Corollary 1.1 on p. 6 of [CorSch20]):

\medskip

{\bf 6.1.1. Theorem.} {\it Assume that $K$ is either totally real, or imaginary quadratic field.
Then for any $K$--invisible $K$--variety $V$ there exists its Zariski open covering $X=\cup_i U_i$
such that invisibility of each $U_i$ is explained by the Brauer--Manin obstruction.}
\medskip

The notions of ``sieves'', ``assemblers'' etc. do not appear explicitly in [CorSch20],
and there are no references to the papers of I. Zakharevich there.  However, Theorem 6.1.1.
leaves no doubt  that the setup of sieves is essential for the better understanding of  their work.
\smallskip

Some details and explanations follow.

\medskip

{\bf 6.2. Functor Obstructions.} (Sec. 8.1 of [Po17]). 
Corwin and Schlank work with definitions of obstructions,
that were systematically developed in [Po17].  Denote by $Sch_K$, resp. $Var_K$, the category
of $K$--schemes, resp. algebraic $K$--varieties.
 Let $F$ be a contravariant functor $Sch_K \to Sets$
 (presheaf of sets on the category of $K$--schemes).
\smallskip

If $V=\roman{Spec}\, L$ where $L$ is a commutative $K$--algebra, we
will write $V(L)$, resp. $F(L)$, in place of $V(\roman{Spec}\,L)$, resp. 
$F(\roman{Spec}\,L)$, so that with this notation, $F$ becomes a {\it covariant}
functor of its argument $L$. 
\medskip

{\bf 6.2.1. Definition.}  {\it Let $\omega$ be a subfunctor of the functor
$F: Var_K \to Sets$.

\smallskip

We will say that $\omega$ is a generalised obstruction to the Hasse principle,
if $V(K) \subseteq \omega (V(A_K))$ for all  $K$--varieties $V$.
}
\smallskip

We will often write $V(A_K)^{\omega}$ in place of  $\omega (V(A_K))$,
as in the notation for the Brauer--Manin obstruction above.

\smallskip

From the definitions, one easily sees that if $V$ is a proper variety, then
$$
V(A_K) = V(\prod_v K_v)= \prod_{v} V(K_v).
$$
Generally, we have only inclusion
$$
V(A_K) \subseteq V(\prod_v K_v)= \prod_{v} V(K_v).
$$
\smallskip
Let now  $F$ be a contravariant functor from $Sch_K$ to $Sets$,
and $L$ a $K$--algebra.

\smallskip

In [Po17], pp. 227-- 228, it is shown how to define for each $A\in F(V)$, and
each $K$--variety $V$, the ``evaluation'' map
$$
ev_A : V(L) \to F(L).
$$
Namely, an $L$--point  $x\in V(L)$ considered as morphism 
of $K$--schemes, induces a map $F(V) \to F(L)$. The image
of  $x$ in $F(L)$ is denoted $ev_A(x)$, or simply $A(x)$. 

 \medskip
 
 From the definitions, it follows that we have a commutative diagram
 $$
 \xymatrix{V(K) \ar[r] \ar[d]^{ev_A} & V(A_K) \ar[d]^{ev_A}\\
 F(K)\ar[r]  & \prod_v F(K_v)} 
 $$ 
The upper horizontal arrow is an embedding. Denote by
$V(A_K)^A$, for $A\in F(V)$, the subset of $V(A_K)$ consisting
of elements whose image in $F(A_K)$ lies 
in the image of $F(K)$.
\smallskip
Finally, put
$$
V(A_K)^F := \cap_{A\in F(V)} V(A_K)^A.
$$
Clearly, $V(K)\subseteq V(A_K)^F$.
\smallskip

We can similarly define the respective subsets in 
$V(\prod K_v)) = \prod V(K_v)$ and, using
$\prod F(K_v)$, define $V(\prod K_v)^F$.

\medskip

{\bf 6.2.2. Definition.}  {\it If
$$
V(\prod K_v) \neq \prod V(K_v)^F := \prod F(K_v) ,
$$
we will say that this inequality 
is explained by an $F$--obstruction for weak approximation.
}
\medskip
In the proof of Theorem 6.1.1 above, Corwin and Schlank use also versions
of these definitions involving ad\`eles
with restricted denominators.  We will briefly sketch parts of their
arguments below.

\smallskip

Using notations of Subsection~2.1 above,
choose a subset of places $S\subset \Omega_K$, and denote (as
in Section~1.1 of [CorSch20]) by $A_{K,S}$
the restricted product
$A_{K,S}=\prod_{v\in S}^\prime K_v$ where the components belong 
to $\Cal{O}_v$ for all but finitely many $v\in S$. In particular, one
considers the case where $S=\Omega_{K,f}$ is the set of finite places.

\smallskip

Generally, a triple $(K,S,\omega )$ as above is called {\it an obstruction datum}
in Definition~2.3 of [CorSch20].

\medskip

{\bf 6.2.3. Definition.}  ([CorSch20], Def. 2.4). {\it A $K$--variety
$V$ satisfies VSA (very strong approximation principle) for $(K,S,\omega )$,
if $V(K) = V(A_{K,S})^{\omega}$.
}
\medskip

As in Sec. 4.2 above, we can check   that the family of $K$--varieties
satisfying VSA for $(K,S,\omega )$ (or simply $VSA$--varieties) is closed with respect to:
\medskip

 (i) {\it  direct products; }
 \smallskip
(ii) {\it passage to locally closed subvarieties.}

\medskip

{\bf 6.3. Theorem.}  ([CorSch20], Theorem 4.3.)   {\it Assume that there is an open subvariety
of $\bold{P}^1_K$ satisfying VSA for $(K,S,\omega )$.
\smallskip

Then for any $K$--variety $V$  there exists a finite affine open cover
$V=\cup_i \, U_i$ such that each $U_i$ satisfies VSA for  $(K,S, \omega )$. 
}

\smallskip

{\it Sketch of proof.} Any closed point of $V$ belongs to an open affine subvariety
of $V$. Because of stability of VSA varieties with respect to direct products
and locally closed embeddings, we can find an affine neighbourhood of each
closed point satisfying VSA, and then to choose from all such neighbourhoods
a finite open cover. $\blacksquare$

\bigskip

 \centerline{\bf 7. Assemblers and spectra for Grothendieck rings with exponentials}

\medskip

Starting with this section, we investigate the homotopy theoretic framework underlying another
aspect of the height function and the associated height zeta function, in the framework of the
{\it motivic height zeta function} introduced by Chambert-Loir and Loeser in [ChamLoe15].
In this section we focus on the relevant Grothendieck rings, namely the {\it Grothendieck rings with exponentials}
defined in [ChamLoe15] and their  localizations. We construct associated assembler categories and homotopy
theoretic spectra, in the sense of [Za17a], [Za17b], [Za17c].

\medskip

{\bf 7.1. Localization of the Grothendieck ring.}
Consider the localization $\cM_K$ of the Grothendieck ring of varieties $K_0(\cV_K)$ 
obtained by inverting the Lefschetz motive $\bL=[\bA^1]$ and all $\bL^n-1$ for
$n\geq 1$. This localization can be identified with the Grothendieck ring of algebraic
stacks $K_0(\cS_K)$ (see [BeDh07], [Ek09]) by expressing it as the
localization with respect to the classes $$[ GL_n ]=(\bL^n-1)(\bL^n-\bL)\cdots (\bL^n-\bL^{n-1}).$$

\smallskip

The Grothendieck ring of stacks $K_0(\cS_K)$ is generated as a group by the 
isomorphism classes $[X]$ of algebraic stacks $X$ of finite 
type over the field $K$ with the property that all their automorphism 
group schemes are affine, modulo relations
$$[X] =[Y]+[U]$$ for a closed substack $Y$ and its complement $U$, and 
$$[E]=[X\times \bA^n]$$ for a vector bundle $E\to X$ of rank $n$. The multiplication
 is given, as in the case of varieties, by the product $$[X\times Y]=[X]\cdot [Y].$$

\smallskip

There is a natural map $K_0(\cV_K) \to K_0(\cS_K)$. 
The classes $[ GL_n ]$ are invertible in $K_0(\cS_K)$. This can be seen by
first observing that, if $X \to Y$ is a 
$GL_n$-torsor of algebraic stacks of finite type over $K$, then $[X]=[GL_n]\cdot [Y]$
in $K_0(\cS_K)$. This fact follows from the relation $[E]=\bL^n \cdot [S]$ for rank $n$
vector bundles over $S$, repeatedly applied to $X =X_n \to X_{n-1}\to \cdots \to Y$ 
with $X_k$ the stack of $k$ linearly independent vectors in the vector bundle determined
by $X$, and viewing each $X_k \to X_{k-1}$ as a complement of a vector subbundle in
a vector bundle (see Proposition~1.1 of [Ek09]). The invertibility then follows by
considering the $GL_n$-torsor $Spec(K)\to BGL_n$ that gives $1=[GL_n]\cdot [BGL_n]$
in $K_0(\cS_K)$.  

\smallskip

This implies that the map $K_0(\cV_K) \to K_0(\cS_K)$ induces
a map from the localization $$\cM_K \to K_0(\cS_K).$$ In order to produce
 the inverse map, one assigns to the
stack given by a global quotient $[X/GL_n]$, the class $[X]/[GL_n]$
in $\cM_K$, where the resulting class is independent of the description of initial one as a 
global quotient, because for the general case one can stratify a stack by global quotients
and add the resulting classes (see Theorem~1.2 of [Ek09]). 

\smallskip

Observe also that the second relation $[E]=\bL^n [Y]$ in $K_0(\cS_K)$ for a rank
$n$ vector bundle $E\to Y$ is in fact equivalent to assuming the relation 
$[X]=[GL_n]\cdot [Y]$ for a $GL_n$-torsor $X\to Y$. One direction was explained
above, while in the other direction it suffices to observe that one can realize a rank $n$
vector bundle $E\to Y$ as $E=V\times_{GL_n} X$, where $V$ is an $n$-dimensional vector
space and $X\to Y$ is a $GL_n$-torsor. The quotient map $$V\times X \to V\times_{GL_n} X$$
is also a $GL_n$-torsor, hence we have $$[GL_n] \cdot [E]=[V] \cdot [X]=[V]\cdot [GL_n] \cdot [Y],$$ which gives
$[E]=\bL^n [Y]$ since the class $[GL_n]$ is invertible. 

\smallskip

Notice, moreover, that it suffices to consider the class $\cS^\cZ$ of stacks where the automorphism 
group schemes in the class $\cZ$ are connected finite type group schemes for which torsors over
any finitely generated extension of $K$ are trivial. Indeed, it is shown in Proposition~1.4
of [Ek09] that $K_0(\cS_K)=K_0(\cS^\cZ_K)$, and that for $G$ a group
scheme in $\cZ$ any $G$-torsor $X\to Y$, with $X,Y \in \cV_K$, has a stratification where the
torsor is trivial on each stratum, so that $[X]=[G][Y] \in K_0(\cV_K)$. If $K$ is an
algebraically closed field, the group schemes in $\cZ$ are affine and with the property that
$GL_n \to GL_n/G$ is a Zariski locally trivial fibration ([Ek09], Remark after
Proposition~1.4). 

\smallskip

We can define an assembler category underlying the localization 
$\cM_K=K_0(\cS_K)$ and the associated homotopy theoretic spectrum in the
following way.

\medskip

{\bf 7.2. Proposition.} {\it
Let $\cC^{\cS^\cZ_K}$ be the category whose objects are algebraic stacks $X$ of finite 
type over the field $K$ with automorphism group schemes in $\cZ$, and whose  morphisms
are locally closed embeddings. A Grothendieck topology on $\cC^{\cS^\cZ_K}$ is 
generated by the families 
$$\{ Y \hookrightarrow X, U \hookrightarrow X \}$$ where 
$Y \hookrightarrow X$ a closed substack and $U$ its complement.  The category
$\cC^{\cS^\cZ_K}$ is an assembler, and the associated spectrum $K \cC^{\cS^\cZ_K}$
satisfies the relation
$$
\pi_0(K \cC^{\cS^\cZ_K})=K_0(\cS^\cZ_K).
$$
}

\medskip

{\it Proof.} The proof that this is an assembler category uses similar arguments
as in [Za17a]: the empty set is the initial object; finite disjoint covering families are
given by 
$$
X_i \hookrightarrow X$$ with $$X_i=Y_i\smallsetminus Y_{i-1}
$$ 
for
a chain of embeddings 
$$
\emptyset =Y_0\hookrightarrow Y_i\hookrightarrow \cdots \hookrightarrow Y_n=X.
$$
Any two finite disjoint coverings have a common refinement because the category has pullbacks,
as in [Za17a]. 

\smallskip

The spectrum $K \cC^{\cS^\cZ_K}$ determined by the assembler $\cC^{\cS^\cZ_K}$ has 
$\pi_0(K \cC^{\cS^\cZ_K})$ generated by the objects of $\cC^{\cS^\cZ_K}$ with the
scissor-congruence relations determined by the disjoint covering families. These include
the relations $[X]=[Y]+[U]$ for closed substacks $Y\hookrightarrow X$ and their 
complements $U\hookrightarrow X$, as well as the relations $[X]=[GL_n] [Y]$ for 
$GL_n$-torsors $X\to Y$ through the existence of a disjoint covering family determined
by a stratification where the torsor is trivial on each stratum as discussed above. Thus
we obtain that $\pi_0(K \cC^{\cS^\cZ_K})=K_0(\cS^\cZ_K)$.
$\blacksquare$

\medskip

{\bf 7.3. Grothendieck rings with exponentials.}
The {\it Grothendieck ring with exponentials} of [ChamLoe15] is designed to
describe a motivic version of exponential sums of the form
$$
 \sum_{x\in X(\F_q)} \chi(f(x)), 
 $$
where $X$ a variety over a finite field, endowed with a function $f: X\to \bA^1$ and
a fixed non-trivial character $\chi: \F_q \to \C^*$. 

\smallskip

In fact, we will see in Section~7.8 that
this heuristics can be made precise, in the sense that
assigning such an exponential sum to a class in the Grothendieck ring with exponentials
defines a ring homomorphism.

\smallskip

In the formulation of  [ChamLoe15], the Grothendieck ring with exponentials
$KExp_K$ is defined as follows.

\smallskip

{\bf 7.3.1. Definition.} {\it The Grothendieck ring with exponentials $KExp_K$
is generated by isomorphism classes of pairs $(X,f)$ of a $K$-variety $X$ 
and a morphism $f: X\to \bA^1$, where two such pairs $(X_1,f_1)$ and $(X_2,f_2)$ 
are isomorphic if there is an isomorphism $u: X_1\to X_2$ of $K$-varieties such
that $f_1=f_2\circ u$. The relations are given by
$$
 [X,f]=[Y,f|_Y] + [U,f|_U], 
 $$
for a closed subvariety $Y\hookrightarrow X$ and its open complement $U$, and
$$
 [X\times \bA^1, \pi_{\bA^1}]=0 
 $$
where $\pi_{\bA^1}: X\times \bA^1 \to \bA^1$ is the projection on the second factor.
The ring structure is given by the product
$$
 [X_1,f_1]\cdot [X_2,f_2]= [X_1\times X_2, f_1 \circ \pi_{X_1} + f_2 \circ \pi_{X_2}] 
 $$
where $f_1 \circ \pi_{X_1} + f_2 \circ \pi_{X_2}: (x_1,x_2)\mapsto f_1(x_1)+f_2(x_2)$.}

\smallskip

The relation $[X\times \bA^1, \pi_{\bA^1}]=0$ corresponds to the fact that, for
an exponential sum as above, one has $\sum_{a\in \F_q} \chi(a)=0$.

\smallskip

There is an embedding of the ordinary Grothendieck ring of varieties $K_0(\cV_K)$ in
the Grothendieck ring with exponentials $KExp_K$ that maps $[X]\mapsto [X,0]$.
This map is compatible with localizations and it induces an embedding of the
localization $\cM_K$ in the corresponding localization $Exp\cM_K$ of $KExp_K$.

\smallskip

Given the identification $\cM_K \simeq K_0(\cS_K)$ of the localization $\cM_K$
of the Grothendieck ring of varieties and the Grothendieck ring of algebraic stacks,
one can similarly obtain a description of the corresponding localization 
$Exp\cM_K$ of $KExp_K$ of the Grothendieck ring of varieties with exponentials
in terms of a {\it Grothendieck ring of algebraic stacks with exponentials}. 

\smallskip

{\bf 7.3.2. Definition.} {\it The Grothendieck ring of algebraic stacks with exponentials
$K_0(Exp\cS_K)$ is generated by the isomoprhism classes $[X,f]$ of pairs $(X,f)$ 
of an algebraic stack $X$ and a morphism $f: X \to \bA^1$, where an isomorphism
$(X_1,f_1)\simeq (X_2,f_2)$ is given by an isomorphism $u: X_1 \to X_2$
with $f_2\circ u=f_1$.  The relations in $K_0(Exp\cS_K)$ are given by 
$[X,f]=[Y,f|_{Y}]+[U,f|_{U}]$ for a closed substack $Y$ and its open complement; 
$[E,f]=\bL^n \cdot [S,f_S]$ for $\pi: E\to S$ a vector bundle and $f: E\to \bA^1$
with $f_S\circ \pi=f$, with $\bL^n=[\bA^n,0]$; and with the further relation
$[X\times \bA^1,\pi_{\bA^1}]=0$. The ring structure is defined as in $KExp_K$.}

\smallskip

The argument of [Ek09] recalled above, for the identification of $K_0(\cS_K)$ with
the localization $\cM_K$ of $K_0(\cV_K)$ extends to the case with exponentials.

\medskip

{\bf 7.3.3. Lemma.} {\it The natural morphism $KExp_K \to K_0(Exp\cS_K)$
factors through the localization $Exp\cM_K$ and induces an isomorphism
$Exp\cM_K \simeq K_0(Exp\cS_K)$.}

\smallskip

{\it Proof.} As in the original argument of [Ek09], the vector bundle relation
implies that, if $X\to Y$ is a $GL_n$-torsor, endowed with compatible maps
$f_X,f_Y$ to $\bA^1$, then $[X,f_X]=[GL_n]\cdot [Y,f_Y]$. The invertibility
of $[GL_n]=[GL_n,0]$ follows from its invertibility in $K_0(\cS_K)$. This
implies that the homomorphism $KExp_K \to K_0(Exp\cS_K)$ factors through
a homomorphism $Exp\cM_K \to K_0(Exp\cS_K)$. To construct the inverse map
we assign to a class $[X/GL_n, f]$, where $X/GL_n$ is a global quotient,
the class $[GL_n]^{-1}\cdot [X,\tilde f]$ where $\tilde f=f\circ \pi$ for the
quotient map $\pi: X \to X/GL_n$. This map extends 
an arbitrary $[X,f]$ by stratifying the stack $X$ by global quotients
with the restriction of the morphism $f$ to the strata.
$\blacksquare$

\medskip

We can then consider assembler and homotopy theoretic spectra associated to
these Grothendieck rings. This can be done using the setting of [Za17a]
for subassemblers and cofiber sequences. 

\medskip

{\bf 7.4. Theorem.} {\it 
A simplicial assembler $\cC^{KExp}_K$ with $\pi_0 K(\cC^{KExp}_K)=KExp_K$, the 
Grothendieck ring with exponentials, is obtained as the cofiber of the 
morphism of assemblers $\Phi: \cC \to \cC$, where the objects of $\cC$ 
are  pairs $(X,f)$ of a $K$-variety $X$ 
and a morphism $f: X\to \bA^1$, and morphisms are locally closed embeddings
of subvarieties with compatible morphisms to $\bA^1$. Finally,
 $$
 \Phi(X,f) :=
(X\times \bA^1, f\circ \pi_X + \pi_{\bA^1}).
$$
}

\medskip

{\it Proof.} Endow the category $\cC$ with the Grothendieck topology 
generated by the families 
$$
\{ (Y,f|_Y) \hookrightarrow (X,f), (U,f|_U) \hookrightarrow (X,f) \}
$$ 
where $Y \hookrightarrow X$ is a closed subvariety, and $U$ its open complement. The category $\cC$
is an assembler, since  the category has 
pullbacks, hence finite disjoint families have a common refinement, 
and morphisms are compositions of embeddings and therefore monomorphisms. 
The associated spectrum $K(\cC)$ has $\pi_0 K(\cC)$ generated by the isomorphism
classes $[X,f]$ with relations 
$$
[X,f]=[Y,f|_Y]+[U,f|_U].
$$

\smallskip

Consider then the endofunctor $\Phi: \cC \to \cC$ that assigns to an object
$(X,f)$ the object $$\Phi(X,f)=(X\times \bA^1, f\circ \pi_X + \pi_{\bA^1})$$ for 
$\pi_X,\pi_{\bA^1}$ the projections onto the two factors of $X\times \bA^1$.
This functor is a morphism of assemblers and 
the induced map on $\pi_0 K(\cC)$ is given by multiplication by
the class $[\bA^1, id]$. Note that this is not the multiplication by the
Lefschetz motive, as in this setting $\bL=[\bA^1,0]$. We introduce the
notation $\bY:=[\bA^1,id]$ and we write the map on $\pi_0 K(\cC)$
as $\cdot \bY$.

\smallskip

By the localization theorem (Theorem C of [Za17a]), the morphism $\Phi:  \cC \to \cC$
of assemblers (seen as a morphism between constant simplicial assemblers) has an associated
simplicial assembler $\cC/\Phi$ with a morphism $\iota: \cC \to \cC/\Phi$ that gives a cofiber sequence
$$ 
\xymatrix{ K(\cC) \ar[r]^{K(\Phi)} & K(\cC) \ar[r]^{K(\iota)} & K(\cC/\Phi) }. 
$$
There is an associated exact sequence
$$
 \xymatrix{ \pi_1 K(\cC) \ar[r] & \pi_1 K(\cC/\Phi) \ar[r] & \pi_0 K(\cC) \ar[r]^{\cdot \bY} & \pi_0 K(\cC)  \ar[r] & \pi_0 K(\cC/\Phi)}. 
 $$
Thus, the cokernel $\pi_0 K(\cC/\Phi)$ of the map 
$$
\pi_0 K(\cC) \to \pi_0 K(\cC)
$$
mapping 
$$
 [X,f] \mapsto [X,f]\cdot \bY
 $$
 can be identified with the Grothendieck ring
generated by the classes $[X,f]$ and the relations 
$$
[X,f]=[Y,f|_Y]+[U,f|_U]
$$ 
and
$$
[X\times \bA^1, \pi_{\bA^1}]=0.
$$
Note that $\pi_0 K(\cC/\Phi)$ has in fact all relations of the form 
$[X\times \bA^1,f\circ\pi_X+\pi_{\bA^1}]$, but these are all zero 
because $[\bA^1,\pi_{\bA^1}]=0$ and the product is well defined.
$\blacksquare$

\medskip

In a similar way, we can treat  the localization $Exp\cM_K$ of the
Grothendieck ring with exponentials $KExp_K$, by describing it in terms of
algebraic stacks with exponentials as in Lemma~7.3.3. 

\smallskip

The following statement can be proved 
in the same way, as
 Theorem 7.4. 

\medskip

{\bf 7.5. Proposition.} {\it 
Let $\cC^{KExp\cS}_K$ denote the cofiber of the morphism of assemblers
$$\Phi: \cC^\cS \to \cC^\cS$$ where $\cC^\cS$ has objects given by pairs $(X,f)$
of an algebraic stack $X$ and a morphism $f: X \to \bA^1$, and morphisms
given by locally closed embeddings of substacks with compatible maps to $\bA^1$, and
$$\Phi(X,f)=(X\times \bA^1, f\circ \pi_X + \pi_{\bA^1}).$$
The localization $Exp\cM_K$ of the
Grothendieck ring with exponentials $KExp_K$ is identified with
$$Exp\cM_K=\pi_0 K(\cC^{KExp\cS}_K).$$
}

\medskip

One can also consider the relative version of the Grothendieck ring with
exponentials and its localization, as in Section~1.1.5 of [ChamLoe15]. 

\smallskip

{\bf 7.5.1. Definition.} {\it 
Let $S$ be a $K$-variety. Consider pairs $(X,f)$ where $X$ is an $S$-variety (a
 variety endowed with a morphism $u:X\to S$)
and $f: X \to \bA^1$ is a morphism. The relative Grothendieck ring with
exponentials $KExp_S$ is defined as in Definition~7.3.1, generated by isomorphism classes
$[X,f]_S$ with the relations $$[X,f]_S=[Y,f|_Y]_S+[U,f|_U]_S$$ for $Y\hookrightarrow X$
a closed embedding of $S$-varieties (with morphism $u|_Y: Y \to S$) and 
$$[X,f]_S \cdot [\bA^1, \pi_{\bA^1}]_S=0.$$ }

\smallskip

A morphism $\varphi: S\to T$ induces a
ring homomorphism $$\varphi^*: KExp_T \to KExp_S $$ given by $$\varphi^*[X,f]_T=
[X\times_T S, f\circ \pi_X]_S.$$ 
There is an embedding $$K_0(\cV_S)\to KExp_S$$
given by $$[X]_S \mapsto [X,0]_S,$$ with the localizations satisfying
$$ \cM_S = \cM_K \otimes_{K_0(\cV_K)} K_0(\cV_S), \ \ \ \ \ Exp\cM_S=Exp\cM_K \otimes_{KExp_K} KExp_S. $$

\smallskip

If we interpret   Grothendieck rings with exponentials as
abstract motivic functions, the relative case corresponds to a motivic
generalization of functions that for finite base fields take the form
$$ 
\Psi: S(\F_q) \to \C, \ \ \ \  \Psi(s)=\sum_{x\in X_s(\F_q)} \chi(f(x)) .
$$
Here $X$ is a variety over $S$ with fibers $X_s$ over $s\in S$, and $\chi$ is a fixed character $\chi: \F_q \to \C^*$.

\smallskip

We then have the following statements whose proof is analogous to the previous theorem.

\medskip

{\bf 7.6. Proposition.} {\it 
Let $\cC^{KExp}_S$ denote the cofiber of the morphism of assemblers
$$
\Phi: \cC_S \to \cC_S .
$$ 
Here objects of $\cC_S$ are
triples $(X,u,f)$, where $X$ an $S$--variety, $u: X\to S$ its structure morphism, and
 $f: X \to \bA^1$. Morphisms between such objects are locally closed embeddings of
$S$--subvarieties compatible with the maps to $\bA^1$. Finally,
$\Phi (X,u,f)$ is the $S$--variety obtained by taking a product with $\bA^1$,
endowed with the morphism 
$$
f\circ \pi_X^* +\pi_{\bA^1}
$$ 
to $\bA^1$. Then
$$
\pi_0 K(\cC^{KExp}_S)=KExp_S.
$$

}

\medskip

{\bf 7.7. Proposition.} {\it 
The case of the localization $Exp\cM_S$ can be treated similarly in terms of an assembler
$\cC^\cS_S$ with objects given by algebraic stacks with morphisms to $S$
and to $\bA^1$ and the cofiber $\cC^{KExp\cS}_S$ of the morphism
$$\Phi: \cC^\cS_S \to \cC^\cS_S$$ that multiplies a stack $X$ by $\bA^1$ with the 
morphism $$f\circ \pi_X^* +\pi_{\bA^1}$$ to $\bA^1$, so that
$$\pi_0 K(\cC^{KExp\cS}_S)=Exp\cM_S.$$
}

\medskip

{\bf 7.8. Grothendieck ring with exponentials and zeta functions.}
The Hasse-Weil zeta function of a variety $X$ over a finite field $\F_q$,
$$ \zeta^{HW}_X(t) =\exp \left( \sum_{m\geq 1} \frac{\# X(\F_{q^m})}{m} t^m \right) $$
has an Euler product expression of the form of a product over closed points of $X$,
$$ \zeta^{HW}_X(t) = \prod_{x\in X} (1-t^{\deg(x)})^{-1}. $$
One can see this by observing that the latter expression can be written in the form
$$ \prod_{x\in X} (1-t^{\deg(x)})^{-1}= \prod_{r\geq 1} (1-t^r)^{-a_r}, $$
where 
$$ a_r:= \# \{ x\in X\,|\, [k(x): \F_q] =r \} $$
with $k(x)$ the residue field at the point $x$ with $deg(x)=[k(x): \F_q]$, and that the counting
$$ N_m(X):= \# X(\F_{q^m}) $$
is given by
$$ N_m(X) = \sum_{r|m} r \cdot a_r . $$

\smallskip

One can also reformulate the expression of the Hasse-Weil zeta function as a product over points
in the form of sum over effective zero-cycles of a given degree, by expanding the product, 
$$ \zeta^{HW}_X(t) =\sum_{n\geq 0} \# \{ \text{effective 0-cycles of degree n on X} \} t^n. $$
The latter expression motivated the introduction of Kapranov's motivic zeta function
[Kap00],
$$ Z_X(t)=\sum_{n\geq 0} [Sym^n X]\, t^n, $$
where $[Sym^n X]\in K_0(\cV)$ are the classes in the Grothendieck ring of
varieties of the symmetric products $Sym^n X=X^n /\Sigma_n$ of $X$. The symmetric
product $Sym^n X$ indeed parameterizes the effective zero-cycles of degree $n$ on $X$.

\smallskip

In general, by a motivic measure we mean here a ring homomorphism $$\mu: K_0(\cV) \to R$$
from the Grothendieck ring of varieties to a commutative ring $R$. 

\smallskip

As above, consider Kapranov's motivic zeta function
$$ Z_X(t)=\sum_{n\geq 0} [Sym^n X]\, t^n, $$
with $Sym^n X$ the symmetric products. For a motivic measure $\mu: K_0(\cV) \to R$
one obtains an associated zeta function by taking
$$
 \zeta_\mu (X, t) :=\sum_{n=0}^\infty \mu(Sym^n X)\, t^n .
$$

\smallskip

Consider the case of varieties over a finite field $K=\F_q$, and the Grothendieck
ring with exponentials $KExp_K$. The choice of a character $\chi: \F_q \to \C^*$
determines a motivic measure, namely a ring homomorphicm $\mu_\chi : KExp_K \to \C$ by 
$$ \mu_\chi [X,f]=\sum_{x\in X(\F_q)} \chi(f(x)). $$
This is indeed a ring homomorphism as $$\mu_\chi ([Y,f|_Y]+[U,f|_U])=\mu_\chi [Y,f|Y]+\mu_\chi [U,f|_U]$$ 
and $$\mu_\chi ([X_1,f_1]\cdot [X_2,f_2])=\mu_\chi ([X_1\times X_2, f_1\circ\pi_{X_1} + f_2\circ \pi_{X_2})= $$ $$
\sum_{(x_1,x_2)\in X_1(\F_q)\times X_2(\F_q)} \chi(f_1(x_1)) \chi(f_2(x_2))=\mu([X_1,f_1]\cdot \mu[X_2,f_2]).$$

\smallskip

By precomposition with the embedding $K_0(\cV_K) \hookrightarrow KExp_K$ given by
$[X]\mapsto [X,0]$, the homomorphims $\mu_\chi$ induces a ring homomorphism
$\mu: K_0(\cV_K)\to \C$, that is, a motivic measure in the usual sense.
One has $\mu_\chi [X,0]=\# X(\F_q)$, hence this induced motivic measure is independent 
of the character $\chi$, takes values in $\Z$, and is just given 
by the usual counting function for varieties over finite fields, whose
associated zeta function $\zeta_\mu(t)$ is the Hasse--Weil zeta function. 

\smallskip

In this setting one considers the Kapranov motivic
zeta function in $KExp_K$ and the zeta function $\zeta_{\mu_\chi}$ is obtained by
composing it with the given motivic measure. Namely,
given a class $[X,f]\in KExp_K$, one considers the symmetric products
$$ Sym^n[X,f]:= [Sym^n X,  f^{(n)}], $$
where the morphism $f^{(n)}: Sym^n(X) \to \bA^1$ is given by
$$  f^{(n)} [x_1, \ldots, x_n]=f(x_1)+\cdots+ f(x_n), $$
where $[x_1, \ldots, x_n]$ is the class in $Sym^n(X)=X^n/\Sigma_n$ of
$(x_1, \ldots, x_n)\in X^n$.

\smallskip

Note that this shows, in particular, that there is a unique way of interpreting the 
term $\chi(f(x))$ for a character $\chi:\F_q\to \C^*$
when the point $x$ has degree $r>1$, namely as $\chi$ evaluated at the 
trace of $f(x)$. This follows from the description of such $x$ as an $\F_q$-point
of $Sym^n(X)$.

\smallskip

One can then define an analog in $KExp_K[[t]]$ of the Kapranov motivic zeta function as
$$ Z_{(X,f)}(t)=\sum_{n\geq 0} [Sym^n X, f^{(n)}]\, t^n, $$
and zeta functions associated to motivic measures $\mu: KExp_K \to R$, with values
in a commutative ring $R$ with values in $R[[t]]$ defined by $\zeta_\mu: KExp_K \to R[[t]]$, 
$$ \zeta_\mu((X,f),t):= \sum_{n\geq 0} \mu [Sym^n X, f^{(n)}]\, t^n . $$

\smallskip

In particular, for a choice of a character $\chi: \F_q \to \C^*$, we can consider the zeta
function associated to the corresponding motivic measure $\mu_\chi: KExp_{\F_q} \to \C$,
$$ \zeta_{\mu_\chi} ((X,f),t) = \sum_{n\geq 0} \mu_\chi [Sym^n X, f^{(n)}]\, t^n. $$
$$ = \sum_{n\geq 0} N_\chi(Sym^n X, f^{(n)})\, t^n, $$
This zeta function is a generalization of the Hasse--Weil zeta function
to which it restricts in the case $f=0$.

\smallskip

For a given pair $(X,f)$ of an $\F_q$-variety and a morphism $f: X\to \bA^1$, and a character
$\chi: \F_q \to \C^*$, consider for $r\in \N$ and $\alpha\in \C$ the sets
$$ X_{\alpha,r} :=\{ x\in X \, |\, \deg(x)=r\, \text{ and } \chi(f(x))=\alpha \}, $$
with $k(x)$ the residue field and $\deg(x)=[k(x): \F_q]$. Let $a_{\alpha,r}:=\# X_{\alpha,r}$.

\smallskip

{\bf 7.8.1. Proposition.} {\it For $(X,f)$ of an $\F_q$-variety and a morphism $f: X\to \bA^1$
and $\chi: \F_q \to \C^*$ a character as above, the zeta function $\zeta_{\mu_\chi} ((X,f),t)$
has an Euler product expansion of the form
$$ \zeta_{\mu_\chi} ((X,f),t) =\prod_\alpha \prod_{r\geq 1} (1-\alpha t^r)^{-a_{\alpha,r}}, $$
and can be written equivalently in the form
$$ \zeta_{\mu_\chi} ((X,f),t) = \exp\left( \sum_{m\geq 1} N_{\chi,m}(X,f) \, \frac{t^m}{m} \right), $$
with coefficients
$$ N_{\chi,m}(X,f)=\sum_\alpha  \sum_{r|m} r \, a_{\alpha,r} \, \alpha^{\frac{m}{r}}. $$ }

\medskip

{\it Proof.}

Note that, if we identify points $[x_1,\ldots, x_n]\in Sym^n X$ with effective divisors $D=x_1+\cdots +x_n$
on $X$ of degree $n$, then we can write
$$  \mu_\chi [Sym^n X, f^{(n)}]=\sum_D \chi(f^{(n)}(D)), $$
where
$$ \chi(f^{(n)}(D)) =\prod_{i=1}^n \chi(f(x_i)). $$
Thus we have
$$ \zeta_{\mu_\chi} ((X,f),t) =\sum_n \sum_{D\,:\, deg(D)=n} \chi(f^{(n)}(D))  \, t^n . $$

\smallskip

Let $X_{\alpha,r}$ be the level sets defined as in the statement above.
We also denote by $X_\alpha=\{ x\in X \, |\, \chi(f(x))=\alpha \}=\cup_r X_{\alpha,r}$.
Instead of considering the
integer numbers $a_r=\# X_r$ we now consider $a_{\alpha,r}:= \# X_{\alpha,r}$.

\smallskip

Observe then that we have
$$ \sum_n \sum_{D\,:\, deg(D)=n} \chi(f^{(n)}(D))  \, t^n = \sum_n \sum_{D=x_1+\cdots + x_n} \prod_{i=1}^n \chi(f(x_i)) t^{\deg(x_i)} $$
$$ = \prod_\alpha \prod_{x\in X_\alpha} (1+\alpha t^{\deg(x)} + \alpha^2 t^{2\deg(x)} + \cdots )= $$
$$ = \prod_\alpha \prod_{r\geq 1} \prod_{x\in X_{\alpha,r}} (1-\alpha t^r)^{-1} = \prod_\alpha \prod_{r\geq 1} (1-\alpha t^r)^{-a_{\alpha,r}}. $$
This gives the Euler product expansion of $\zeta_{\mu_\chi} ((X,f),t)$. 

\smallskip

Moreover, we have 
$$ \log \zeta_{\mu_\chi} ((X,f),t) = - \sum_\alpha \sum_{r\geq 1}  a_{\alpha,r} \log (1-\alpha t^r) $$
$$ = \sum_\alpha \sum_{r\geq 1} a_{\alpha, r} \sum_\ell \alpha^\ell t^{r\ell} $$
$$ = \sum_{m\geq 1}  \sum_{r|m} \sum_\alpha r \, a_{\alpha,r} \, \alpha^{m/r} \, \frac{t^m}{m}, $$
where the sum over alpha ranges over the non-empty level sets $X_{\alpha,r}\neq \emptyset$.
Consider then the sequence $N_{\chi,m} =N_{\chi,m}(X,f)$,  for $m\geq 1$, defined as
$$ N_{\chi,m}:= \sum_\alpha \sum_{r|m} r \, a_{\alpha,r} \, \alpha^{m/r}. $$
We have
$$ \log \zeta_{\mu_\chi} ((X,f),t) = \sum_{m\geq 1} N_{\chi,m}(X,f) \, \frac{t^m}{m}. $$
$\blacksquare$

\smallskip

Motivic Euler product expansions were considered, for instance, in [Bour09].
Our main focus in Proposition~7.8.1 above is on introducing a version for 
the case with exponentials.

\smallskip

In the next two subsections we discuss another categorification of the
Grothendieck ring with exponentials, which instead of using the formalism
of assemblers is based on Nori's diagrams and Nori's categories.

\medskip

{\bf 7.9. Nori's Tannakian formalism.}
We have presented in Theorem~7.4 and Propositions~7.5, 7.6, and 7.7
a categorification and spectrification of the Grothendieck ring with
exponentials based on the categorical formalism of assemblers and on
the associated spectra. There is, however, another possible approach to
categorifying the Grothendieck ring with exponentials, via an appropriate
category of motives. We discuss this other approach in this and the
following subsection, where we show that the appropriate category of
motives is provided by the exponential motives of Fres\'an and Jossen, [FreJo20],
constructed through the general Tannakian formalism of Nori motives,
which we review in this subsection. Thus, assemblers and Nori diagrams
can be viewed as two complementary paths to the categorification of
Grothendieck rings, one that leads naturally to the homotopy-theoretic
world and the other to the motivic. The possible interactions between
homotopy-theoretic and motivic settings appear to be most promising
for future developments.

\smallskip

We recall here briefly Nori's formalism, constructing Tannakian categories
associated to diagrams and their representations, and the application of
this formalism to the construction of the category of Nori motives.
As we recall below, a Nori diagram is like a quiver 
and the Nori formalism makes it possible to construct from representations of
Nori diagrams in categories $Mod_R$ of modules
an abelian (and under suitable circumstances Tannakian)
category that satisfied a universal property with respect to such
representations. We recall here the main steps of
this construction. 
A main reference for the material we review in this subsection is the book
[HuMu-St17]. What we review here will be useful in the next subsection,
where we describe another way of categorifying the Grothendieck ring
of varieties with exponentials in terms of a category of exponential motives,
due to Fres\'an and Jossen [FreJo20].

\smallskip

{\it A category of Nori diagrams} is defined as follows (Definition~7.1.1 of [HuMu-St17]).
{\it A diagram $D$} consists of a family $V(D)$ of vertices and a family $E(D)$ of
edges, with a boundary map $\partial: E(D) \to V(D)\times V(D)$, where
$\partial(e)=(\partial_{out}(e), \partial_{in}(e))$ means source and target of the
oriented edge. {\it A morphism} $D_1\to D_2$ of diagrams consists of
two maps $V(D_1)\to V(D_2)$ and $E(D_1)\to E(D_2)$ compatible with
orientations and boundaries. {\it A diagram with identities} is a diagram $D$
where for each $v\in V(D)$ there is a unique oriented edge $id_v$ with $\partial(id_v)=(v,v)$.
In the case of diagrams with identities, morphisms are require to map identity
edges to identity edges. 

\smallskip

To a category $\cC$ one can associate the diagram  $D(\cC)$ with
$V(D(\cC))=Obj(\cC)$ and $E(D(\cC))=Hom_\cC$.  {\it A representation
of a diagram $D$ in a category $\cC$} is a morphism of diagrams $T: D\to D(\cC)$.
One considers representations in categories $Mod_R$ of modules
over a commutative ring $R$, and in particular representations in categories of vector spaces. 

\smallskip

Given a diagram $D$ and a representation $T$ in $Mod_R$, for some commutative ring $R$, 
one defines the ring $End(T)$ as
$$ End(T)=\left\{ (\phi_v)\in \prod_{v\in V(D)} End_R(T(v))\,\biggr|\, \phi_{\partial_{out}(e)} \circ T(e) = T(e) \circ \phi_{\partial_{in}(e)}, \ \forall e\in E(D) \right\}. $$

\smallskip

Nori's {\it diagram category} $\cC(D,T)$ is then obtained in the following way. If the diagram $D$ is finite then
$\cC(D,T)$ is the category $Mod_{End(T)}$ of finitely generated $R$--modules with an $R$--linear
action of $End(T)$. If the diagram $D$ is infinite, one considers all finite subdiagrams $D_F$ and
constructs the corresponding categories $\cC(D_F, T|_{D_F})$. The category $\cC(D,T)$  has
as objects the union of all the objects of the $\cC(D_F, T|_{D_F})$ for all finite $D_F\subset D$.
An inclusion $D_F \subset D_F'$ determines a morphism $End(T|_{D_F'}) \to End(T|_{D_F})$
by projecting the product $\prod_{v\in V(D_F')} End_R(T|_{D'_F}(v))$ onto $\prod_{v\in V(D_F)} End_R(T|_{D_F}(v))$.
This morphism induces a functor from $Mod_{End(T|_{D_F})}$ to $Mod_{End(T|_{D_F'})}$. Morphisms
in $\cC(D,T)$ are then defined as colimits of morphisms in  $\cC(D_F, T|_{D_F})$ under these extensions. 
The category $\cC(D,T)$ constructed in this way is $R$--linear abelian with an $R$--linear faithful exact forgetful
functor $f_T: \cC(D,T) \to Mod_R$. The representation $T: D \to Mod_R$ factors as $T= f_T \circ \tilde T$ with
a representation $\tilde T: D \to \cC(D,T)$. 
We refer the reader to [HuMu-St17] pp.~140--144 for more details. 

\smallskip

Diagram categories satisfy the following universal property:

\smallskip

{\it  Any representation $F: D \to A$ where $A$
an $R$--linear abelian category with an $R$--linear faithful exact functor $f: A \to Mod_R$ factors
through a faithful exact functor $L(F): \cC(D,T)\to A$, where $T=f\circ F$, compatibly with the
decomposition $T=f_T\circ \tilde T$.}  (See [HuMu-St17], pp.~140--141).

\smallskip

Tensor structures on diagram categories can be obtained through the notion of
graded diagrams with a commutative product with unit, in
the sense of Definition~8.1.3 of [HuMu-St17].  Here are some details.

\smallskip

{\it A graded diagram $D$} is a diagram endowed with a map $\deg: V(D) \to \Z/2\Z$
extended to $\deg: E(D)\to \Z/2\Z$ by $\deg(e)=\deg(s(e))-\deg(t(e))$. 
The product $D\times D$  is the diagram with vertices the pairs $(v,w)\in V(D)\times V(D')$
and edges of the form $(e,id)$ or $(id,e')$. {\it A product structure} on $D$ is a map of
graded diagrams (a degree preserving map of directed graphs) $D\times D \to D$
together with a choice of edges
$$ 
\alpha_{v,w}: v\times w \to w\times v, \ \ \forall v,w\in V(D) 
$$
$$ 
 \beta_{v,w,u}: v\times (w\times u) \to (v\times w)\times u, 
 $$
 $$
 \beta'_{v,w,u}: (v\times w)\times u \to v\times (w\times u),  . 
 $$
 for all  $v,w,u\in V(D)$.
 A unit is a vertex $\bold{ 1}$ with $\deg(\bold{ 1})=0$ and edges $u_v: v \to \bold{ 1}\times v$ for
 all $v\in V(D)$.
 
\smallskip

One can then consider those 
representations of $D$ that are compatible with the grading and with the commutative
product above, where the compatibility is expressed as the existence of isomorphisms 
(see Definition 8.1.3 of [HuMu-St17])
$$ 
\xymatrix{\tau_{v,w}: T(v\times w) \ar[r]^{\simeq} & T(v)\otimes T(w) }
$$ 
for all $v,w\in V(D)$, with the following properties:
$$
\xymatrix{T(v)\otimes T(w) \ar[r]^{\tau_{v,w}^{-1}} & T(v\times w) \ar[r]^{T(\alpha_{v,w})}  &  T(w\times v) \ar[r]^{\tau_{w,v}} 
& T(w)\otimes T(v)}
$$
is equal to multiplication by $(-1)^{\deg(v) \deg(w)}$; the $\beta$--maps satisfy $T(\beta_{v,w,u})^{-1}=T(\beta'_{v,w,u})$, 
and moreover
$$
\tau_{v,w'}\circ T(1,e) = (id \otimes T(e))\circ \tau_{v,w} :\, T(v\times w) \to T(v)\otimes T(w'),
$$
$$
\tau_{v',w}\circ T(e,1) = (T(e)\otimes id)\circ \tau_{v,w} :\, T(v\times w)\to T(v')\otimes T(w),
$$
$$ \xymatrix{ T(v \times (w\times u)) \ar[r]^{T(\beta_{v,w,u})}
\ar[d]^{\tau\circ \tau} & T((v\times w)\times u) \ar[d]^{\tau\circ\tau} \\ 
T(v)\otimes (T(w)\otimes T(u)) \ar[r]^{\simeq} & (T(v)\otimes T(w))\otimes T(u)) } 
$$
and similarly for the inverse $T(\beta'_{v,w,u})$.

\smallskip

In Theorem~7.1.12 of [HuMu-St17] it  is shown that if the representation $T$ is valued in
the subcategory of $Mod_R$ of finite projective modules, and $R$ is a Dedekind domain,
the Nori diagram category $\cC(D,T)$ is equivalent to the
category of finitely generated comodules over the coalgebra $\cA(D,T)$ given
by the colimit 
$$ 
\cA(D,T) = colim_{D_F} End(T |_{D_F})^\vee 
$$
over finite sub-diagrams $D_F$. 

\smallskip

In Sec. 7.5.1 of [HuMu-St17] it is ahown, that if
$R$ is a Dedekind domain, then for the $R$--algebra $E=End(T|_{D_F})$
with $D_F$ a finite diagram,
the $R$--dual $E^\vee=Hom_R(E,R)$ has the property that the canonical map
$E^\vee \otimes_R E^\vee \to Hom(E,E^\vee)\simeq (E\otimes_R E)^\vee$
is an isomorphism. If an $E$-module  is finitely generated projective as
an $R$--module then the same is true for comodule $E^\vee$.  The coalgebra
$\cA(D,T)$ also carries an algebra structure induced by the monoidal structure of
$\cC(D,T)$ discussed in Sections 7.1.4 and 8.1 of [HuMu-St17], so that $\cA(D,T)$
determines a pro--algebraic monoid scheme $Spec(\cA(D,T))$ (see 
Section 7.1.4 of [HuMu-St17]). This is the general form of Nori's Tannakian formalism.

\smallskip

More specifically, for $K$ a subfield of $\C$, the category of Nori motives (which we here call Nori {\it classical}
motives, not to be mixed with {\it exponential} motives we will review in the next subsection)
is constructed by considering the Nori diagram of {\it effective pairs} (see [HuMu-St17], pp.~207--208)
where vertices are of the form $(X,Y,n)$ with $X$ a $K$--variety, $Y\subseteq X$ a closed embedding,
and $n$ an integer, and non--identity edges are of the following types:

\smallskip

(a) Let $(X,Y)$ and $(X^{\prime}, Y^{\prime})$ be two pairs
of closed embeddings.
Every morphism $f: X\to X^{\prime}$ such that $f(Y)\subset Y^{\prime}$
produces  functoriality edges $f^*$ (or rather $(f^*,n)$) going from $(X^{\prime}, Y^{\prime},n)$
to $(X, Y, n)$.

\smallskip

(b) Let $(Z\subset Y\subset X)$ be a pair of closed embeddings. Then
it defines coboundary edges  $\partial$ from $(Y,Z,n)$ to $(X,Y, n+1)$.

\smallskip

A representation of the diagram $D(Pairs^{eff})$ obtained in this way in the
category $Mod_\Z$ (or in the category $Vect_\Q$) is given by relative singular cohomology
$H^n(X(\C), Y(\C), \Z)$ (respectively, $H^n(X(\C), Y(\C), \Q)$). 
The category of effective Nori motives is given by $\cC(D(Pairs^{eff}), H^*)$ and the
category of Nori motives is obtained as the localization at $(\bG_m, \{ 1 \}, 1)$.
A tensor structure is obtained as discussed above, after restricting to a subcategory
of {\it good pairs}, see Sections~8.1 and 9.3 of [HuMu-St17]. One obtains a Tannakian
category of classical Nori motives. 

\smallskip

{\bf 7.9.1. Spectra from Nori diagrams.} The categorification through assemblers
discussed in the previous sections leads to the construction of an associated
spectrum, obtained through a $\Gamma$-space, as discussed in Section~3.3.2. It is
then natural to ask whether the categorical construction outlined here above can also
have an associated homotopy-theoretic spectrum. 

\smallskip

As recalled in Section~1.5, the formalism of $\Gamma$-spaces provides a general
method for the construction of spectra from categories, through the nerve of the associated
category of summing functors. 

\smallskip

Given a Nori category $\cC(D,T)$, obtained as above from a Nori
diagram and a representation $T$, let $F_{\cC(D,T)}$ be the associated
$\Gamma$-space that maps a finite pointed set $X$ to the pointed simplicial
set given by the nerve $\Cal{N}(\Sigma_{\cC(D,T)}(X))$ of the category $\Sigma_{\cC(D,T)}(X)$
of summing functors $\Phi: P(X)\to \cC(D,T)$, with the notation of Section~1.5.
This $\Gamma$-space determines then a spectrum $\bS(F_{\cC(D,T)})$, obtained by
promoting the functor $F_{\cC(D,T)}$ to an endofunctor of the category
of pointed simplicial sets and applying it to spheres, see Section~1.5.

\smallskip

Since the category $\cC(D,T)$ is abelian, its higher algebraic $K$-theory groups
[Qui73] are the homotopy groups of an infinite loop space $K(\cC(D,T))$. The
spectrum $\bS(F_{\cC(D,T)})$ constructed as above provides a delooping of
this infinite loop space, in the sense discussed in [Carl05]. 

\medskip

{\bf 7.10. Exponential motives.} As mentioned in the previous section,
another possible way to categorify the notion of Grothendieck 
ring with exponentials is through a category of ``motives with
exponentials". A general treatment of such ``exponential motives"
was presented by Fres\'an and Jossen in [FreJo20], where the relations between the
resulting category of exponential motives, 
the Grothendieck ring of varieties
with exponentials of [ChamLoe15], the formalism of
Nori motives [HuMu-St17], and the ``exponential periods" 
of Kontsevich--Zagier [KoZa01] is also explained. 

\smallskip

We assume here that $K$ is a subfield of $\C$. 
The category $MotExp(K)$ of exponential motives over a field $K$ constructed in [FreJo20] is based on Nori diagrams
where the vertices are tuples $(X,Y,f,n,i)$ with $X$ a variety over $K$ and $Y\hookrightarrow X$
a closed subvariety, together with a morphism $f: X \to \bA^1$, with the restriction $f|_Y : Y \to \bA^1$,
and integers $n,i$, respectively referred to as degree and twist; the edges are of the following types:

\smallskip

(1) a morphism of varieties $h: X \to X'$ with $h(Y)\subseteq Y'$ and $f'\circ h =f$ determines an
edge $h^*: (X',Y',f',n,i) \to (X,Y,f,n,i)$;

\smallskip

(2) a pair of closed immersions $Z\subseteq Y \subseteq X$ determines an edge
$$ \partial: (Y,Z,f|_Y, n-1,i) \to (X,Y,f,n,i); $$

(3) edges of the form
$$ (X\times \bG_m, (Y\times \bG_m)\cup (X\times \{ 1 \}), f \boxtimes 0 , n+1,i+1) \to (X,Y,f,n,i) $$ 

\smallskip

A representation of these Nori diagrams in the category of vector spaces
$Vect_\Q$ is given in [FreJo20] in terms of {\it rapid decay cohomology}
$$ 
H^n(X,Y,f)(i) 
$$
where $(i)$ denotes the Tate twist given by the tensor product with tensor powers of $H^1(\bG_m,\Q)$ 
and the rapid decay cohomology $H^n(X,Y,f)$.

\smallskip

{\it Rapid decay cohomology} is constructed in the following way. Let $X$ be a complex variety with
a closed subvariety $Y\subseteq X$ and a regular function $f: X \to \bA^1$. For a real number $r\in \R$,
denote by $S_r=\{ z\in \C\,|\, Re(z)\geq r \}$ the closed half-plane with boundary the vertical line $Re(z)=r$.
One then defines
$$ 
H_n(X,Y,f) =\varinjlim_{r\to \infty} H_n(X,Y\cup f^{-1}(S_r)), \ \ \ \  H^n(X,Y,f) =\varprojlim_{r\to \infty} H^n(X,Y\cup f^{-1}(S_r)), 
$$
as limits in $Vect_\Q$, with respect to the system of inclusions $f^{-1}(S_t) \subseteq f^{-1}(S_r)$ for $t\geq r$.
This cohomology theory originates in the study of differential equations with irregular singularities, [DeMaRa07].

\smallskip

The reader should also look back at the definition of convergence boundaries in Sec. 2.3.

\smallskip

A property of rapid decay cohomology that is directly relevant to exponential motives and the Grothendieck
ring of varieties with exponentials is the fact that one has
$$ H^n(X\times \bA^1, Y\times \bA^1, \pi_{\bA^1})=0, $$
which suggests that this is indeed the right cohomology theory that provides a realization compatible
with the relations in the Grothendieck ring of varieties with exponentials that we discussed above. 

\smallskip

A more general construction of rapid decay cohomology is given in [FreJo20] in terms of perverse sheaves.
For our purposes we only review the more elementary definition above, and we refer the
reader to [FreJo20] for the more general treatment. 

\smallskip

The category $MotExp(K)$ of exponential motives is the abelian $\Q$--linear category obtained by applying
the Nori formalism to these diagrams and representations. A tensor product on the category of
exponential motives is also introduced in [FreJo20] following an analogous construction for
Nori motives recalled in the previous subsection, and it is shown that each object admits a dual. 
The general Nori formalism then
shows that the resulting category $MotExp(K)$ of exponential motives is Tannakian, as
discussed in the previous subsection (see [HuMu-St17] for more details). 

\smallskip

The relations between the category $MotExp(K)$ of exponential motives and other motivic theories discussed in [FreJo20]
that are more directly relevant for us here can be
summarized as follows. There is a fully faithful exact functor from the category of classical 
Nori motives to exponential motives, which gives rise to a morphism (in the reverse direction)
of affine group schemes, between the respective motivic Galois groups, which is faithfully flat. 
Moreover, it is shown in Section~5.4 of [FreJo20] that there is a unique ring homomorphism 
$$ \chi:  KExp_K \to K_0(MotExp(K)) $$
such that, for a pair $(X,f)$ of a $K$-variety and a morphism $f: X \to \bA^1$ 
$$ \chi[X,f] = \sum_n (-1)^n [H^n_c(X,f)], $$
where the bracket notation indicates that we view the $H^n_c(X,f)$ on the
right-hand-side as elements in $K_0(MotExp(K))$.

\smallskip

Exponential periods arise from the pairing of de Rham cohomology and rapid decay homology
$$
 H^n_{dR}(X,f) \otimes H_n(X,f) \to \C . 
 $$
 
{\it Rapid decay homology} describes cycles that are possibly non--compact but unbounded in directions where
$Re(f)\to \infty$, so that the exponential periods are given by the pairing
$$
 \int_\gamma e^{-f} \omega = \lim_{r\to \infty} \int_{\gamma_r} e^{-f} \omega 
 $$
where $\gamma$ is a direct limit of compact cycles $\gamma_r$ for $r\in \R$ with 
$\partial \gamma_r \subset f^{-1}(S_r)$. 

An example of a number that is not expected to be a period of a classical motive but that
is an exponential period (for $f(x)=x^2$) is
$$
 \sqrt\pi = \int_\R e^{-x^2} dx. 
 $$
 
Exponential motives provide the motivic framework for exponential periods.

\medskip

{\bf 7.11. Motivic Fourier transform.} 
Let $K$ denote an algebraically closed field of characteristic zero,
and $V$ a finite--dimensional linear space over $K$.
A motivic Fourier transform on the Grothendieck ring $KExp_V$ of varieties with exponentials
over $V$, can be defined as follows (see Section 7.1 of [CluLoe10], Section 1.2 of [ChamLoe15]
and Section~7.12 and Definition~2.2 of [Wy17]):
$$ \cF: KExp_V \to KExp_{V^\vee} $$
$$
 \cF( [X,f]_V) := [X \times V^\vee, f \circ \pi_X + \langle u \circ \pi_X , \pi_{V^\vee} \rangle ]_{V^\vee},
 $$
where $V^\vee=Hom_K(V,K)$ is the dual linear space,   $\langle\cdot,\cdot\rangle: V \times V^\vee \to K$
the natural pairing, and $u: X \to V$ is the structure morphism of $X$ as a $V$--variety. 
This motivic Fourier transform satisfies the relation
$$
\cF\circ \cF [X,f]_V = \bL^{\dim V} \cdot j^* [X,f]_V,
$$
where $j^*: KExp_V\to KExp_V$ is the pullback induced by the map $j: V \to V$ given by 
multiplication by $-1$, see [CluLoe10]. 

\smallskip

We consider here again the assemblers $\cC^{KExp}_S$ of Proposition~7.6 and
the associated homotopy--theoretic spectra $K(\cC^{KExp}_S)$ underlying 
the relative Grothendieck ring $KExp_S$ of varieties with exponentials over a base scheme $S$, 
as we discussed in Theorem~7.4 and Propositions~7.5 and 7.6.

\medskip

{\bf 7.12. Theorem.} {\it
The motivic Fourier transform $\cF$ lifts to an morphism of assemblers $\cF: \cC^{KExp}_V \to \cC^{KExp}_{V^\vee}$
and induces a morphism of the associated spectra.
The second iterate $\cF\circ \cF$ of the Fourier transform determines a covering family
$$
 \{ (\bA^d,0)\times (Z,f|_Z) \hookrightarrow \cF\circ \cF (X,f), (Z^c\times V^\vee, h|_{Z^c}) \hookrightarrow
\cF\circ \cF (X,f) \}
 $$
 in the assembler $\cC^{KExp}_V$,
 where $Z=\{ (x,v)\in X\times V\,|\, u(x)+v =0 \}$ with $u: X \to V$ is the structure morphism of $X$ as a $V$--variety.
 The term $(Z^c\times V^\vee, h|_{Z^c})$ in the range of the endofunctor 
$\Phi: \cC^{KExp}_V \to \cC^{KExp}_V$ of the assembler $\cC^{KExp}_V$. This family implements
the relation 
$$
\cF\circ \cF [X,f]_V = \bL^{\dim V} \cdot j^* [X,f]_V 
$$
in $KExp_V=\pi_0 K(\cC^{KExp}_V/\Phi)$.
}

\medskip  

{\it Proof.} The argument is similar to Section~2 of [Wy17].
For $(X,f)_V$ an object in $\cC^{KExp}_V$,
the second iterate of the Fourier transform is given by
$$ \cF \circ \cF (X,f) = (X\times V \times V^\vee, f\circ \pi_X + \langle u\circ \pi_X + \pi_V, \pi_{V^\vee} \rangle), $$
where $u: X \to V$ is the structure morphism of $X$ as a $V$--variety.

\smallskip

Consider the variety 
$$ Z=\{ (x,v)\in X\times V\,|\, u(x)+v =0 \}. $$
The structure morphism of $Z$ as a $V$--variety given by the second projection $\pi_V: Z \to V$
fits in the following commutative diagram, with $j(v)=-v$:
$$
 \xymatrix{ Z \ar[r]^{\pi_X} \ar[d]^{\pi_V} & X \ar[d]^u \\ V \ar[r]^j & V } 
 $$ 

\smallskip

Consider the embeddings
$$ 
(Z \times V^\vee, f|_Z \circ \pi_Z) \hookrightarrow (X\times V \times V^\vee, h) 
\hookleftarrow (Z^c\times V^\vee, h|_{Z^c}), 
$$
where $Z^c$ is the complement of $Z$ in $X\times V$ and $h: X\times V\times V^\vee \to K$ is given by 
$$
 h(x,v,v^\vee) = \langle u(x) + v, v^\vee \rangle + f(x). 
 $$
We have 
$$
(Z \times V^\vee, f|_Z \circ \pi_Z)=(\bA^d,0)\times (Z,f|_Z).
$$ 

On the other hand, for the object
$(Z^c\times V^\vee, h|_{Z^c})$ of the assembler $\cC^{KExp}_V$ we can see that, for 
any $(x,v)\in Z^c$, the object
$$
(\{ (x,v) \} \times V^\vee, h_{(x,v)} )$$ with $$h_{(x,v)}(v^\vee)=f(x)+\langle u(x) + v, v^\vee \rangle
$$
is in the range of endofunctor $\Phi: \cC^{KExp}_V \to \cC^{KExp}_V$
of the assembler $\cC^{KExp}_V$. This in fact follows from the following general observation.

\smallskip

If $W$ is a finite dimensional $K$--vector space and $\lambda: W \to K$ is a linear map, then
the pair $(W,\lambda)$ is either $(\bA^d, 0)$ when $\lambda$ is trivial, or a product 
$$
(W'\times \bA^1, \lambda\circ \pi_W + \pi_{\bA^1}),
$$
where the $\bA^1$ factor is spanned by a vector $w_0$ in $W$ such that $\lambda(w_0)=1$ in $K$, so
that writing $w=w'+t v_0$ with $w'\in W'$ gives $\lambda(w)=\lambda(w')+t$. Thus, for $\lambda$ nontrivial, 
such $(W,\lambda)$ is
in the range of the endofunctor $\Phi: \cC^{KExp}_K \to \cC^{KExp}_K$
of the assembler $\cC^{KExp}_K$. We can then apply this to the case where 
$(W,\lambda)$ is given by 
$$
(\{ (x,v) \} \times V^\vee, h_{(x,v)})
$$ 
for a fixed $(x,v)$.
This shows that the map
$\varphi: Z^c \to Obj(\cC^{KExp}_K)$ with 
$$
\varphi(x,v)=(\{ (x,v) \} \times V^\vee, h_{(x,v)} )
$$
has image contained in the range of $\Phi$. 

\smallskip

This implies the existence of some object $(W_{(x,v)}, h_{(x,v)}|_W)$ in $\cC^{KExp}_K$, obtained
as in the general observation above, such that
$$
\varphi(x,v)=(W_{(x,v)}, h_{(x,v)}|_W) \times (\bA^1, id),
$$ 
for a subspace $W_{(x,v)} \subset V^\vee$. More precisely,
 $$
 Z^c\times V^\vee=\cup_{(x,v)\in Z^c} W_{(x,v)} \times \bA^1
 $$ 
 with compatible morphisms to $\bA^1$,
and the decomposition given by identifying $\bA^1$ with the span of a vector 
$w_{(x,v)}^\vee\in V^\vee$ such that 
$$
f(x)+\langle u(x) + v, w_{(x,v)}^\vee \rangle=1
$$ 
in $\bA^1$.
The bundle $W$ over $Z^c$ constructed in this way is locally trivial, hence by Noetherian descent induction 
we then have a finite decomposition $\{ Z_i \}$ with locally closed $Z_i\subset Z^c$ such that 
$$
(Z_i\times V^\vee, h|_{Z_i\times V^\vee})
$$ 
is isomorphic to a product 
$$ 
\xymatrix{ \xi: Z_i\times V^\vee \ar[r]^{\simeq} & Z_i \times W \times \bA^1 } 
$$
with, for $g=h\circ \xi$,
$$ 
(Z_i \times W \times \bA^1, g\circ \pi_{Z_i\times W}+ \pi_{\bA^1}). 
$$ 
Thus $(Z^c\times V^\vee, h|_{Z^c})$ is itself in the range of the endofunctor 
$\Phi: \cC^{KExp}_V \to \cC^{KExp}_V$ of the assembler $\cC^{KExp}_V$. 

\smallskip

This induces the relation
$$
\cF\circ \cF [X,f]_V = \bL^{\dim V} \cdot j^* [X,f]_V 
$$
in $KExp_V=\pi_0 K(\cC^{KExp}_V/\Phi)$.
$\blacksquare$

\medskip

{\bf 7.13. Motivic Bruhat--Schwartz functions and Poisson summation.}
In this subsection we review the Hrushovski--Kazhdan motivic Poisson summation formula,
following [ChamLoe15] and [HruKaz09].

\smallskip

As in Section~2.1 above, with slightly changed notations, consider a global
field  $F$ be a global field and denote by $A_F$ its   ad\`eles group. For  a place $v$ of $F$,
let $F_v$ denote the completion of $F$ at $v$. The space $\cS(F_v)$ of Bruhat--Schwartz 
functions on $F_v$ consists of rapidly decaying complex valued
functions when $v$ is archimedean and of locally constant and compactly supported 
complex valued functions when $v$ is archimedean.  Bruhat--Schwartz functions
$\varphi \in \cS(A_F)$ are linear combinations of complex valued functions of the form $\prod_v \varphi_v$, 
where $v$ ranges over the places of $F$ and $\varphi_v\in \cS(F_v)$, with the property that $\varphi_v =1_{\cO_v}$,
the characteristic function of the ring of integers $\cO_v\subset F_v$ at all but finitely many 
of the non--archimedean places. 

\smallskip

Since $F$ is a discrete cocompact subgroup of the ad\`eles group $A_F$, the Poisson
summation formula gives
$$ 
\sum_{x\in F} \varphi(x) =\frac{1}{\mu(A_F/F)} \sum_{y\in F} \hat \varphi(y), 
$$
where $\hat \varphi$ is the Fourier transform of  $\varphi$:
$$
 \hat \varphi(y) :=\int_{A_F} \varphi(x) \chi(xy) d\mu(x) 
 $$
for a nontrivial character $\chi: A_F \to \C^*$ and  the Haar measure $\mu$.
This classical fact 
has a motivic counterpart (for function fields), due to Hrushovski and Kazhdan ([HruKaz09]), 
where Bruhat--Schwartz functions are replaced by elements of a relative Grothendieck ring with
exponentials, and motivic Fourier transform is used. 

\smallskip

We will focus here on the case of function fields $F$ for a curve $C$ over a finite field,
so there are only non--archimedean places corresponding to the points of the curve. 
As discussed in Section~5.2 of [Bilu18], Section~1.2 of [ChamLoe15], Section~4 of [HruKaz09], 
the Bruhat--Schwartz functions
$\varphi\in \cS(F_v)$ can be organized according to their level, labelled by two integers $(M,N)$, respectively 
measuring ``support" and ``invariance". Namely, for a function $\varphi\in \cS(F_v)$ there are integers $M,N\geq 0$ 
such that $\varphi\equiv 0$ outside of $t^{-M}\cO_v$ and $f$ is invariant modulo the subgroup $t^N \cO_v$,
where $t$ is a uniformizer. Thus, $\varphi$ descends to a function on the quotient $t^{-M}\cO_v/t^N \cO_v$, identified with
an $M+N$ dimensional vector space $\bA_{k_v}^{M+N}$ over the residue field $k_v$. This vector space is denoted by
$\bA_{k_v}^{(M,N)}$, to keep track of the levels. 

\smallskip

The motivic version of Bruhat--Schwartz functions of level $(M,N)$ on a nonarchimedean local field 
with residue field $k$ are then defined (see Section~5.2 of [Bilu18], Section~1.2 of [ChamLoe15], Section~4 of [HruKaz09]) 
as elements in the relative localized Grothendieck ring with exponentials $Exp\cM_{\bA^{(M,N)}_k}$.

\smallskip

This definition is motivated by regarding, as we discussed earlier, classes in the Grothendieck 
ring with exponentials $KExp_{\bA^{(M,N)}_k}$
as motivic versions of functions $\Psi_{[X,f],\chi}: \bA^{(M,N)}(k) \to \C$ of the form
$$ 
\varphi(s):=\Psi_{[X,f],\chi}(s)=\sum_{x\in X_s(k)} \chi(f(x)). 
$$

\smallskip

Consider the canonical morphism $\pi: \bA^{(M,N)}_k \to Spec(k)$ and the
induced morphism of additive groups $\pi_! : Exp\cM_{\bA^{(M,N)}_k} \to Exp\cM_k$ (summation
along the fibers) 
$$
 \pi_!: [X,f]_{\bA^{(M,N)}_k} \mapsto  [X,f]=[X,f]_k . 
 $$
 Note that in this case there is only one fiber since the base is $Spec(k)$.

\smallskip

One then defines, for an element $\varphi =[X,f]_{\bA^{(M,N)}_k} \in Exp\cM_{\bA^{(M,N)}_k}$,
$$ 
\int \varphi(s) := \bL^{-N} \pi_! \varphi  \, \in Exp\cM_k. 
$$

\smallskip

The dependence on the level $(M,N)$ is regulated by the behaviour under the
inclusions $\iota: \bA_k^{N-M}\hookrightarrow \bA_k^{N-(M-1)}$:
$$ 
\iota: (x_M, \ldots, x_{N-1}) \mapsto (0, x_M, \ldots, x_{N-1}) 
$$
and projections $\rho: \bA_k^{(N+1)-M} \to \bA_k^{(N-M}$:
$$
 \rho: (x_M, \ldots, x_N) \mapsto (x_M, \ldots, x_{N-1}), 
 $$
with $\iota^* \iota_! =Id$ and $\rho_! \rho^* =\bL \cdot $.
The motivic Bruhat--Schwartz functions of all levels for a non--archimedean local field with residue field $k$ are then defined as
$$
 \cS^{mot} := \varinjlim_{M, \iota_!} \varinjlim_{N, \rho^*} Exp\cM_{\bA^{(M,N)}_k}. 
 $$

\smallskip

For a smooth projective curve $C$ over $k$, with $F=k(C)$, and for a place $v\in C(k)$ with
residue field $k_v$, let $Res_{k_v/k}$ denote the restriction of scalars functor, satisfying
$$
 Res_{k_s/k}(\bA^m_{k_v})=\bA^{m\,[k_v:k]}_k . 
 $$
The ad\`elic version is then defined as
$$
 \cS^{mot} (A_F):= \varinjlim_{B \subset C(k)} \cS^{mot}(\prod_{v\in B} F_v), 
 $$
where the limit is taken over all finite subsets $B$ of the set of points $C(k)$ and 
$\cS^{mot}(\prod_{v\in B} F_v)$ is obtained in the following way. Given a finite set $B$
of places, and given levels $M_v, N_v$, define 
$$
\bA^{(M_B,N_B)}_k := \prod_{v\in B}  Res_{k_s/k}(\bA^{(M_v,N_v)}_{k_v}) = \prod_{v\in B}  \bA_k^{(M_v [k_v:k],N_v [k_v:k])},
 $$
and consider the Grothendieck ring with exponentials $Exp\cM_{\bA^{(M_B,N_B)}_k}$. One then sets
$$
 \cS^{mot}(\prod_{v\in B} F_v):= \varinjlim_{M_B,\iota_!} \varinjlim_{N_B, \rho^*} Exp\cM_{\bA^{(M_B,N_B)}_k}. 
 $$
For $B'\supset B$, taking products with the unit elements $1_{\cO_v}$ for $v\in B'\smallsetminus B$ gives the maps 
$$ \cS^{mot}(\prod_{v\in B} F_v) \to \cS^{mot}(\prod_{v\in B'} F_v) $$ that determine the directed system
computing $\cS^{mot} (A_F)$.

\smallskip

The motivic Fourier transform discussed in the previous subsection adapts to this setting of 
motivic ad\`elic Bruhat--Schwartz functions. One proceeds as follows (see Section~1.3 of [ChamLoe15]
and Section~5 of [HruKaz09]).
Consider a nontrivial $k$--linear map $r_v : F_v \to k$ of conductor $c$, meaning the smallest integer $c$ such
that $r$ vanishes on $t^c \cO_v$. Such linear maps can be obtained geometrically from the residue  $res_v: \Omega_{F_v/k}\to k$ 
at the point $v\in C(k)$, of a chosen meromorphic differential form $\omega\in \Omega_{F_v/k}$ by 
setting $r_\omega: x \mapsto res_v(x\omega)$, with $c$ the order of pole of $\omega$ at $v$. Such a $k$-linear map $r_v : F_v \to k$
induces linear morphisms $r^{(M,N)}: \bA^{(M,N)}_k \to \bA^1_k$, for $N\geq c$ (see Section~1.2.6 of [ChamLoe15]).
The multiplication in $F_v$ induces morphisms
$$
 \bA_k^{(M,N)}\times \bA_k^{(M',N')} \to \bA_k^{(M+M',N'')}, \ \ \ N''=\min\{ M'+N, M+N' \}. 
 $$
For $N''\geq c$, composing with $r^{(M+M',N'')}$ gives a morphism
$$ 
\bA_k^{(M,N)}\times \bA_k^{(M',N')} \to \bA^1_k, 
$$
which we write as 
$$
 (v,v')\mapsto r(vv'). 
 $$

\smallskip

Then, for $\varphi=[X,f]_{\bA_k^{(M,N)}}$
in $Exp\cM_{\bA^{(M,N)}_k}$, one defines the motivic Fourier transform as
$$ 
\cF_\omega \varphi := \bL^{-N} [X\times_{\bA^{(M,N)}_k} \bA^{(M,N)}_k \times \bA_k^{(M',N')}, f\circ \pi_X +  \langle   
u\circ \pi_1, \pi_2\rangle_\omega ]_{\bA_k^{(M',N')}}, 
 $$
with $(M',N')=(c-N, c-M)$, and where 
$$ 
X\times_{\bA^{(M,N)}_k} \bA^{(M,N)}_k \times \bA_k^{(M',N')}
 $$
is the fibered product over the structure morphism $u: X \to \bA^{(M,N)}_k$ and the first projection of
$\bA^{(M,N)}_k \times \bA_k^{(M',N')}$, with the structure of $\bA_k^{(M',N')}$--variety given by the
projection onto the $\bA_k^{(M',N')}$ factor. The morphism to $\bA^1_k$ is given by
$$
 f(x) +  \langle   v, v'\rangle_\omega,
  $$
for $x\in X$ and $(v,v')\in \bA^{(M,N)}_k \times \bA_k^{(M',N')}$, with
$$
 \langle   v, v'\rangle_\omega := r_\omega(vv'). 
 $$

\smallskip

With the integral notation discussed above, the local motivic Fourier transform is also written in the form 
$$ \cF_{\omega_v} \varphi =\int \varphi(x) e(xy) dx, $$
where $e(xy)$ stands for the motivic Fourier kernel described explicitly here above.

\smallskip

As discussed earlier, this motivic Fourier transform satisfies a Fourier inversion formula, which in this
case takes the form
$$ 
\cF_\omega \circ \cF_\omega [X,f]_{\bA_k^{(M,N)}} = \bL^{-c} j^* [X,f]_{\bA_k^{(M,N)}}, 
$$
where $j^*$ is the pullback on $KExp_{\bA_k^{(M,N)}}$ of the map on $\bA_k^{(M,N)}$ given by multiplication by $-1$.

\smallskip

The motivic Fourier transform defined as above extends from the local to the global case. Given a finite
set $B \subset C(k)$ of places of the global function field $F=k(C)$, consider an ad\`elic motivic Bruhat--Schwartz
function in $\varphi \in \cS^{mot}(\prod_{v\in B} F_v)$, defined as discussed above. 

\smallskip
Let $\omega \in \Omega_{F/k}$ be a meromorphic
differential form,  and $D=\sum_v c_v \, v$ its divisor,  with $\deg(D)=\chi(C)=2-2g$. 
This defines linear maps $r_v=r_{\omega_v}: F_v \to k$ by $r_v(x)=res_v(x\omega)$, the residue at $v\in C(k)$.

\smallskip

For $\varphi=[X,f]_{\bA_k^{(M_B,N_B)}}$ in $Exp\cM_{\bA^{(M_B,N_B)}_k}$, the local motivic Fourier transform
extends to this semi--local case, and we can write it in the form
$$ 
\cF_\omega\varphi = \bL^{-\sum_v N_v [k_v: k]} [ X\times_{\bA^{(M_B,N_B)}_k} \bA^{(M_B,N_B)}_k \times \bA^{(M'_B,N'_B)}_k, f\circ \pi_X + \langle u\circ \pi_1, \pi_2 \rangle_\omega ]_{\bA^{(M'_B,N'_B)}_k}, 
$$
with the Fourier inversion formula
$$ 
\cF_\omega \circ \cF_\omega = \bL^{-\sum_v [k_v: k] c_v} j^* . 
$$

\smallskip

By considering finite sets of places $B$ that include the zeros and poles of $\omega$, one can further
extend this from the semi--local to the global ad\`elic case, with Fourier inversion formula
$$ \cF_\omega \circ \cF_\omega = \bL^{-\chi(C)} j^* , $$
with $\chi(C)=2-2g$ for $g=g(C)$ the genus, see Section~1.3 of [ChamLoe15] for more details. 

\smallskip

The motivic analog of summing a Bruhat--Schwartz function over the discrete $F\subset A_F$ is defined
in the following way. 
For a divisor $D$ on $C$, let $\cL(D)$ denote the set of nontrivial rational functions $h$ with 
$div(h)+D \geq 0$ together with $0$, so that it
forms a finite dimensional $k$-vector space, $\cL(D)= H^0(C,\cO(D))$.  
Given a finite set $B\subset C(k)$ and a class $\varphi=[X,f]_{\bA_k^{(M_B,N_B)}}$,
one defines 
$$ 
\sum_{x\in F} \varphi(x) := [ \cL(D) \times_{\bA_k^{(M_B,N_B)}} X, f\circ \pi_2 ]   \in Exp\cM_k , 
$$
where $D$ is the divisor $D=-\sum_v M_v \, v$, and 
the fiber product is taken over the structure morphism $u: X \to \bA_k^{(M_B,N_B)}$ of $X$
as a $\bA_k^{(M_B,N_B)}$--variety and the morphism $\alpha: \cL(D) \to \bA_k^{(M_B,N_B)}$ 
with components determined by mapping $\cL(D)$ into $t^{M_v} \cO_v$ via the inclusions $F \hookrightarrow F_v$.

\smallskip

With this summation notation understood, the Hrushovski--Kazhdan motivic Poisson summation formula
takes the form
$$
 \sum_{x\in F} \varphi(x) = \bL^{1-g} \, \sum_{y\in F} \cF\varphi(y). 
 $$
 
 \medskip
 
 {\bf 7.14. Categorical aspects of Poisson summation.}
 We now return to our point of view based on the assembler category $\cC_k^{KExp\cS}$ and the associated
 spectrum $K(\cC_k^{KExp\cS})$ underlying the Grothendieck ring $Exp\cM_k=\pi_0 K(\cC_k^{KExp\cS})$, 
 where the Hrushovski--Kazhdan motivic Poisson summation formula takes place. Our goal here is to
 rephrase the Poisson summation at the level of objects and morphisms in the category $\cC_k^{KExp\cS}$
 discussed in Proposition~7.5 above. 
 
 \smallskip
 
 The identity in the Grothendieck ring with exponentials $Exp\cM_k$
$$ 
\sum_{x\in F} \varphi(x) = \bL^{1-g} \, \sum_{y\in F} \cF\varphi(y) 
$$
given by Poisson summation can be viewed as the following statement. Consider the classes
$$
 \sum_{x\in F} \varphi(x) = [ \cL(D) \times_{\bA_k^{(M_B,N_B)}} X, f\circ \pi_X ], 
 $$
where $D=-\sum_v M_v v$ and
$$ 
\sum_{y\in F} \cF\varphi(y) =  
$$
$$
 \bL^{-\sum_v N_v [k_v:k]}  [ \cL(D') \times_{\bA_k^{(M'_B,N'_B)}} (
X\times_{\bA^{(M_B,N_B)}_k} \bA^{(M_B,N_B)}_k \times \bA^{(M'_B,N'_B)}_k), f\circ \pi_X + \langle u\circ \pi_1, \pi_2 \rangle_\omega], 
$$
with $D'=-\sum_v M'_v v=\sum_v N_v v - \sum_v c_v v$. We view these as classes
in the Grothendieck ring $\pi_0 K(\cC^\cS)$ of the assembler $\cC^\cS$ of algebraic stacks with 
exponentials (see Proposition~7.5). The Poisson summation formula means that these classes
satisfy
$$  
 \sum_{y\in F} \cF\varphi(y)- \bL^{g-1} \sum_{x\in F} \varphi(x) \,  \in Range( \bY \cdot ), 
 $$
where $\bY\cdot : \pi_0 K(\cC^\cS) \to \pi_0 K(\cC^\cS)$, given by multiplication by $\bY=[\bA^1, id]$ is
the morphism induced by the endofunctor $\Phi: \cC^\cS \to \cC^\cS$ of Proposition~7.5 with
$$
 Exp\cM_k = \pi_0(\cC_k^{KExp\cS}) =Coker (\bY\cdot) 
 $$
with $\cC_k^{KExp\cS}=\cC^\cS/\Phi$ the cofiber of the morphism of assemblers. 

\smallskip

As shown in Section~1.3.2 of [ChamLoe15], it suffices to consider the case where the motivic ad\`elic function
$\varphi$ is a {\it simple function}, that is, the motivic analog of a product of characteristic functions of balls for
ordinary Bruhat--Schwartz functions. These are defined by assigning a finite set $B\subset C(k)$ of points, an
element $a=(a_v)_{v\in B}$ in $\prod_{v\in B} F_v$ and levels $(M_v,N_v)$ with $ord(a_v)\geq M_v$ for all $v\in B$,
with $\bA_k^{(M_B,N_B)}$ as above. A simple function is a class of the form $[Spec(k),0]_{\bA_k^{(M_B,N_B)}}$,
where the morphism $u_a: Spec(k)\to \bA_k^{(M_B,N_B)}$ assigns the $t_v$-expansion of $a_v$ for each $v\in B$.
This is regarded as the motivic counterpart of the characteristic function of the product of balls with centers $a_v$ and 
radii $N_v$ in $F_v$. More general motivic ad\`elic functions can be identified with families of simple functions
$\varphi_{a(z)}$ parameterized by $z\in Z$, where $Z$ is a $k$--variety.

\smallskip

The following result shows that the Poisson summation formula (viewed as a relation in the
ring $Exp\cM_k$) comes from a relation at the categorical level, described in terms
of a covering family in the assembler $\cC^\cS_k$.

\medskip

{\bf 7.14.1. Theorem.}
{\it For a motivic ad\`elic Bruhat--Schwartz function 
$$\varphi=[X,f]_{\bA_k^{(M_B,N_B)}} \in Exp\cM_{\bA_k^{(M_B,N_B)}},
$$  
denote by $\cQ_\varphi$, resp.  $\cF\cQ_\varphi$, the objects in the category $\cC^\cS_k$ of  Proposition~7.5 given by
$$
 \cQ_\varphi  :=(\cL(D)\times_{\bA_k^{(M_B,N_B)}} X, f\circ \pi_X) \  \in Obj(\cC^\cS_k),
 $$
$$ 
\cF\cQ_\varphi := (\cL(D') \times_{\bA_k^{(M'_B,N'_B)}}
X\times_{\bA^{(M_B,N_B)}_k} \bA^{(M_B,N_B)}_k \times \bA^{(M'_B,N'_B)}_k, f\circ \pi_X + \langle u\circ \pi_1, \pi_2 \rangle_\omega). 
$$
Then there is a covering family in $\cC^\cS_k$
$$ 
((\bA^{\deg(D)+g-1},0)\times \cQ_{\varphi,1}\hookrightarrow \cF\cQ_\varphi,
\cF\cQ_{\varphi,2} \hookrightarrow \cF\cQ_\varphi), 
$$
where $\cQ_{\varphi,1}$ is a family of simple functions with $a(z)\in \cL(div(\omega)+D)^\perp$ (the orthogonal with
respect to the Serre duality pairing), and $\cF\cQ_{\varphi,2}$ is in the range of the functor $\Phi: \cC^\cS_k \to \cC^\cS_k$
of Proposition~7.5. This covering family lifts to the level of the assembler $\cC^\cS_k$ the motivic
Poisson summation formula in $Exp\cM_k$.}

\medskip

{\it Proof.} The definition of $\sum_{x\in F} \varphi(x)$ can be equivalently described as the composition
of the pushforward to $Exp\cM_k$ of the class in $Exp\cM_{\cL(D)}$ given by the pullback of the
class $\varphi$ in $Exp\cM_{\bA_k^{(M_B,N_B)}}$, along the morphism $\cL(D)\hookrightarrow \bA_k^{(M_B,N_B)}$.

\smallskip

Using this description, we can refer to the Theorem~1.3.10 of [ChamLoe15] which shows that if $\varphi$ a simple function,
the summation $\sum_{y\in F} \cF \varphi(y)$ is zero unless $ord_v(y_v)+ord_v(D)\geq 0$. Moreover,
if  $ord_v(y_v)+ord_v(D)\geq 0$, then he computation reduces to the case of a linear function as morphism to $\bA^1$, where (as we recalled in the proof
of Theorem~7.12 above) the resulting class in $\pi_0 K(\cC^\cS_k)$ is in the range of multiplication by $\bY=[\bA^1, id]$ 
if the linear map is nonzero and is a power of $\bL$ if it is zero. 

\smallskip

This vanishing condition is satisfied when the
element $a=(a_v)$, that determines the simple function $\varphi$, belongs to the orthogonal of $\cL(div(\omega)+D)$
with respect to the Serre duality pairing, and in that case the resulting class is equal to $\bL^{-\deg(D)+\dim \cL(div(\omega)+D)}$.
Thus, for the simple function $\varphi=\varphi_a$ the summation of the motivic Fourier transforms, seen as a class
in $\pi_0 K(\cC^\cS_k)$ is given by
$$
 \sum_{y\in F} \cF \varphi(y) = \bL^{-\deg(D)+\dim \cL(div(\omega)+D)} 
 $$
when $a=(a_v)_{v\in B}\in \cL(div(\omega)+D)^\perp$, and is in the range of multiplication by $\bY=[\bA^1, id]$ otherwise.

\smallskip

By the same argument, one can show that the
left--hand--side $\sum_{x\in F} \varphi(x)$ belongs to the range of multiplication 
by $\bY=[\bA^1, id]$ unless $a=(a_v)_{v\in B}\in \cL(div(\omega)+D)^\perp$, and equal to 
$\bL^{\dim \cL(-D)}$ in that case. The Riemann--Roch formula for curves
$$
 \dim \cL(-D) = \dim  \cL(div(\omega)+D) - \deg(D) + 1-g 
 $$
then shows, that the respective classes satisfy
$$ 
\sum_{y\in F} \cF \varphi(y) - \bL^{g-1} \sum_{x\in F} \varphi(x) \in Range( \bY \cdot ). 
$$

This argument extends from the case of simple functions to the more general case by identifying
general motivic ad\`elic functions with families of simple functions over a parameterising $k$--variety $Z$.

\smallskip

In this case, one starts with a decomposition $Z=Z_1\sqcup Z_2$, where points $z\in Z_1$ have corresponding
simple functions $\varphi_{a(z)}$ with $a(z)\in \cL(div(\omega)+D)^\perp$. We denote by $Z_2=Z\smallsetminus Z_1$
the complementary set where this condition is not satisfied. The corresponding class $\sum_{y\in F} \cF \varphi(y)$
then  decomposes into a $Z_1$--part,  that can be identified as above with the corresponding
$Z_1$--part of $\bL^{\deg(D)+g-1} \sum_x \varphi(x)$, and a $Z_2$--part belonging to $Range( \bY \cdot )$.
 
\smallskip

Let now $X$ be a $\bA_k^{(M_B,N_B)}$--variety, endowed with a morphism $f: X \to \bA^1$, whose class 
$\varphi=[X,f]_{\bA_k^{(M_B,N_B)}}$ belongs to $Exp\cM_{\bA_k^{(M_B,N_B)}}$. 
Consider the objects $\cQ_\varphi$ and $\cF\cQ_\varphi$ of $\cC^\cS_k$ as in the statement, which
correspond respectively to classes $$
 [\cQ_\varphi ]=\sum_{x\in F} \varphi(x) 
 $$ and 
$$ 
[ \cF\cQ_\varphi ]=\bL^{\deg(D)} \sum_{y\in F} \cF \varphi(y) 
$$ 
in $\pi_0 K(\cC^\cS_k)$. 

\smallskip

In particular, for a simple function $\varphi=\varphi_a=[Spec(k),0]_{\bA_k^{(M_B,N_B)}}$ we write simply
$\cQ_a$ and $\cF\cQ_a$ for the corresponding objects in $\cC^\cS_k$ as above. The situation described above for the
Grothendieck classes corresponds to decompositions $[ \cQ_\varphi ]=[\cQ_{\varphi,1}]+[\cQ_{\varphi,2}]$ and
$[ \cF\cQ_\varphi ]=[\cF\cQ_{\varphi,1}]+[\cF\cQ_{\varphi,2}]$, where the first term corresponds to the
subset where the $\cQ_{a(z)}$ and $\cF\cQ_{a(z)}$ have $a(z)\in \cL(div(\omega)+D)^\perp$,  and the second is
the complementary case as above. The argument described above for the classes can then be
rephrased, as we did in Theorem~7.12 above, as the statement that, the object $\cF\cQ_\varphi$ in $\cC^\cS_k$
has a disjoint covering family
$$
 (\cF\cQ_{\varphi,1} \hookrightarrow \cF\cQ_\varphi , \cF\cQ_{\varphi,2} \hookrightarrow \cF\cQ_\varphi ) 
 $$
where one can identify $\cF\cQ_{\varphi,1}$ with $(\bA^{\deg(D)+g-1},0)\times \cQ_{\varphi,1}$ and where
$\cF\cQ_{\varphi,2}$ is in the range of the functor $\Phi: \cC^\cS_k \to \cC^\cS_k$. 
$\blacksquare$

\medskip

A generalisation of the Hrushovski--Kazhdan motivic Poisson summation formula was later obtained 
by Bilu in [Bilu18]. This more general motivic Poisson summation applies also to the case of
motivic test functions that are infinite products. This requires a different new technique based on
motivic Euler products, see Chapters~3 and 5 of [Bilu18]. 

\medskip

{\bf 7.14.2. Question.} {\it 
Is there a categorification of Bilu's motivic Euler products and of the resulting 
motivic Poisson summation formula, in an appropriate assembler category, 
as discussed here, or in the category of exponential Nori motives recalled in
Section~7.9? }

\medskip

{\bf 7.15. Motivic height zeta function} 
After reviewing the Grothendieck ring with exponential and the motivic Poisson summation,
and discussing some of their categorical properties in terms of assembler categories, we
return to the main theme of height functions. In this subsection we summarize and compare
the results of [ChamTsch02], [ChamTsch12] on the height zeta function of equivariant
compactifications of vector groups, and the motivic counterpart of [ChamLoe15] based
on a motivic version of the height zeta function and the motivic Poisson summation
that we recalled in the previous subsections. We refer the reader to [ChamLoe15] and
to the introductory chapter of [Bilu18] for a more detailed overview of this material. 

\smallskip

As already discussed in Section~2 of this paper, given a projective Fano variety (or almost Fano as
in [Pe01]) $X$ over a number field $F$,
endowed with an ample line bundle $L$, one has an associated height function $h=h_L: X(F)\to \R_+$.
For $V\subset X$ and $B>0$, the asymptotic behavior of
$$
 N_{V,h}(B)=\# \{ x\in V(F)\,|\, h(x)\leq B \} 
 $$
when $V(F)$ is infinite and for $B\to \infty$, is expected to be of the form
$$
 N_{V,h}(B) \sim C\, B^a (\log B)^{b-1}, 
 $$
for real numbers $C>0$, $a>0$ and for a half--integer $b\in \frac{1}{2} \Z$, with $b\geq 1$,
where these parameters should carry a geometric interpretation. As we already discussed 
in Section~4, in the case of the height function associated to the anticanonical line bundle,
the expected asymptotics has $a=1$, $b=rk Pic(X)$ and $C$ expressible in terms
of volumes of adelic spaces and a cohomological factor (see [Pe95], [BatTsch98]).

\smallskip

The investigation of the asymptotic behavior of the counting function $N_{V,h}(B)$ can be reformulated
in terms of questions about the associated height zeta function
$$
 \zeta_{V,h}(s)=\sum_{x\in V(F)} h(x)^{-s}. 
 $$
In particular, one can ask whether there is a dense open subvariety $V\subset X$ such that $\zeta_{V,h}(s)$
is absolutely convergent for $Re(s)>1$ with meromorphic continuation to some strip $Re(s)> 1-\delta$
and a unique pole of order $r=rk\, Pic(X)$ at $s=1$. 

\smallskip

The case of toric varieties was treated in [BatTsch98] using Poisson summation, and this approach
was generalized to varieties with an action of an algebraic group with an open dense orbit (equivariant 
compactifications of algebraic groups). The case of vector groups, that is, equivariant compactifications
of $\bG_a^n$, was studied in [ChamTsch02] and [ChamTsch12] by considering the height 
zeta function
$$
 \zeta_{G,h}(s)=\sum_{x\in G(F)} h(x)^{-s}
  $$
and using Poisson summation for $G(F)$ as a discrete subgroup of $G(A_F)$, leading to the identity
$$
 \sum_{x\in G(F)} h(x)^{-s}= \sum_{y \in G(F)} \cF(h^{-s})(y). 
 $$
The Fourier transform $\cF(h^{-s})(y)$ is then written as a product of local factors
$$
 \cF(h^{-s})(y) =\prod_v \cF(h_v^{-s})(y_v).
  $$
This can be expressed in terms of local Igusa zeta functions, which makes it then
possible to show the existence of meromorphic continuation and identify the location
and order of the main pole.

\smallskip

The formulation of Manin's problem on the asymptotic behavior of the number of rational
points of bounded height in terms of the height zeta functions extends to the case of function fields. 
For a smooth projective curve $C$ over a finite field $k=\F_q$, with function field $F=k(C)$,
and an almost Fano variety $X$ over $F$ for which the set $X(F)$ is Zariski dense in $X$,
one considers as above the height function $h$ associated to the anti-canonical line bundle
and the height zeta function $\zeta_{V,h}(t)$ with $t=q^{-s}$, and with $V\subset X$ an
open dense subset. Then the above question (absolute convergence in the disk $|t|< q^{-1}$,
meromorphic continuation in some disk $|t|< q^{-1+\delta}$, unique pole of order $r=rk Pic(X)$
at $q^{-1}$) can also be formulated in this setting.

\smallskip

In this function field setting, given a model $\cX$ of $X$ over $C$ with $u: \cX \to C$, the rational
points in $V(F)$ correspond to sections $\sigma: C \to \cX$ of $u$ with $\sigma(\eta_C)\in V(F)$
for the generic point $\eta_C$, and the height with respect to a line bundle $L$ is given by 
$h(\sigma)=q^{\deg \sigma^* L}$.  This makes it possible to rewrite the height zeta function in
the form
$$
 \zeta_{V,h}(s)=\sum_{x\in V(F)} h(x)^{-s}=\sum_{d\geq 0} N_{V,L,d}\, q^{-ds},
 $$
where
$$
 N_{V,L,d} :=\# \{ \sigma: C \to \cX \,|\, \sigma(\eta_C)\in V(F), \,\, \deg \sigma^* L = d \}. 
 $$

\smallskip

This more geometric formulation then suggests the existence of a motivic version of
the same problem, where the counting function $N_{V,L,d}$ is replaced by a class in
the Grothendieck ring of varieties. This is the approach introduced by Chambert--Loir and
Loeser in [ChamLoe15].  In this reformulation, one needs to show that the sections
$\sigma: C \to \cX$ satisfying the conditions $\sigma(\eta_C)\in V(F)$ and $\deg \sigma^* L = d$
form a moduli space $M_{V,L,d}$ that is a quasi-projective $k$--scheme, so that
it makes sense to consider the Grothendieck classes $[ M_{V,L,d} ]\in K_0(\cV_k)$ and
one can form a motivic height zeta function
$$ 
Z_{V,L}(T):= \sum_{d\in \Z} [  M_{V,L,d} ] \, T^d, 
$$
seen as a formal series in $K_0(\cV_k)[[T]][T^{-1}]$.  The analog of the height problem in
this motivic setting becomes showing that the series $Z_{V,L}(T)$, viewed as an element
in $\cM_k[[T]][T^{-1}]$ with $\cM_k$ the localization of the Grothendieck ring, has the
property that
$$ 
(1-\bL^a T^a)^b Z_{V,L}(T) = P(T)
 $$
where $P(T)$ is an element in the subring of $\cM_k[[T]][T^{-1}]$ generated by the inverses of
the polynomials $1-\bL^\alpha T^\beta$ for $\beta > \alpha \geq 0$, with value $P(\bL^{-1})$
at $T=\bL^{-1}$ given by a nontrivial effective class. 

\smallskip

This problem can be formulated more generally for the data of a smooth projective 
curve $C$ over a field $k$ and a variety $X$ with a line bundle $L$ and an open
subset $V\subset X$, with structure morphism $u: X \to C$,
for which the moduli space $M_{V,L,d}$ of sections exists as a constructible set
so that the Grothendieck class $[ M_{V,L,d} ]$ can be considered.

\smallskip

This motivic problem was solved affirmatively by Chambert--Loir and
Loeser in [ChamLoe15] in the case of equivariant compactifications of $\bG_a^n$. 
The argument is conceptually similar to the original case of 
 [ChamTsch02] and [ChamTsch12] for the ordinary height zeta function.
 The Hrushovski--Kazhdan motivic Poisson summation formula replaces
 the ordinary Poisson summation formula, by writing the motivic height
 zeta function in the form
 $$ 
 Z_{V,L}(T) = \sum_{x\in G(F)} \prod_{v\in B} (\sum_{m\in \Z} \varphi_{v,m}(x) T^m ), 
 $$
 where the $\varphi_{v,m}(x)$ are motivic Bruhat--Schwartz functions as recalled in
 the previous subsections, vanishing for $m<<0$. Applying termwise in $T$ the
 Hrushovski--Kazhdan motivic Poisson summation formula, we get
 $$ 
 Z_{V,L}(T) = \bL^{(1-g)n} \sum_{y\in G(F)} \hat Z_{V,L}(Y,y),
 $$
 where
 $$ 
  \hat Z_{V,L}(T,y)  =\prod_v  \hat Z_{V,L,v}(T,y), \ \ \ \text{ with } \ \  
 \hat Z_{V,L,v}(T,y)=  \sum_m \cF_v\varphi_{v,m}(y)\, T^m. 
 $$
 Here the Hrushovski--Kazhdan motivic Poisson summation formula
 is used in the multidimensional version 
 $$
  \sum_{x\in F^n} \varphi(x) = \bL^{(1-g)n} \sum_{y\in F^n} \cF \varphi(y),
   $$
 which is proved as in the case $n=1$ we recalled in the previous subsections 
 (see Theorem~1.3.10 of [ChamLoe15]).  The classes $[ M_{V,L,d} ]$ in
 the coefficients of $Z_{V,L}(T)$ are classes in $K_0(\cV_k)$. After
 considering them as classes in the localization $\cM_k$, and using
 the embedding $\cM_k \hookrightarrow Exp\cM_k$, the 
 Hrushovski--Kazhdan motivic Poisson summation formula applies
 and the resulting terms $\cF_v\varphi_{v,m}(y)$ give classes in the
 Grothendieck ring with exponentials $Exp\cM_k$. 
 
 \smallskip
 
 The $y=0$ term is identified as the leading term responsible for the
 main pole. This term $\hat Z_v(T,0)$ is a motivic version of
 the local Igusa zeta function and is understood in terms of
 motivic integration. This method is used in  [ChamLoe15] to
 show both the rationality and identify the leading pole of the
 motivic height zeta function.

 \medskip
 
{\bf 7.16. The lifting problem for zeta functions.}
 In the context of categorical structures underlying Grothendieck rings, one can
one can ask, which particular zeta functions of arithmetic origin  may be lifted to the categorical level. 
In the case of zeta functions associated to exponentiable motivic measures that satisfy
conditions of rationality and factorization, a lift to the categorification of the Witt ring was
constructed in [LMM19]. We review this construction here. We then outline the question
of a possible categorification for the motivic height zeta function.

\medskip

{\bf 7.16.1. Witt rings and exponentiable motivic measures.}
Let $R$ be an associative and commutative
ring. Let $End_R$ be the category of endomorphisms of projective $R$--modules 
of finite rank. The objects of this category are pairs $(E,f)$, where
$f \in End_R(E)$. With the direct sum and the tensor product defined componentwise on the objects, 
the Grothendieck group $K_0(End_R)$ also acquires a commutative ring structure. 
Let $K_0(R)$ be the ideal generated by the pairs of the form $(E, f = 0)$. Then one defines 
$$ 
W_0(R) = K_0(End_R)/K_0(R).  
$$ 

\smallskip

The ring $W_0(R)$ embeds as a dense subring of the big Witt ring $W(R)$
via the map 
$$ 
L: (E,f) \mapsto \det(1-t \, M(f))^{-1}, 
$$
where $M(f)$ is the matrix representing $f\in End_R(E)$, and
$\det(1-t \, M(f))^{-1}$ is viewed as an element in $\Lambda(R)=1+ t R[[t]]$.
The subring $W_0(R) \hookrightarrow W(R)$ 
consists of the rational Witt vectors
$$ 
W_0(R)= \left\{ \frac{1+a_1 t + \cdots + a_n t^n}{1+b_1 t+ \cdots + b_m t^m} \, | \,\, a_i,b_i\in R,\,\, n,m\geq 0 \right\}. 
$$

\smallskip

We introduce now into this picture  {\it a motivic measure}, as discussed in Section~7.8, namely a ring homomorphism from
either the Grothendieck ring of varieties or that of varieties with exponentials, 
$\mu: K_0(\cV_k)\to R$ or $\mu: KExp_k \to R$. We also consider, as in Section~7.8,
the associated zeta function $\zeta_\mu$ defined by applying the motivic measure $\mu$ to
the Kapranov motivic zeta function
$$ 
\zeta_\mu(X,t)= \sum_n \mu(Sym^n(X)) t^n 
$$
in the first case, and
$$ 
\zeta_\mu(X,f,t)= \sum_n \mu(Sym^n(X), f^{(n)}) t^n 
$$
in the case with exponentials. 

\smallskip

We can regard $\zeta_\mu (X, t)$ as defining an element in the
Witt ring $W(R)$. 

\smallskip

The addition in $K_0(\cV)$
is mapped by the zeta function to the addition in $W(R)$, which is the usual 
product of the power series,
$$
 \zeta_\mu (X \sqcup Y, t)= \zeta_\mu (X, t) \cdot \zeta_\mu (Y, t) =\zeta_\mu (X, t) +_{W(R)} \zeta_\mu (Y, t),
$$
and similarly in the case with exponentials.

\smallskip

A motivic measure $\mu: K_0(\cV) \to R$ is called {\it exponentiable} (see [Ram15], [RamTab15]),
if the zeta function $\zeta_\mu (X, t)$ determines a {\it ring homomorphism}
$$
 \zeta_\mu : K_0(\cV) \to W(R), 
 $$
that is, if it satisfies $\zeta_\mu (X \sqcup Y, t)=\zeta_\mu (X, t) +_{W(R)} \zeta_\mu (Y, t)$ as above
and it also satisfies
$$
 \zeta_\mu (X \times Y, t)= \zeta_\mu (X, t) \star_{W(R)} \zeta_\mu (Y, t) . 
$$
with the product $\star_{W(R)}$ of the Witt ring. 

\smallskip

Similarly, a motivic measure $\mu: KExp_k \to R$ is exponentiable if the associated zeta function $\zeta_\mu$
is a ring homomorphism $\zeta_\mu: KExp_k \to W(R)$. (Note that the term ``exponential" here has
two different meanings: as exponential sums in $KExp_k$ and as exponentiability of measures in the properties
of $\zeta_\mu$. The context should help the reader to avoid confusion.)

\smallskip

A motivic measure $\mu: K_0(\cV) \to R$ or $\mu: KExp_k \to R$  is {\it rational} if the zeta function
$\zeta_\mu : K_0(\cV) \to W(R)$ takes values in the subring $W_0(R)$ of the Witt ring $W(R)$. 

\smallskip

Moreover, a motivic measure $\mu: K_0(\cV) \to R$ or $\mu: KExp_k \to R$  is called {\it factorizable}
if it is rational and it admit a factorization into linear factors
$$
 \zeta_\mu(X,t)= \frac{\prod_i (1-\alpha_i t)}{\prod_j (1-\beta_j t)}. 
 $$
The latter expression given by a ratio of polynomials can also be written as a difference in the Witt ring
$$ 
\zeta_{\mu,+}(X,t)\,  -_W\,\, \zeta_{\mu,-}(X,t) ,
$$
where $\zeta_{\mu,+}(X,t)=\prod_j (1-\beta_j t)^{-1}$ and $\zeta_{\mu,-}(X,t)=\prod_i (1-\alpha_i t)^{-1}$. 

\smallskip

The motivic measure on $K_0(\cV_K)$ given by the counting function
is an exponentiable motivic measure, in the sense of [Ram15], [RamTab15].

\smallskip

One can consider the same question for the motivic measure $\mu: KExp_K \to \C$
discussed in Section~7.8, given by the exponential sum.

\smallskip

In general, a useful feature in the theory of zeta functions that makes it possible to write
them in the form of ratios of polynomials, is based on a simple identity relating
the generating series of traces of powers of an endomorphism $\sigma$ and its characteristic
polynomial,
$$ 
\exp \left( \sum_{m=1}^\infty (tr\, \sigma^m) \frac{t^m}{m} \right) = \frac{1}{\det (1-\sigma t)}. 
$$

\smallskip

Exponential sums are usually very difficult to evaluate explicitly. Cases when this is possible
typically reduce to expressing the exponential sums as traces of Frobenius on certain $\ell$--adic
sheaves. 

\medskip
{\bf 7.16.2. Categorification of Witt vectors and lifting of zeta functions.}
There are different ways to obtain a categorification and spectrification of Witt vectors.
A spectrification of the ring $W(R)$ of Witt vectors was introduced in [Hess97].
We will describe here a different categorification and spectrification of $W_0(R)$ obtained in Section~6.2 of 
[LMM19], based on its description in terms of the $K_0$ of the endomorphism category $End_R$
and of $R$, and the formalism of Segal Gamma--spaces. 

\smallskip
 
 Let $\cP_R$ denote the category
of finite projective modules over a commutative ring $R$ 
with unit, and let $End_R$ be the endomorphism
category as above. By the Segal construction, we obtain associated 
$\Gamma$--spaces $F_{\cP_R}$ and $F_{End_R}$
and spectra that we write as $F_{\cP_R}(\bS)=K(R)$ (the $K$--theory spectrum of the ring $R$), 
and $F_{End_R}(\bS)$ (the spectrum of the endomorphism category) respectively.

\smallskip
 
The spectrum $\bW(R)$ is then defined as the cofiber
$\bW(R):=F_{End_R}(\bS)/F_{\cP_R}(\bS)$ obtained from these $\Gamma$--spaces. It is
induced by the inclusion of the category $\cP_R$ of finite projective modules as the subcategory of the
endomorphism category. The spectrum $\bW(R)$ has $\pi_0 \bW(R) =W_0(R)$. 
(We refer the reader to Section~6 of [LMM19] for a more detailed discussion.)

\smallskip

It is useful to consider also a variant of the above construction with
$\cP_R^\pm$ and $End_R^\pm$ the categories of $\Z/2\Z$--graded
finite projective $R$--modules and $\Z/2\Z$--graded endomorphism category with objects
given by pairs $\{ (E_+,f_+), (E_-,f_-) \}$. Writing objects as $(E_\pm, f_\pm)$, the
morphisms are given by morphisms $\phi: E_\pm \to E'_\pm$ of $\Z/2\Z$--graded finite 
projective modules that commute with $f_\pm$. 

\smallskip

The map $\delta: K_0(End_R^\pm) \to K_0(End_R)$ given by $[E_\pm,f_\pm]\mapsto [E_+,f_+]-[E_-,f_-]$
induces a  ring homomorphism
$$ K_0(End_R^\pm) /K_0(\cP_R^\pm) \to K_0(End_R)/K_0(R) \simeq W_0(R). $$

\smallskip

As above, the categories $\cP_R^\pm$ and $End_R^\pm$ have associated $\Gamma$-spaces 
$F_{\cP_R^\pm}$ and $F_{End_R^\pm}$ and spectra. We write 
$\bW^\pm(R)=F_{End_R^\pm}(\bS)/F_{\cP_R^\pm}(\bS)$ for the cofiber
of $F_{\cP_R^\pm}(\bS) \to F_{End_R^\pm}(\bS)$.

\smallskip 

A factorizable motivic measure $\mu: K_0(\cV)\to R$ determines a
functor $\Phi_\mu:  \cC_\cV \to End_R^\pm$ where $\cC_\cV$ is the assembler category 
connected to 
the Grothendieck ring $K_0(\cV)$ of varieties (or the assembler category underlying
the Grothendieck ring with exponentials in the case of a motivic measure
$\mu: KExp_k\to R$) and $End_R^\pm$ is the $\Z/2\Z$-graded endomorphism category
described above. 

\smallskip

As shown in Sections~6.2 and 6.3 of [LMM19], this is obtained, starting with a factorization
$$ 
\zeta_\mu(X,t)= \frac{\prod_{i=1}^n (1-\alpha_i t)}{\prod_{j=1}^m (1-\beta_j t)}, 
$$
by considering $E_+^{X,\mu} =R^{\oplus m}$ and $E_-^{X,\mu}=R^{\oplus n}$ with
endomorphisms $f_\pm^{X,\mu}$ respectively given in matrix form by 
$M(f_+^{X,\mu})=diag(\beta_j)_{j=1}^m$
and $M(f_-^{X,\mu})=diag(\alpha_i)_{i=1}^n$. 
The pair $(E_\pm^{X,\mu}, f_\pm^{X,\mu})$ is an object of the endomorphism 
category $End_R^\pm$. Embeddings
$Y \hookrightarrow X$ correspond to multiplicative decompositions of the zeta function,
and the factorization of each term then determines the associated morphism in the
endomorphism category.

\smallskip

In Section~6.3 of [LMM19] it is then shown, that the functor $\Phi_\mu:  \cC_\cV \to End_R^\pm$ 
induces a map of $\Gamma$--spaces and of associated spectra $\Phi_\mu: K(\cV) \to F_{End_R^\pm}(\bS)$.
The induced maps on the homotopy groups have the property that the composition
$\delta\circ \Phi_\mu: K(\cV) \to K_0(End_R)$ followed by the quotient map 
$K_0(End_R) \to K_0(End_R)/K_0(R)=W_0(R)$, is given by the zeta function 
$\zeta_\mu: K_0(\cV)\to W_0(R)$. Thus, the functor $\Phi_\mu$ and the induced map of
spectra can be regarded as the appropriate categorification and spectrification of the zeta function.

\medskip
{\bf 7.17. The lifting problem for height zeta functions.}
In the same vein as the previous discussion of categorification and spectrification of
certain classes of zeta functions, one can ask whether a form of 
categorification and spectrification may be possible for height
zeta function as well. The case of the motivic height zeta function reviewed in
the previous subsection is especially interesting, because it lies outside of the
class of zeta functions discussed above, hence a different approach would
be needed. 

\medskip

{\bf 7.17.1. Question.}
{\it Is there a categorification and spectrification of the motivic height zeta function?}

\medskip

While we do not provide an answer to this question in the present paper, 
through the rest of this section we discuss some general aspects of 
this question and we highlight what we consider to be the most relevant and useful
properties of the motivic height zeta functions for approaching this problem.

\smallskip

Already by considering the
case of the ordinary height zeta function
$$
 \zeta_{U,h}(s)=\sum_{x\in U(F)} h(x)^{-s}
$$
rather than their motivic counterparts, one can identify
two main problems with respect to the possible construction of
a categorification and spectrification along the lines discussed in
the previous subsections.

\smallskip

The first problem is the fact that 
these zeta functions in general do not exhibit the same nice behaviour
of other arithmetic cases such as the Hasse--Weil zeta function, hence  
do not have the properties required for the categorification of [LMM19] 
reviewed in the previous subsection. The second main problem is
the fact that one would like to be able to describe them in terms of
a motivic measure on the $\pi_0 K(\cC)$ of the spectrum
$K(\cC)$ of an assembler category $\cC$.

\smallskip

For the first problem, 
the Poisson summation approach initiated in 
[BaTsch98], and applied in [ChamTsch02], [ChamTsch12] to equivariant compactifications
of vector groups, suggests that 
the object of interest for a possible lifting of the zeta function should be the
Fourier transform $\cF(h^{-s})$ and its local factors, rather than directly
working with the height zeta function itself, since the local Igusa zeta functions
that express these local factors are better behaved functions.  The reformulation
in motivic terms can help highlighting some useful properties with respect to the
second problem mentioned above.

\medskip

{\bf 7.17.2. Grothendieck classes and symmetric products.}
We first recall a convenient formalism for symmetric products introduced in [Bilu18]
to the purpose of studying motivic Euler product decompositions. Consider again
the Kapranov motivic zeta function
$$ Z_X(t)=\sum_{n\geq 0} [Sym^n X]\, t^n, $$
as recalled in Section~7.8. 

\smallskip

Given a positive integer $n$, let $\Pi(n)$ denote the set of partitions of $n$.
For $\pi\in \Pi(n)$ a partition, one denotes by  $Sym^\pi X$ the locally closed subset of $Sym^n X$
consisting of those effective zero--cycles of degree $n$ that realise the partition $\pi$.
Equivalently, for a partition $\pi$ of the form $n=\sum_i n_i$, one obtains
$Sym^\pi X$ as the quotient of the complement of the diagonals in $X^n =X^{\sum_i n_i}=\prod_i X^{n_i}$ 
by the action of the product of symmetric groups $\prod_i \Sigma_{n_i}$. 

\smallskip

Consider then the case of a $k$--variety $X$ and a family $\cX$ of $X$--varieties $\cX=(X_i)_{i\in \N}$
with structure morphisms $u_i: X_i \to X$. Given a partition $\pi=\{ n_i \}$ 
with $\sum_i n_i=n$, consider the space $\prod_i X_i^{n_i}$ and the 
complement of the diagonals $(\prod_i X^{n_i})\smallsetminus \Delta$.
As is customary in the case of configuration spaces, here $\Delta$ stands for the union of
all the diagonals (that is, the locus where two or more of the coordinates coincide), rather than just the deepest
diagonal where all coordinates coincide. We can construct the fibered product 
$$
\prod_i X_i^{n_i}\times_{\prod_i X^{n_i}} ((\prod_i X^{n_i})\smallsetminus \Delta) 
$$
with the maps $u=(u_i)$ and the inclusion. 
One defines $Sym^\pi \cX$ as the quotient of this fibered product by the group $\prod_i \Sigma_{n_i}$.
There is a natural map $Sym^\pi \cX \to Sym^\pi X$. One can define a multivariable Kapranov
zeta function
$$ Z_{\cX}(\underline{t}):=\sum_\pi [Sym^\pi \cX] \underline{t}^\pi, $$
with $\underline{t}^\pi=\prod_i t_i^{n_i}$. 

\smallskip

Given a closed embedding $Y\subset X$ with open complement $U=X\smallsetminus Y$, one
has the relation between Grothendieck classes of symmetric products
$$ 
[Sym^n X]=\sum_{r=0}^n [Sym^r Y]\cdot [Sym^{n-r} U], 
$$
which extends to the case of partitions as
$$
 [Sym^\pi X]=\sum_{\pi'\leq \pi} [Sym^{\pi'} Y]\cdot [Sym^{\pi-\pi'} U]. 
 $$
Similarly, for a family $\cX$ of $X$--varieties with $\cY\subset \cX$, $\cY=(Y_i)_{i\in \N}$, and their
complements $\cU=(U_i=X_i\smallsetminus Y_i)_{i\in \N}$, one has
$$ 
[Sym^\pi(\cX)]=\sum_{\pi'\subset \pi} [Sym^{\pi'}(\cY)]\cdot [Sym^{\pi-\pi'}\cU], 
$$
which implies the factorization
$$ 
Z_\cX(\underline{t})=Z_\cY(\underline{t})\cdot Z_\cU(\underline{t}). 
$$

The reason for introducing these multivariable Kapranov
zeta functions (see [Bilu18]) is summarized in the next subsection.

\medskip
{\bf 7.17.3. Summary of the geometric setting.}
In order to better identify the underlying geometry, relevant to the possible
construction of a related assembler category, we recall the following general
setting from [Bilu18].  As above we have a smooth projective curve $C$ over $k$
with function field $F=k(C)$. We also consider a smooth equivariant compactification
$X$ of $G=\bG_a^n$ as a smooth projective $F$--scheme, with $X\smallsetminus G=D_G$ 
a strict normal crossings divisor of components $(D_\alpha)_{\alpha \in \cA}$, as well as a 
partial compactification given by
a $G$--invariant quasi--projective scheme $U\subset X$, with strict normal crossings divisor
$X\smallsetminus U=D \subset D_G$, where $-K_X(D)$ is the log--anticanonical class of $U$.

\smallskip

This geometry is also described in terms of a good model given by a $k$--scheme $\cX$ with $u: \cX\to C$ with
an open $\cU \subset \cX$, with $\cX\smallsetminus \cU$ a strict normal crossings divisor, 
together with a line bundle $\cL$ over $\cX$ that restricts to $-K_X(D)$
on the generic fiber $X$. For $v\in C(k)$, if $\cB_v$ denotes the set of irreducible components
of $u^{-1}(v)$ and $\cB=\sqcup_{v\in C(k)} \cB_v$, the line bundle $\cL$ can be written as
$\sum_\alpha \lambda_\alpha \cL_\alpha$, for integers $\lambda_\alpha$ and for line bundles 
$\cL_\alpha =\cD_\alpha+\sum_{\beta\in \cB} c_{\alpha,\beta} E_\beta$ with $E_\beta$ the 
component corresponding to $\beta\in \cB$.

\smallskip

Given such a setting, one can consider the moduli spaces $M_{\underline{n},U}$, for
$\underline{n}=(n_\alpha)_{\alpha\in \cA}$ of sections 
$\sigma: C \to \cX$ with $\sigma(\eta_C)\subset G$, satisfying the integrality condition $\sigma(C_0)\subset \cU$ for
$C_0\subset C$ a subset of places, and satisfying $\deg(\sigma^* \cL_\alpha)=n_\alpha$. 
The moduli spaces $M_{n,U}$ recalled earlier in this section are decomposed as a disjoint union
of the $M_{\underline{n},U}$ for $n=\sum_\alpha \lambda_\alpha n_\alpha$. As shown in
Section~6 of [Bilu18], one can further decompose the $M_{\underline{n},U}$ into level sets of
intersection indices of the sections with the components of the boundary divisor. Namely, 
given an element $g\in G(F)$ and the corresponding section $\sigma_g : C \to \cX$ that extends it, one has
$$
 \deg (\sigma_g^* (\cD_\alpha))=\sum_{v\in C(k)} (g,\cD_\alpha)_v, 
 $$
where $(g,\cD_\alpha)_v$ denote the intersections indices. Similarly, one has intersection indices
$(g,E_\beta)_v$, which have value $1$ for a single $E_\beta$ and zero otherwise ([Bilu18], Sec.6.2.2 and
[ChamLoe15], Sec.3.3). One can define the level sets of these intersection indices as
$$
 G(\underline{m}_v,\beta)_v=\{ g\in G(F_v)\,|\, (g,E_\beta)_v=1 \ \text{ and } \  (g,\cD_\alpha)_v=m_{v,\alpha} \}, 
 $$
where $\underline{m}_v=(m_{v,\alpha})_{\alpha\in \cA}$ and $\sum_{v\in C(k)} \underline{m}_v=\underline{n}=(n_\alpha)_{\alpha\in \cA}$. When taking also into account the property $\sigma(C_0)\subset \cU$, these points
$\underline{m}_v=(m_{v,\alpha})_{\alpha\in \cA}$ are parametrised by 
$$ 
Sym^{\underline{n}'}(C\smallsetminus C_0)\times Sym^{\underline{n}''}(C) .
$$ 
Here $\underline{n}'$
consists of those $n_\alpha$, for which $D_\alpha$ is a component of $X\smallsetminus U$, and $\underline{n}''$
consists of the $n_\alpha$ of the remaining components of $X\smallsetminus G$. The sets
$H(\underline{m},\beta)_v$ are defined as $G(\underline{m}_v,\beta)_v$ when $(\underline{m},\beta)$
satisfy the integrality condition and $\emptyset$ otherwise. It is shown in Section~6.2.6 of [Bilu18] that
the Grothendieck classes of the subsets $M_{\underline{n},\beta}$ of the moduli spaces of sections $M_{n,U}$
with assigned intersection indices decompose as
$$ [ M_{\underline{n},\beta} ]=\sum_{\underline{m}\in Sym^{\underline{n}^\beta}(C)} [ H(\underline{m},\beta)\cap G(F) ], $$
where $\underline{n}^\beta=(n_\alpha^\beta)_{\alpha\in \cA}$, where 
$\deg(\sigma^* \cD_\alpha)=n_\alpha^\beta := n_\alpha -\sum_v c_{\alpha,\beta_v}$.
Thus, for the multivariable version of the motivic height zeta function one obtains

$$
\eqalign{ 
Z_{U,L}(\underline{T}) & = \sum_{\underline{n}\in \N^\cA} [M_{\underline{n}.U}] \underline{T}^{\underline{n}} \cr
& = \sum_{\underline{n},\beta} \sum_{\underline{m}\in Sym^{\underline{n}^\beta}(C)} [ H(\underline{m},\beta)\cap G(F) ] \underline{T}^{\underline{n}}. }
$$
The Poisson summation formula, in the more general form proved in [Bilu18],  can then be applied to
$$
 \sum_{\underline{m}\in Sym^{\underline{n}^\beta}(C)} [ H(\underline{m},\beta)\cap G(F) ] =
\sum_{\underline{m}\in Sym^{\underline{n}^\beta}(C)} \sum_{x\in F^n} 1_{H(\underline{m},\beta)}(x), 
$$
and analyzed in terms of the properties of the motivic Fourier transforms $\cF (1_{H(\underline{m},\beta)})$,
where $1_{H(\underline{m},\beta)}$ denotes the family of motivic Bruhat--Schwartz functions parameterized
by $Sym^{\underline{n}^\beta}(C)$. This family can also be described (Section~6.2.5 of [Bilu18]) in terms
of symmetric products $Sym^{\underline{n}}(\cH_\beta)$, where
$\cH_\beta=(H_{\underline{m},\beta})_{\underline{m}\in \N^\cA}$ is the family over $C$ with $H(\underline{m},\beta)$
over $v\in C(k)$.
Thus, we can consider the associated multivariable zeta functions as in Section~7.12.2 above,
$$ Z_{\cH_\beta}(\underline{T})=\sum_{\underline{n}\in \N^\cA} [Sym^{\underline{n}}(\cH_\beta)] \underline{T}^{\underline{n}}. $$

\medskip
{\bf 7.17.4. Poisson summation and rationality.}
Another step towards achieving a suitable setting for a lifting of the motivic height zeta function is
provided by the following result of [ChamLoe15] and [Bilu18]. Using motivic Poisson summation,
the motivic height zeta function is rewritten in the form
$$ \bL^{(1-g)N} \sum_{\xi\in F^n} Z(T,\xi), $$
where the $Z(T,\xi)=\prod_v Z_v(T,\xi)$ are sums of the motivic Fourier transforms of the form
$\cF (1_{H(\underline{m},\beta)})$
mentioned above. It is then shown (see Proposition~5.3.4 of [ChamLoe15]) that for each $\xi$ the multivariable
$Z(\underline{T},\xi)$ are rational functions with denominators given by products of terms of the form $1-\bL^n \underline{T}^{\underline{m}}$. with $n\in \N$ and $\underline{m}\in \N^\cA$.

\medskip
{\bf 7.17.5. Scissor congruences and assembler.} 
We finally discuss briefly the problem of finding a suitable scissor congruence relations
group (or ring) $\pi_0 K(\cC)$ and an
underlying assembler $\cC$ and spectrum $K(\cC)$, that makes it possible to realize the zeta
functions one wants to lift as a morphism of additive groups (or possibly of rings) from
$\pi_0 K(\cC)$ to $Exp\cM_k[[T]][T^{-1}]$ endowed with the product of Laurent series as 
addition (and a Witt--ring multiplication). The compatibility with addition means the requirement
that the zeta function splits multiplicatively under scissor congruence relations in $\pi_0 K(\cC)$.

\smallskip

As mentioned in Sections~7.17.2 and 7.17.3 above, using the zeta functions
$Z_{\cH_\beta}(\underline{T})$ briefly described at the end of Section~7.17.3, has
the advantage that these do satisfy a scissor congruence relation as recalled in 
Section~7.17.2. 

\smallskip

One expects to obtain a suitable assembler $\cC=\cC_{X,G}$, for $X$ over $F=k(C)$,
an equivariant compactification of $G=G_F$ with $X\smallsetminus G=D_G$ a strict normal 
crossings divisor. This may be built with objects given by choices of data $(C_0,D)$ with 
$C_0\subset C$ and $D=X\smallsetminus U$ with $U$ a partial compactification,
with morphisms given by inclusions and with a Grothendieck topology generated by
covering families of the form
$$ (C_0',D') \hookrightarrow (C_0,D) \hookleftarrow (C_0\smallsetminus C_0',D\smallsetminus D') . $$

\smallskip

In fact, this idea for a construction of an assembler $\cC=\cC_{X,G}$ may be further refined by taking into
consideration the structure of the Clemens complex considered in [ChamLoe15] and [Bilu18]
that is reflected in the structure of the $Z(T,\xi)$ as rational functions. This is the simplicial complex
with vertices the components of a normal crossings divisor $D$ and a $n$--simplex for each irreducible
component of a nonempty $(n+1)$-fold intersection of components of $D$. The combinatorics of
this complex influences the structure of the local factors, where for instance one has (Proposition~4.3.2
of [ChamLoe15])
$$ Z_v(\underline{T},0) =\sum_A P_{v,A}(\underline{T}) \prod_{\alpha\in A} \frac{1}{1-\bL^{\rho_\alpha-1}T_\alpha} $$
with $A$ ranging over the set of maximal faces of the Clemens complex.

\bigskip

{\bf Acknowledgements.}  M.  Marcolli  acknowledges support
from NSF grants DMS--1707882 and DMS--2104330 and from  
NSERC grants RGPIN--2018--04937 and
RGPAS--2018--522593.

\smallskip

Yu. Manin acknowledges the excellent scientific environment of
the Max Planck Institute for Mathematics in Bonn and permanent
support of its administration and of the Max Planck Society.

\smallskip

We thank the three anonymous referees for a very careful reading of the
paper and for providing many detailed comments and suggestions
that greatly improved the paper. 

\bigskip

\centerline{\bf References}

\medskip
[BaMa90] V.  Batyrev,  Yu.  Manin.  {\it Sur le nombre des points rationnels de hauteur born\'ee
des vari\'et\'es alg\'ebriques.} Math. Ann. 286 (1990), pp. 27--43. 

\smallskip

[BatTsch98] V.~Batyrev, Yu.~Tschinkel, {\it Manin's conjecture for toric varieties}, 
Journ. of Alg. Geom., 7 (1998) no. 1, 15--53.

\smallskip

[BeDh07] K.~Behrend, A.~Dhillon, {\it On the motivic class of the stack of bundles}, Adv. Math.
212 (2007), no. 2, 617--644. 

\smallskip

[Bilu18] M.~Bilu, {\it Motivic Euler products and motivic height zeta function}, \newline arXiv:1802.06836.

\smallskip

[BlBrDeGa20] V.  Blomer,  J.  Br\"udern,  U. Derenthal, G.  Gagliardi.
{\it The Manin--Peyre conjecture for smooth spherical Fano varieties 
of semisimple rank one.} arXiv:
2004.09357, 76 pp. 

\smallskip

[Bour09] D.~Bourqui, {\it 
Produit eul\'erien motivique et courbes rationnelles sur les vari\'et\'es toriques}. 
Compos. Math. 145 (2009), no. 6, 1360--1400.

\smallskip

[Carl05] G.~Carlsson, {\it Deloopings in algebraic $K$-theory}. Handbook of $K$-theory. 
Vol.~1, 2, pp.~3--37, Springer, 2005.

\smallskip

[Cham10] A.~Chambert-Loir, 
{\it Lectures on height zeta functions: at the confluence of algebraic geometry, algebraic number theory, and analysis}. Algebraic and analytic aspects of zeta functions and L-functions, 17--49,
MSJ Mem., 21, Math. Soc. Japan, Tokyo, 2010.

\smallskip

[ChamLoe15] A.~Chambert-Loir, F.~Loeser, {\it Motivic height zeta function}, 
Amer. J. Math. 138 (2016), no. 1, 1--59. arXiv:1302.2077v4.

\smallskip

[ChamTsch02] A.~Chambert-Loir, Yu.~Tschinkel, {\it On the distribution of points of bounded height on equivariant compactifications of vector groups}, Invent. Math. 148 (2002), 421--452.

\smallskip

[ChamTsch12] A.~Chambert-Loir, Yu.~Tschinkel, {\it Integral points of bounded height on partial equivariant compactifications of vector groups}, Duke Math. Journal 161 (2012) N.15, 2799--2836. 

\smallskip

[CluLoe10] R.~Cluckers, F.~Loeser, {\it 
Constructible exponential functions, motivic Fourier transform and transfer principle},
Ann. of Math. (2) 171 (2010), no. 2, 1011--1065.

\smallskip

[CThSk19]  J.-L. Colliot-Th\'elene, A. Skorobogatov. {\it The Brauer--Grothendieck 
group.} imperial.ac.uk 2019, 360 pp. 

\smallskip

[CoPaSk16] J.-L.~Colliot-Th\'elene, A.~P\'el, A.~Skorobogatov, {\it
Pathologies of the Brauer-Manin obstruction}. 
Math. Z. 282 (2016), no. 3-4, 799--817.

\smallskip

[CoCons16] A.~Connes,  C. Consani. {\it Absolute algebra and Segal's $\Gamma$--rings
au dessous de  $\overline{\roman{Spec}\,(\bold{Z})}$.} J. Number Theory, 162 (2016), pp. 518--551.  

\smallskip

[CorSch20] D.  Corwin,  T. Schlank.  {\it Brauer and etale homotopy obstructions to rational
points on open covers.} arXiv:2006.11699v3 [math.NT], 48 pp.  

\smallskip

[DeMaRa07] P.~Deligne, B.~Malgrange, J-P.~Ramis, {\it Singularit\'es irr\'eguli\`eres}, Documents Math\'ematiques (Paris), vol. 5,
Soci\'et\'e Math\'ematique de France, Paris, 2007, Correspondance et documents.

\smallskip

[DePi19] U. Derenthal,  M.  Pieropan. {\it The split torsor method for Manin conjecture.}
arXiv:1907.09431, 36 pp. 

\smallskip

[Ek09] T.~Ekedahl, {\it The Grothendieck group of algebraic stacks}, arXiv:0903.3143v2. 

\smallskip

[EKMM97] A.D.~Elmendorf, I.~Kriz, M.A.~Mandell, J.P.~May, {\it Rings, modules, and 
algebras in stable homotopy theory}. With an
appendix by M. Cole, Mathematical Surveys and Monographs, Vol.47, American Mathematical Society, 
1997, xii+249 pp.

\smallskip

[FMa20] F.  Manin. {\it Rational homotopy type and computability.}
arXiv:2007.10632, 18 pp. 

\smallskip

[FreJo20] J. Fres\'an,  P.  Jossen.  {\it Exponential motives.} 
Manuscript available at http://javier.fresan.perso.math.cnrs.fr/expmot.pdf 

\smallskip

[FrMaTsch89] J.  Franke, Yu.  Manin,  Yu.  Tschinkel.  {\it Rational points
of bounded height on Fano varieties.} Invent. Math., 95, no. 2 (1989),
pp. 421--435.  

\smallskip

[GeMa03] S. Gelfand, Yu. Manin. {\it Methods of homological algebra.}
Second Ed., Springer Monograph in Math., (2003), xvii + 372 pp.  

\smallskip

[Hess97] L.~Hesselholt, {\it Witt vectors of non-commutative rings and topological cyclic homology}, Acta
Math., Vol. 178 (1997) N.1, 109--141.

\smallskip

[HoShiSmi00] M.~Hovey, B.~Shipley, J.~Smith, {\it Symmetric spectra}, J. Amer. Math. Soc. 13 (2000), 149--208.

\smallskip

[HruKaz09] E. Hrushovski, D. Kazhdan. {\it Motivic Poisson summation.} arXiv:0902.0845

\smallskip

[HuMu-St17] A.~Huber, S.~M\"uller-Stach, {\it Periods and Nori motives}, Ergebnisse der Mathematik und ihrer Grenzgebiete, vol.~65, Springer, 2017.  

\smallskip

[Ig78] J.-I. Igusa. {\it Lectures on forms of higher degree.} Tata Inst.  of Fund. Research, 
Bombay., Narosa Publ. House New Delhi (1978), 169 pp.  

\smallskip

[Kap00] M.~Kapranov, {\it The elliptic curve in the $S$-duality theory and Eisenstein series for Kac-Moody groups}, arXiv:math/0001005.

\smallskip

[KoZa01] M.~Kontsevich, D.~Zagier, {\it Periods}, in ``Mathematics unlimited -- 2001 and beyond", 
Springer, 2001, pp.~771--808.

\smallskip
[La76] R.  Langlands. {\it On the functional equations satisfied by Eisenstein series.}
Lecture Notes in Math. vol. 544 (1976), Springer Verlag, Berlin--Heidelberg--New York,
337 pp.   

\smallskip

[LeSeTa18] B.  Lehman, A.  Sengupta,  S.  Tanimoto. {\it Geometric consistency
of Manin's conjecture.} arXiv:1805.10580 

\smallskip

[LiXu15] Qing Liu, Fei Xu. {\it Very strong approximation for certain algebraic
varieties.}  Math. Ann. 364 (3--4) (2015), pp. 701--731.  

\smallskip

[LMM19]  J.~Lieber, Yu.I.~Manin, M.~Marcolli, {\it Bost--Connes systems and $\F_1$-structures in Grothendieck
rings, spectra, and Nori motives}, to appear in ``Facets of Algebraic Geometry: A Volume in 
Honour of William Fulton’s 80th Birthday", Cambridge University Press, 2020.  arXiv:1901.00020

\smallskip

[Ly99]  M.  Lydakis. {\it Smash products and $\Gamma$--spaces.} Math. Proc.
Cambridge Phil. Soc., 126 (1999), pp. 311--328.     

\smallskip

[MaMar18]  Yu. Manin, M. Marcolli. {\it Homotopy types and geometries below $\roman{Spec}\, \bold{Z}$.} In: Dynamics: Topology and Numbers. Conference to the Memory of Sergiy Kolyada. Contemp. Math., AMS 744 (2020). arXiv:math.CT/1806.10801.

\smallskip

[Pe95] E.~Peyre, {\it Hauteurs et mesures de Tamagawa sur les vari\'et\'es de Fano}, 
Duke Math. J. 79 (1995), No 1, 101--218. 

\smallskip

[Pe01] E. Peyre.  {\it Points de hauteur born\'ee et g\'eom\'etrie des vari\'et\'es (d'apr\`es 
Y. Manin et al.)}, S\'eminaire Bourbaki, Vol. 2000/2001, Ast\'erisque No. 282 (2002), Exp. No. 891, ix, 323--344.

\smallskip

[Pe03] E. Peyre.  {\it Points de hauteur born\'ee, topologie ad\'elique
et mesures de Tamagawa.} J. Th\'eor. Nombres Bordeaux 15 (2003), pp. 319--349.

\smallskip

[Pe17] E. Peyre. {\it Libert\'e et accumulation.} Documenta Math., 22 (2017),
pp. 1616--1660.

\smallskip

[Pe18] E. Peyre. {\it Beyond heights: slopes and distribution of rational
points.} arXiv:1806.11437, 53 pp.  

\smallskip

[Pe19] E. Peyre. {\it Les points rationnels.} SMF Gazette, 159  

\smallskip

[Po17] B.  Poonen.  {\it Rational points on varieties.} Graduate studies in Math.,
186 (20017), AMS, Providence, RI.  

\smallskip

[Qu15] G.  Quick. {\it Existence of rational points as a homotopy limit problem.}
J. Pure Appl. Algebra 219:8 (2015), pp. 3466--3481. 

\smallskip

[Qui73]  D.~Quillen, {\it Higher algebraic K-theory. I}. Algebraic K-theory, I: Higher
K-theories (Proc. Conf., Battelle Memorial Inst., Seattle, Wash., 1972), pp.~85--147. 
Lecture Notes in Math., Vol. 341, Springer, 1973.

\smallskip

[Ram15] N.~Ramachandran, {\it Zeta functions, Grothendieck groups, and the Witt ring}, 
Bull. Sci. Math. Soc. Math. Fr. 139 (2015) N.6, 599--627. arXiv:1407.1813

\smallskip

[RamTab15] N.~Ramachandran, G.~Tabuada, {\it 
Exponentiable motivic measures},  J. Ramanujan Math. Soc. 30 (2015), no. 4, 349--360. 
arXiv:1412.1795

\smallskip

[Sa20] W.  Sawin.  {\it Freeness alone is insufficient for Manin--Peyre.}  \newline
arXiv:math.NT/2001.06078. 7 pp.  

\smallskip

[Scha79] S.~Schanuel, {\it Heights in number-fields}. Bull. Soc. Math. France 107 (1979), 433--449.

\smallskip

[Sc01] S.  Schwede.  {\it $\bold{S}$--modules and symmetric spectra.}
Math. Ann. 319 (2001), pp. 517--532. 

\smallskip

[Sc12] S.~Schwede {\it Symmetric spectra}, book manuscript, 
{\tt http://www.math.uni-bonn.de/people/schwede/SymSpec-v3.pdf}

\smallskip

[Se74] G. Segal. {\it Categories and cohomology theories.} Topology,
vol. 13 (1974), pp. 293--312 .   

\smallskip

[Sko99] A. Skorobogatov {\it Beyond the Manin obstruction.}  Invent.  Math. 135:2 (1999),
pp. 399--424. 

\smallskip
[Sko09] A.  Skorobogatov.  {\it Descent obstruction is equivalent to \'etale Brauer--Manin
obstruction.} Math. Ann. 344:3 (2009), pp.  501--510. 

\smallskip

[SpWh53] E. H. Spanier, J. H. C.  Whitehead.  {\it A first approximation to
homotopy theory.} Proc. Nat. Ac.  Sci. USA, 39 (1953), pp. 656--660.   

\smallskip

[Ta19] Sh. Tanimoto. {\it On upper bounds of Manin type.} arXiv:1812.03423v2, 29 pp. 

\smallskip

[Tho95] R.W.~Thomason, {\it Symmetric monoidal categories model all connective spectra}, 
Theory Appl. Categ. 1 (5) (1995) 78--118.

\smallskip

[ToVa09] Bertrand T\"oen and Michel Vaqui\'e, {\it Au-dessous de $\roman{Spec} \bold{Z}$}, 
J. K-Theory 3 (2009), no. 3, 437--500.

\smallskip

[Wy17] D.~Wyss, {\it Motivic classes of Nakajima quiver varieties}, 
International Mathematics Research Notices (2017) N.~22, 6961--6976. 

\smallskip

[Za17a] I. Zakharevich. {\it The $K$--theory of assemblers.} Adv.  Math. 304 (2017),
pp. 1176--1218 .    

\smallskip

[Za17b] I. Zakharevich. {\it On $K_1$ of an assembler.} J. Pure Appl. Algebra
221, no.7 (2017),  pp. 1867--1898 .    

\smallskip

[Za17c] I. Zakharevich. {\it
The annihilator of the Lefschetz motive}, Duke Math. J. 166 (2017), no. 11, 1989--2022

\bigskip

{\bf Yuri I. Manin, Max-Planck-Institut f\"ur Mathematik, Vivatsgasse 7, 53111 Bonn, Germany}

\medskip

{\bf Matilde Marcolli, Math. Department, Mail Code 253-37, Caltech, 1200 E.California Blvd.,
 Pasadena, CA 91125, USA}

\enddocument